\documentclass[12pt]{amsart}

\usepackage{amsfonts, amstext, amsmath, amsthm, amscd, amssymb, enumitem, url}
\usepackage[table]{xcolor}
\usepackage{tikz-cd}
\usepackage{pdfpages}
\usepackage{pdflscape}
\usepackage{mathabx}
\usepackage{mathtools}

\usepackage{geometry}
\usepackage[parfill]{parskip}
\usepackage{microtype}

\usepackage{tabularx}

\newcolumntype{x}{>{\centering\arraybackslash}X}
\usepackage{makecell}
\usepackage{longtable}
\usepackage{hhline}
\usepackage{multirow}
\usepackage{subfig}
\usepackage{adjustbox}
\usepackage{verbatim}

\usetikzlibrary{arrows.meta}

\usepackage[colorlinks,citecolor=cyan,linkcolor=magenta]{hyperref}
\usepackage{zref-clever} 
\zcsetup{cap = true, nameinlink = true}
\newcommand{\Cref}[1]{\zcref{#1}}

\usepackage{thmtools} 

\newtheorem*{thm:flowreducetogeneralizedht}{\Cref{thm:flowreducetogeneralizedht}}

\newtheorem{thm}{Theorem}[section]
\zcRefTypeSetup{thm}{Name-sg = Theorem}

\usepackage{xparse} 

\NewCommandCopy{\newtheoremcopy}{\newtheorem}
\RenewDocumentCommand{\newtheorem}{m O{thm} m}{ 
\newtheoremcopy{#1}[#2]{#3}
\AddToHook{env/#1/begin}{\zcsetup{countertype={thm=#1}}}
\zcRefTypeSetup{#1}{Name-sg = #3}}

\newtheorem{lemma}[thm]{Lemma}
\newtheorem{prop}[thm]{Proposition}
\newtheorem{cor}[thm]{Corollary}

\newtheorem{conj}[thm]{Conjecture}

\RenewDocumentCommand{\newtheorem}{m O{thm} m}{ 
\newtheoremcopy{#1}[#2]{#3}
\AddToHook{env/#1/begin}{\zcsetup{countertype={thm=#1}}}
\zcRefTypeSetup{#1}{Name-sg = #3}
\AtBeginEnvironment{#1}{\pushQED{\hfill $\lozenge$}}
\AtEndEnvironment{#1}{\popQED}}

\theoremstyle{definition}
\newtheorem{defn}[thm]{Definition}
\newtheorem{rmk}[thm]{Remark}
\newtheorem{constr}[thm]{Construction}

\newtheorem{quest}[thm]{Question}

\numberwithin{equation}{section}

\newcounter{notes}

\usepackage{todonotes}

\renewcommand{\epsilon}{\varepsilon}
\newcommand{\vphi}{\varphi}

\usepackage{stmaryrd}

\newcommand{\std}{\mathrm{std}}

\DeclareMathOperator{\intr}{\mathrm{int}}

\DeclareMathOperator{\wind}{\mathrm{wind}}
\DeclareMathOperator{\Mod}{\mathrm{Mod}}
\DeclareMathOperator{\PMod}{\mathrm{PMod}}
\DeclareMathOperator{\genus}{\mathrm{genus}}

\let\div\relax
\DeclareMathOperator{\div}{\mathrm{div}}

\allowdisplaybreaks

\begin{document}
\title[Only one Anosov flow]{There is only one transitive Anosov flow \\
{\tiny (in dimension 3, with orientable foliations, up to bicontact padding and Goodman surgeries)}} %

\author[T.~Barthelm\'e]{Thomas Barthelm\'e}
\address{Queen's University, Kingston, Ontario}
\email{thomas.barthelme@queensu.ca}
\urladdr{https://sites.google.com/site/thomasbarthelme}

\author[F.~Salmoiraghi]{Federico Salmoiraghi}
\email{federico.salmoiraghi@gmail.com}

\author[C.~C.~Tsang]{Chi Cheuk Tsang}
\address{Tongji University, Shanghai 200092, China}
\email{chicheuk@tongji.edu.cn}
\urladdr{https://chicheuktsang.github.io}

\begin{abstract}
We introduce a new surgical operation on Anosov flows, which we call bicontact padding surgery.
The operation is defined using the machinery of strongly adapted bicontact structures, and involves inserting a `padding' layer where one of the contact structures rotates many times around the Reeb flow of the other.
We show that every skew Anosov flow with orientable stable and unstable foliations can be transformed into a finite cover of a geodesic flow using bicontact padding surgery and Goodman surgery.
Together with previous results, this implies that every transitive Anosov flow with orientable stable and unstable foliations can be transformed into a suspension flow using bicontact padding surgery and Goodman surgery.
This gives some progress towards the Christy--Fried--Ghys conjecture, which predicts that the same statement is true using only Goodman surgery.
\end{abstract}

\maketitle

\section{Introduction} \label{sec:intro} 

An \textbf{Anosov flow} on a closed 3-manifold is a flow that contracts along a stable direction and expands along an unstable direction. These flows were introduced by Anosov in the 1960s as objects of dynamical interest. Through work of Thurston, Ghys, Barbot, Fenley, and many others, it was later realized that Anosov flows also interact with 3-manifold topology in deep and interesting ways.

In recent years, much progress has been made on some fundamental questions about Anosov flows: 
We know that there exist 3-manifolds that do not admit Anosov flows (see, for instance \cite{PT72, Ghy84, Yu23}),
and there exist 3-manifolds that admit arbitrarily many non-orbit equivalent Anosov flows \cite{BBY17, BM22, CP25, BY24}. 
Following earlier work of Barthelmé--Bowden--Mann \cite{BBowM24}, Gabai--Li \cite{GL26} and Bowden--Colin \cite{BC26} have recently shown that every atoroidal 3-manifold admits at most finitely many Anosov flows up to orbit equivalence.
Moreover, we have a complete classification of Anosov flows on some simple 3-manifolds \cite{BF13, BP26, Yu23}.

In this paper, we address a different perspective on the classification problem:
In \cite{Goo83}, Goodman showed that one can construct new Anosov flows from old ones by performing Dehn surgery at a closed orbit and modifying the flow suitably. This operation is now commonly referred to as \textbf{Goodman surgery}, or \textbf{Goodman--Fried surgery}.
By performing this surgery on \textbf{Anosov suspension flows}, i.e., suspension flows of Anosov diffeomorphisms on the torus, she constructed the first examples of Anosov flows on hyperbolic 3-manifolds.

A question of Fried \cite{Fri83} effectively asks whether, barring some obvious necessary conditions, every Anosov flow arises this way. Ghys also asked a version of this question (orally) in the 1990s and the question also appears in the second version of Kirby's problem list \cite[Problem 3.54]{Kir97}, with an attribution to Christy.
The modern formulation of the question is the following conjecture.

\begin{conj}[Christy, Fried, Ghys] \label{conj:friedchristyghys}
Every transitive Anosov flow with orientable foliations can be transformed into an Anosov suspension flow via Goodman surgery.
\end{conj}

Here, we say that an Anosov flow has \textbf{orientable foliations} if its stable and unstable line fields are orientable.%

In spirit, and in particular as promoted by Ghys, this conjecture is a classification problem in the following sense: While the full classification of transitive Anosov flows up to orbit equivalence is bound to be extremely complex, due to the huge variety of known examples, the Christy--Fried--Ghys conjecture says that, if one further allows surgeries, then the classification problem becomes almost trivial.

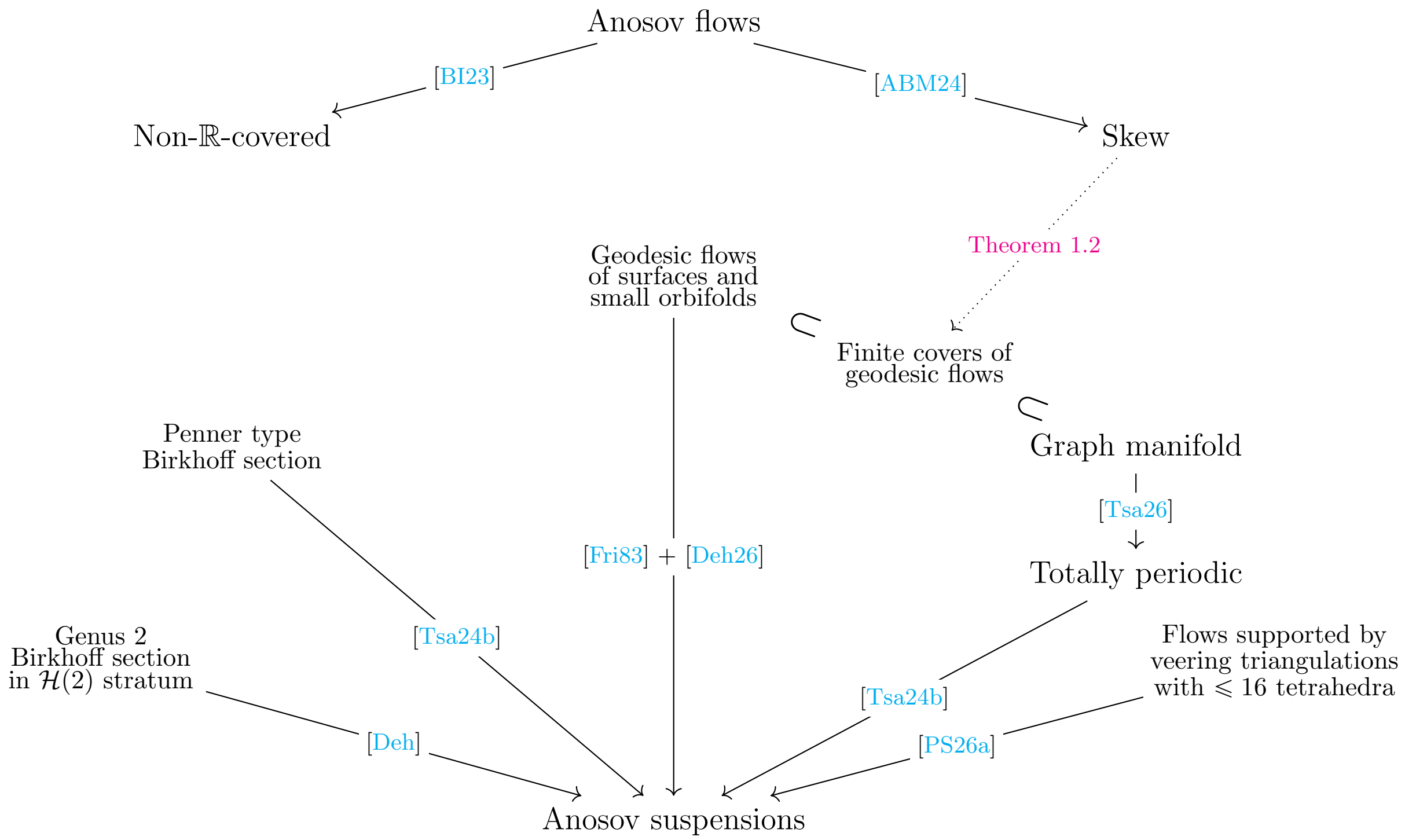
\begin{figure}

\centering

\adjustbox{scale=0.85,center}{
\begin{tikzcd}[column sep = tiny, row sep = tiny, labels = {font=\scriptsize}]
 &  &  & \text{Anosov flows} \arrow[lldd, "{\text{\cite{BI23}}}"{anchor=center,fill=white}] \arrow[rrdd, "{\text{\cite{ABM24}}}"{anchor=center,fill=white}] &  &  &  \\
 &  &  & \phantom{extra row} &  &  &  \\
 & \text{Non-$\mathbb{R}$-covered} &  &  &  & \text{Skew} \arrow[lddd, dotted, "{\text{\Cref{thm:flowreducetogeneralizedht}}}"{anchor=center,fill=white}] &  \\
 &  & \phantom{extraspa} & \phantom{extra row} &  &  &  \\
 &  &  & \text{\large{$\substack{\text{Geodesic flows} \\ \text{of surfaces and} \\ \text{small orbifolds}}$}} \arrow[rd, phantom, sloped, "\subset"{font=\Large}] \arrow[ddddddddd, "{\text{\cite{Fri83} + \cite{Deh26}}}"{anchor=center,fill=white}] &  &  &  \\
 &  &  &  & 
 \phantom{extra}\mathclap{\text{\large{$\substack{\text{Finite covers of} \\ \text{geodesic flows}}$}}}\phantom{spac} 
 \arrow[rd, phantom, sloped, "\subset"{font=\Large}] &  &  \\
 & \text{\large{$\substack{\text{Penner type} \\ \text{Birkhoff section}}$}} \arrow[rrddddddd, "{\text{\cite{Tsa24b}}}"{anchor=center,fill=white}] &  &  &  & \text{Graph manifold} \arrow[dd, "{\text{\cite{Tsa26}}}"{anchor=center,fill=white}] &  \\
 &  &  & \phantom{extra row} &  &  &  \\
 &  &  &  &  & \text{Totally periodic} \arrow[llddddd, "{\text{\cite{Tsa24b}}}"{anchor=center,fill=white}] &  \\
 &  &  & \phantom{extra row} &  &  &  \\
|[fill=white, overlay, inner sep=40pt]| \mathclap{\text{\large{$\substack{\text{Genus 2} \\ \text{Birkhoff section} \\ \text{in $\mathcal{H}(2)$ stratum}}$}}} \arrow[rrrddd, "{\text{\cite{Deh}}}"{anchor=center,fill=white}] &  &  &  &  &  & |[fill=white, overlay, inner xsep=50pt, inner ysep=5pt]| \mathclap{\text{\large{$\substack{\text{Flows supported by} \\ \text{veering triangulations} \\ \text{with $\leq 16$ tetrahedra}}$}}} \arrow[lllddd, "{\text{\cite{ParSeg26}}}"{anchor=center,fill=white}] \\
 &  &  & \phantom{extra row} &  &  &  \\
 &  &  & \phantom{extra row} &  &  &  \\
 &  &  & \text{Anosov suspensions} &  &  &  
\end{tikzcd}
}

\caption{An overview of progress on the Fried--Christy--Ghys conjecture.}

 \label{fig:friedchristyghysprogress}

\end{figure}

Some progress has been made on \Cref{conj:friedchristyghys} in recent years. Among other results, we choose to highlight the following:
\begin{itemize}
    \item Asaoka--Bonatti--Marty \cite{ABM24} showed that every transitive Anosov flow can be transformed into a \textbf{skew} Anosov flow, i.e., one where the orbit space is a diagonal strip, via Goodman surgery. 
    \item Tsang \cite{Tsa24b} showed that every transitive totally periodic Anosov flow with orientable foliations can be transformed into an Anosov suspension flow via Goodman surgery. 
    Here we shall omit the definition of `totally periodic', since it is unimportant to this paper (see \cite{BF15} or \cite{BarMan25} for the definition).
    \item Fried \cite{Fri83} showed that every geodesic flow of a hyperbolic surface can be transformed into an Anosov suspension flow via Goodman surgery.
    \item Dehornoy \cite{Deh26} showed that every geodesic flow of a hyperbolic orbifold `of small type' can be transformed into an Anosov suspension flow via Goodman surgery. 
    \item Building on Dehornoy's ideas, Tsang \cite{Tsa26} showed that every transitive Anosov flow with orientable foliations on a graph manifold can be transformed into a totally periodic Anosov flow via Goodman surgery.
\end{itemize} 
We indicate these results and some others in \Cref{fig:friedchristyghysprogress}. 
In the figure, a solid arrow $A \to B$ should generally be read as `Every Anosov flow satisfying property $A$ can be transformed into an Anosov flow satisfying property $B$ via Goodman surgery'. 
Also, in the figure and throughout this paper, all Anosov flows will be assumed to be transitive and to have orientable foliations, unless otherwise stated.

The main result of this paper connects the upper and lower parts of \Cref{fig:friedchristyghysprogress} with the help of a new surgery operation.

\begin{thm} \label{thm:flowreducetogeneralizedht}
Every skew Anosov flow with orientable foliations can be transformed into a finite cover of a geodesic flow of a hyperbolic orbifold by \emph{bicontact padding surgery} followed by Goodman surgery.
\end{thm}

In fact, we can assume that the base orbifold of the resulting geodesic flow has at most one cone point, see the proof of \Cref{thm:flowreducetogeneralizedht} in \Cref{subsec:mainthmproof}.

The bicontact padding surgery is a new operation that we introduce in this article.
As the name suggests, it is an operation at the level of bicontact structures carrying Anosov flows, which, from a purely topological point of view, corresponds to Dehn surgery along certain curves.

Combining \Cref{thm:flowreducetogeneralizedht} with the results mentioned above, we have the following corollary, which gives a complete classification of transitive Anosov flows up to surgery operations, in the spirit of \Cref{conj:friedchristyghys}.

\begin{cor} \label{cor:flowreducetosusp}
Every transitive Anosov flow with orientable foliations can be transformed into an Anosov suspension flow via bicontact padding surgery and Goodman surgery. %
\end{cor}

The proof of \Cref{thm:flowreducetogeneralizedht} also gives us the following corollary, which can be thought of as a topological version of \Cref{conj:friedchristyghys}.

\begin{cor} \label{cor:dehnsurgerytosuspension} 
For any closed orientable 3-manifold $M$ and transitive Anosov flow $\phi$ with orientable foliations, there exists a finite sequence of manifolds $M = M_0, \dots, M_k$, and transitive Anosov flows $\phi = \phi_0, \dots, \phi_k$ such that $(M_k,\phi_k)$ is an Anosov suspension flow, and, for each $i$, $M_{i+1}$ is obtained from $M_i$ by doing a Dehn surgery on a link which is isotopic to a union of closed orbits of $\phi_i$.
\end{cor}

In the rest of this introduction, we will explain the new surgery operation and give an outline of the proof of \Cref{thm:flowreducetogeneralizedht}.

\subsection{Bicontact machinery} \label{subsec:introbicontactmachinery}

Let $M$ be an oriented 3-manifold. 
A \textbf{positive contact form} on $M$ is a 1-form $\alpha_+$ satisfying $\alpha_+ \wedge d\alpha_+ > 0$.
The \textbf{Reeb vector field} of $\alpha_+$ is defined to be the unique vector field $R_+$ satisfying $\iota_{R_+} \alpha_+ = 1$ and $\iota_{R_+} d\alpha_+ = 0$.
The \textbf{Reeb flow} of $\alpha_+$ is the flow generated by $R_+$.
Similarly, a \textbf{negative contact form} on $M$ is a 1-form $\alpha_-$ satisfying $\alpha_- \wedge d\alpha_- < 0$.
The \textbf{Reeb vector field} of $\alpha_-$ is defined to be the unique vector field $R_-$ satisfying $\iota_{R_-} \alpha_- = 1$ and $\iota_{R_-} d\alpha_- = 0$.
The \textbf{Reeb flow} of $\alpha_-$ is the flow generated by $R_-$.

\begin{figure}
    \centering
    \begingroup%
  \makeatletter%
  \providecommand\color[2][]{%
    \errmessage{(Inkscape) Color is used for the text in Inkscape, but the package 'color.sty' is not loaded}%
    \renewcommand\color[2][]{}%
  }%
  \providecommand\transparent[1]{%
    \errmessage{(Inkscape) Transparency is used (non-zero) for the text in Inkscape, but the package 'transparent.sty' is not loaded}%
    \renewcommand\transparent[1]{}%
  }%
  \providecommand\rotatebox[2]{#2}%
  \newcommand*\fsize{\dimexpr\f@size pt\relax}%
  \newcommand*\lineheight[1]{\fontsize{\fsize}{#1\fsize}\selectfont}%
  \ifx\svgwidth\undefined%
    \setlength{\unitlength}{100.70331123bp}%
    \ifx\svgscale\undefined%
      \relax%
    \else%
      \setlength{\unitlength}{\unitlength * \real{\svgscale}}%
    \fi%
  \else%
    \setlength{\unitlength}{\svgwidth}%
  \fi%
  \global\let\svgwidth\undefined%
  \global\let\svgscale\undefined%
  \makeatother%
  \begin{picture}(1,1.40594073)%
    \lineheight{1}%
    \setlength\tabcolsep{0pt}%
    \put(0,0){\includegraphics[width=\unitlength,page=1]{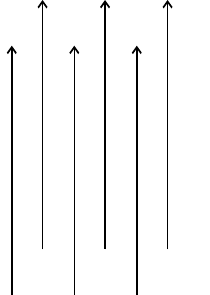}}%
    \put(0.93519575,0.67704338){\color[rgb]{0,0,1}\makebox(0,0)[lt]{\lineheight{1.25}\smash{\begin{tabular}[t]{l}$\xi_+$\end{tabular}}}}%
    \put(0,0){\includegraphics[width=\unitlength,page=2]{contactstructurerudder.pdf}}%
  \end{picture}%
\endgroup%

    \caption{A cooriented plane field is a positive contact structure if and only if it rotates around a nonvanishing tangent vector field in a left-handed way.}
    \label{fig:contactstructurerudder}
\end{figure}

A \textbf{positive/negative contact structure} is a cooriented plane field that can be expressed as the kernel of a positive/negative contact form, respectively.
Equivalently, a cooriented plane field is positive/negative contact if and only if it rotates in a left-/right-handed way around a nonvanishing tangent vector field, respectively.
See \Cref{fig:contactstructurerudder}, and see \Cref{prop:contactstructurerudder} for a more precise statement.

\begin{figure}
    \centering
    \begingroup%
  \makeatletter%
  \providecommand\color[2][]{%
    \errmessage{(Inkscape) Color is used for the text in Inkscape, but the package 'color.sty' is not loaded}%
    \renewcommand\color[2][]{}%
  }%
  \providecommand\transparent[1]{%
    \errmessage{(Inkscape) Transparency is used (non-zero) for the text in Inkscape, but the package 'transparent.sty' is not loaded}%
    \renewcommand\transparent[1]{}%
  }%
  \providecommand\rotatebox[2]{#2}%
  \newcommand*\fsize{\dimexpr\f@size pt\relax}%
  \newcommand*\lineheight[1]{\fontsize{\fsize}{#1\fsize}\selectfont}%
  \ifx\svgwidth\undefined%
    \setlength{\unitlength}{120.60200609bp}%
    \ifx\svgscale\undefined%
      \relax%
    \else%
      \setlength{\unitlength}{\unitlength * \real{\svgscale}}%
    \fi%
  \else%
    \setlength{\unitlength}{\svgwidth}%
  \fi%
  \global\let\svgwidth\undefined%
  \global\let\svgscale\undefined%
  \makeatother%
  \begin{picture}(1,1.3080725)%
    \lineheight{1}%
    \setlength\tabcolsep{0pt}%
    \put(0,0){\includegraphics[width=\unitlength,page=1]{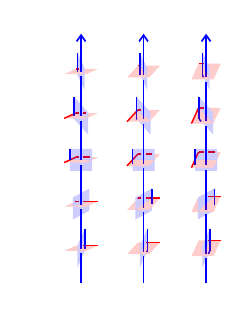}}%
    \put(0.03115543,0.91682039){\color[rgb]{0,0,1}\makebox(0,0)[lt]{\lineheight{1.25}\smash{\begin{tabular}[t]{l}$\xi_+$\end{tabular}}}}%
    \put(0.50066786,1.23255139){\color[rgb]{0,0,1}\makebox(0,0)[lt]{\lineheight{1.25}\smash{\begin{tabular}[t]{l}$R_-$\end{tabular}}}}%
    \put(0.94447856,0.58809301){\color[rgb]{0,0,0}\makebox(0,0)[lt]{\lineheight{1.25}\smash{\begin{tabular}[t]{l}$X$\end{tabular}}}}%
    \put(-0.00616908,0.41782755){\color[rgb]{1,0,0}\makebox(0,0)[lt]{\lineheight{1.25}\smash{\begin{tabular}[t]{l}$\xi_-$\end{tabular}}}}%
    \put(0,0){\includegraphics[width=\unitlength,page=2]{stronglyadaptedbicontact.pdf}}%
  \end{picture}%
\endgroup%

    \caption{A local picture of a strongly adapted bicontact package.}
    \label{fig:stronglyadaptedbicontact}
\end{figure}

For the purposes of this paper, 
we define a \textbf{bicontact package} to be a pair $(\alpha_-,\xi_+)$ consisting of a negative contact \emph{form} $\alpha_-$ and a positive contact \emph{structure} $\xi_+$.
We say that $(\alpha_-,\xi_+)$ is \textbf{strongly adapted} if the Reeb vector field $R_-$ of $\alpha_-$ is tangent to $\xi_+$ everywhere.
We say that $(\alpha_-,\xi_+)$ \textbf{supports} a flow $\phi$ if $\phi$ is generated by a nonvanishing vector field that lies in the intersection of $\xi_- = \ker \alpha_-$ and $\xi_+$ at every point.
See \Cref{fig:stronglyadaptedbicontact} for a local illustration of such a scenario.

A theorem of Hozoori (\cite[Theorem A]{Hoz25}) states that a flow is Anosov with orientable foliations if and only if it is supported by a strongly adapted bicontact package.
Via this theorem, one can define operations on Anosov flows using strongly adapted bicontact packages.
This approach has the advantage that constructing bicontact packages can be done locally, while certifying that a flow is Anosov generally requires a global argument.

For example, Salmoiraghi \cite{Sal25} defined an operation known as \textbf{Legendrian-transverse surgery}:
Given a bicontact package $(\alpha_-,\xi_+)$, a \textbf{Legendrian-transverse curve} is a curve that is \textbf{Legendrian}, i.e., tangent, to $\xi_- = \ker \alpha_-$ and transverse to $\xi_+$, or Legendrian to $\xi_+$ and transverse to $\xi_-$.
One can construct new Anosov flows by modifying a supporting bicontact package in a neighborhood of a Legendrian-transverse curve. See \Cref{constr:legendriantransversesurgery}.
This operation generalizes (a single) Goodman surgery, although it is shown  in \cite{Tsa24a} that Legendrian-transverse surgery on a transitive Anosov flow can be realized by a sequence of Goodman surgeries.

If one allows all Legendrian-transverse surgeries, then one can even remove the transitivity assumption in \Cref{cor:flowreducetosusp}. That is, we will also show: 
\begin{thm}\label{thm:LT_and_bicontact_gives_suspension}
    Let $\phi$ be a (not necessarily transitive) Anosov flow on $M$ with orientable foliations. Then $\phi$ can be transformed into a suspension by Legendrian-transverse and bicontact padding surgeries.
\end{thm}

\subsection{A failed approach to \Cref{conj:friedchristyghys}} \label{subsec:intropaddingmotivation}

One way to motivate our new \textit{bicontact padding surgery} is to explain a failed approach to proving \Cref{conj:friedchristyghys}:

Let $\phi$ be an Anosov flow with orientable foliations. Let $(\alpha_-,\xi_+)$ be a strongly adapted bicontact package that supports $\phi$. By \cite[Corollary 1.2]{CDHR24}, up to a perturbation of $\alpha_-$, we can assume that its Reeb flow $R_-$ admits a \textbf{Birkhoff section}.
This is a compact cooriented immersed surface $S$ satisfying:
\begin{itemize}
    \item the interior of $S$ is embedded and positively transverse to $R_-$,
    \item each boundary component of $S$ lies along a closed orbit of $R_-$, and
    \item every forward and backward trajectory of $R_-$ meets $S$.
\end{itemize}
In particular, the restriction of $R_-$ to the complement of the boundary orbits of $S$ is the suspension flow of some diffeomorphism $f$ on $\intr(S)$.

For each $x \in \intr(S)$, consider the family of lines $\left( (R_-^{-t})_*(\xi_+|_{R_-^t(x)}) \cap T_xS \right)_{t \in \mathbb{R}}$ on $T_x S$. By definition of positive contact structures, these lines rotate clockwise monotonically.
If $\phi$ is positive skew, then, based on Thurston's \textit{slithering} picture, one expects that these lines complete infinitely many rotations as $t \to \pm \infty$. (See also the geodesic flow example in the next subsection.) 

\begin{figure}
    \centering
    \fontsize{8pt}{8pt}\selectfont
    \begingroup%
  \makeatletter%
  \providecommand\color[2][]{%
    \errmessage{(Inkscape) Color is used for the text in Inkscape, but the package 'color.sty' is not loaded}%
    \renewcommand\color[2][]{}%
  }%
  \providecommand\transparent[1]{%
    \errmessage{(Inkscape) Transparency is used (non-zero) for the text in Inkscape, but the package 'transparent.sty' is not loaded}%
    \renewcommand\transparent[1]{}%
  }%
  \providecommand\rotatebox[2]{#2}%
  \newcommand*\fsize{\dimexpr\f@size pt\relax}%
  \newcommand*\lineheight[1]{\fontsize{\fsize}{#1\fsize}\selectfont}%
  \ifx\svgwidth\undefined%
    \setlength{\unitlength}{332.95155202bp}%
    \ifx\svgscale\undefined%
      \relax%
    \else%
      \setlength{\unitlength}{\unitlength * \real{\svgscale}}%
    \fi%
  \else%
    \setlength{\unitlength}{\svgwidth}%
  \fi%
  \global\let\svgwidth\undefined%
  \global\let\svgscale\undefined%
  \makeatother%
  \begin{picture}(1,0.60631376)%
    \lineheight{1}%
    \setlength\tabcolsep{0pt}%
    \put(0,0){\includegraphics[width=\unitlength,page=1]{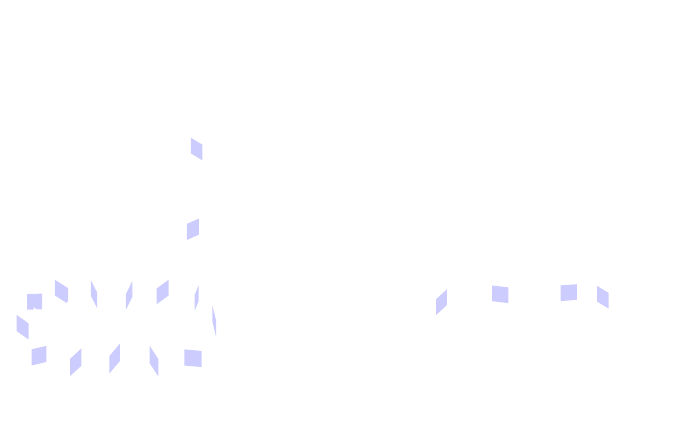}}%
    \put(0.30578108,0.09238606){\color[rgb]{0,0,1}\makebox(0,0)[lt]{\lineheight{1.25}\smash{\begin{tabular}[t]{l}$\xi_+$\end{tabular}}}}%
    \put(0,0){\includegraphics[width=\unitlength,page=2]{homotopeforlegendrian2.pdf}}%
    \put(0.0643497,0.33459566){\color[rgb]{0,0,0}\makebox(0,0)[lt]{\lineheight{1.25}\smash{\begin{tabular}[t]{l}$R_-^{t(x)}(x)$\end{tabular}}}}%
    \put(0.17097209,0.06521211){\color[rgb]{0,0,0}\makebox(0,0)[lt]{\lineheight{1.25}\smash{\begin{tabular}[t]{l}$x$\end{tabular}}}}%
    \put(0.24017985,0.06222376){\color[rgb]{0,0,0}\makebox(0,0)[lt]{\lineheight{1.25}\smash{\begin{tabular}[t]{l}$c$\end{tabular}}}}%
    \put(0,0){\includegraphics[width=\unitlength,page=3]{homotopeforlegendrian2.pdf}}%
  \end{picture}%
\endgroup%

    \caption{We isotope $\xi_+$ along flow lines of $R_-$ so that $c$ is tangent to $\xi_+$.}
    \label{fig:homotopeforlegendrian}
\end{figure}

Now suppose we are given a curve $c \subset \intr(S)$. 
A standard result in contact topology, known as the \textit{Legendrian realization principle}, states that up to isotoping $S$, we can assume that $c$ is Legendrian for $\xi_- = \ker \alpha_-$.
Then for each $x \in c$, since the lines $\left( (R_-^{-t})_*(\xi_+|_{R_-^t(x)}) \cap T_xS \right)_{t \in \mathbb{R}}$ complete infinitely many rotations, there exists a $t(x)$ such that $(R_-^{-t(x)})_*(\xi_+|_{R_-^{t(x)}(x)}) \cap T_xS = T_x c$. 
Hence, one can hope to isotope $\xi_+$ along flow lines of $R_-$, dragging each $\xi_+|_{R_-^{t(x)}(x)}$ down to $x$, until $c$ is tangent to $\xi_+$.
See \Cref{fig:homotopeforlegendrian}.

If possible, this process would determine a smooth family of strongly adapted bicontact packages, supporting a smooth family of Anosov flows, which, by structural stability, stays within one isotopy equivalence class.
At the end of the isotopy, $c$ would be tangent to both $\xi_-$ and $\xi_+$, and thus would be a closed orbit of the Anosov flow.
We would therefore be able to perform Goodman surgery along $c$. In doing so, we change the monodromy of the Birkhoff section $S$ by composing it with a Dehn twist along $c$.

Since the (pure) mapping class group of any surface is generated by Dehn twists, we can hope to repeat this process and reduce the monodromy to the identity. This would force the underlying 3-manifold to be a Seifert fibered space, and hence force the flow to be equivalent to a lift of a geodesic flow by \cite[Theorem 4.9]{Bar96}.
Summarizing, this argument would show that every skew Anosov flow can be transformed into a lift of a geodesic flow via Goodman surgery.

\begin{figure}
    \centering
    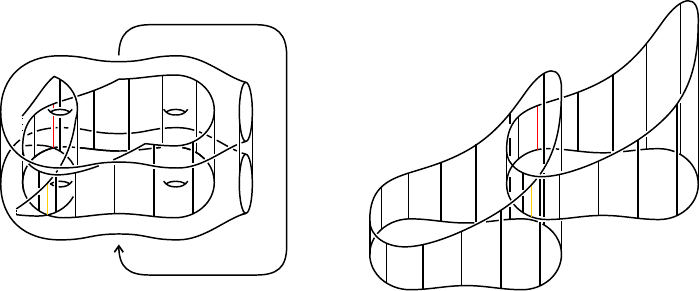
    \caption{In general, the `pullback annulus' $A$ could intersect itself. Left: The self-intersection as seen in $M$. Right: The self-intersection lifted to a cover, for ease of viewing. A `ribbon singularity', i.e., a self-intersection as in the red arc in the figure, obstructs an isotopy of $\xi_+$ along flow lines of $R_-$.}
    \label{fig:annulusribbonsing}
\end{figure}

Unfortunately, the argument has (at least) the following gaps:
\begin{enumerate}
    \item The standard statement of the Legendrian realization principle does not ensure that $S$ remains transverse to $R_-$ after isotopy.
    In more detail, the usual proof involves constructing the isotopy inside an abstract copy of $S \times \mathbb{R}$, then embedding $S \times \mathbb{R}$, and thus the isotopy, inside the 3-manifold. 
    In general, this embedding is only an embedding of contact structures, but not an embedding of contact forms, so the direction of the Reeb flow can be distorted.    
    \item The `pullback time' function $t\colon c \to \mathbb{R}$ can only be continuous if the winding number of $c$ with respect to the line field $\ell = \xi_+ \cap TS$ on $\intr(S)$ is zero, so one can only hope to change the monodromy of $S$ by composing it with Dehn twists along such winding number zero curves.
    On the other hand, such Dehn twists do not, in general, generate the whole mapping class group. 
    For example, using a homological argument, one can show that a Dehn twist along a non-separating curve with non-zero winding number cannot be written as a product of Dehn twists along winding number zero curves.
    \item Consider the set $A$ that is the union of flow lines from each $x \in c$ to $R_-^{t(x)}(x)$. To simplify the discussion, let us assume that $t(x) > 0$ for each $x$, so that $A$ is an annulus. 
    In general, $A$ is only an immersed annulus and not an embedded one, i.e., it could intersect itself. This is more likely to happen if $S$ has a small first return time. See \Cref{fig:annulusribbonsing} for an example where $A$ intersects itself in two arcs, one orange and one red.
    In such a case, it is unclear whether an isotopy dragging each $R_-^{t(x)}(x)$ down to $x$ should exist. In fact, if $A$ has any `ribbon singularities' as in the red arc in \Cref{fig:annulusribbonsing}, then such an isotopy cannot exist, since the isotopy would force a self-intersection of $c$.
\end{enumerate}

\subsection{Bicontact padding surgery} \label{subsec:intropaddingdescription}

Bicontact padding surgery resolves issues (1) and (3). The idea of the surgery is to insert a `padding layer' $P$ along $S$ where $\xi_+$ completes many rotations around $R_-$.
By ensuring that $P$ has enough `height' with respect to $R_-$, we can embed the isotopy given by the Legendrian realization principle while preserving transversality to $R_-$.
Meanwhile, by inserting enough rotations of $\xi_+$ around $R_-$ within $P$, we can arrange for the `pullback annulus' $A$ to lie entirely within $P$, and thus be embedded, ensuring that we can perform the isotopy of $\xi_+$ as desired.

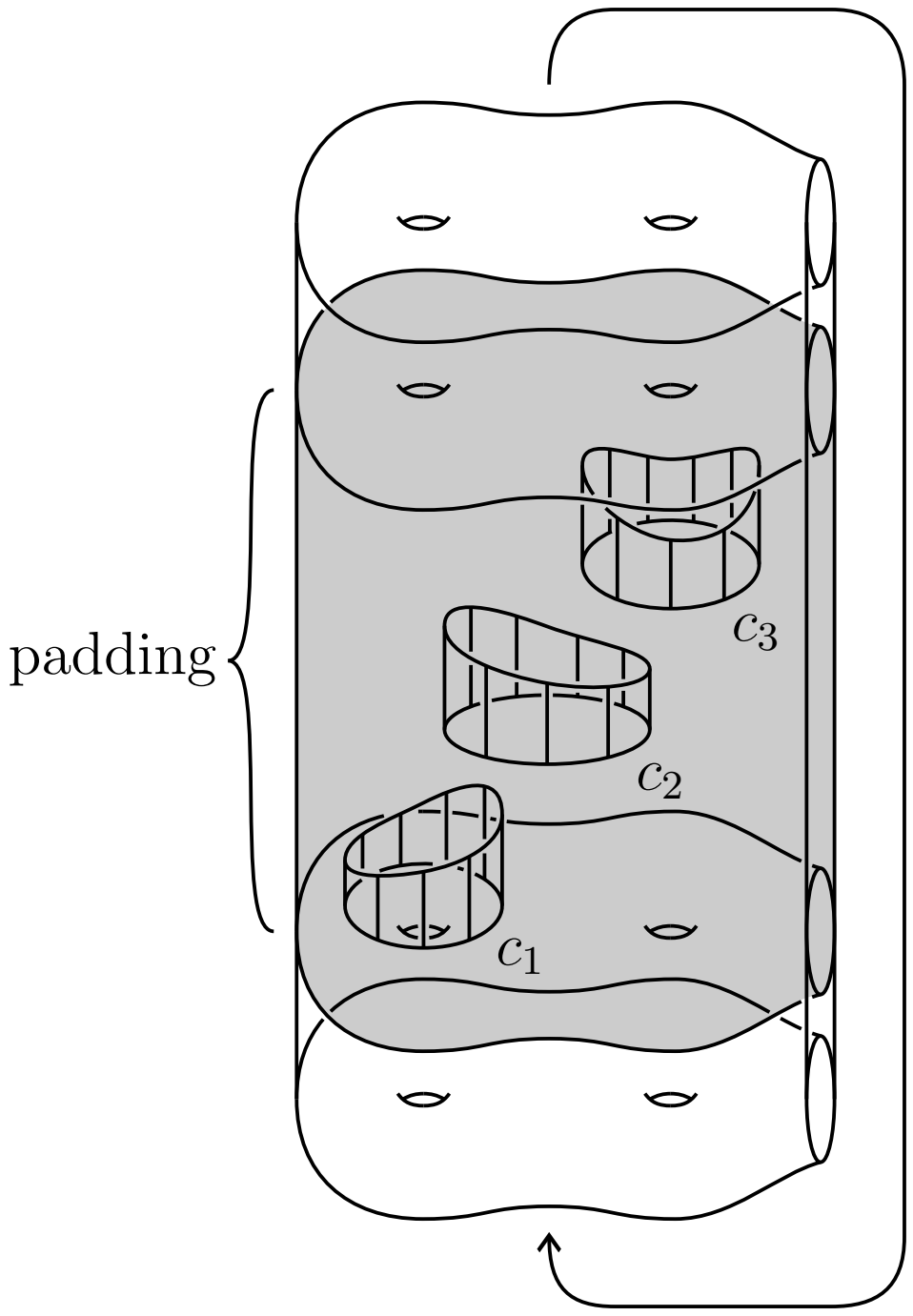
\begin{figure}
    \centering
    \begingroup%
  \makeatletter%
  \providecommand\color[2][]{%
    \errmessage{(Inkscape) Color is used for the text in Inkscape, but the package 'color.sty' is not loaded}%
    \renewcommand\color[2][]{}%
  }%
  \providecommand\transparent[1]{%
    \errmessage{(Inkscape) Transparency is used (non-zero) for the text in Inkscape, but the package 'transparent.sty' is not loaded}%
    \renewcommand\transparent[1]{}%
  }%
  \providecommand\rotatebox[2]{#2}%
  \newcommand*\fsize{\dimexpr\f@size pt\relax}%
  \newcommand*\lineheight[1]{\fontsize{\fsize}{#1\fsize}\selectfont}%
  \ifx\svgwidth\undefined%
    \setlength{\unitlength}{175.81504161bp}%
    \ifx\svgscale\undefined%
      \relax%
    \else%
      \setlength{\unitlength}{\unitlength * \real{\svgscale}}%
    \fi%
  \else%
    \setlength{\unitlength}{\svgwidth}%
  \fi%
  \global\let\svgwidth\undefined%
  \global\let\svgscale\undefined%
  \makeatother%
  \begin{picture}(1,1.45554573)%
    \lineheight{1}%
    \setlength\tabcolsep{0pt}%
    \put(0,0){\includegraphics[width=\unitlength,page=1]{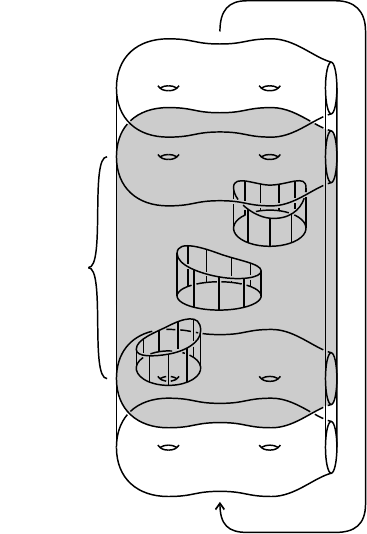}}%
    \put(-0.00483465,0.71003371){\color[rgb]{0,0,0}\makebox(0,0)[lt]{\lineheight{1.25}\smash{\begin{tabular}[t]{l}padding\end{tabular}}}}%
    \put(0.5408906,0.38436194){\color[rgb]{0,0,0}\makebox(0,0)[lt]{\lineheight{1.25}\smash{\begin{tabular}[t]{l}$c_1$\end{tabular}}}}%
    \put(0.69766951,0.58039438){\color[rgb]{0,0,0}\makebox(0,0)[lt]{\lineheight{1.25}\smash{\begin{tabular}[t]{l}$c_2$\end{tabular}}}}%
    \put(0.80462877,0.74705944){\color[rgb]{0,0,0}\makebox(0,0)[lt]{\lineheight{1.25}\smash{\begin{tabular}[t]{l}$c_3$\end{tabular}}}}%
  \end{picture}%
\endgroup%

    \caption{Inserting a padding layer allows us to change the monodromy on a Birkhoff section for $R_-$.}
    \label{fig:proofstratpadding}
\end{figure}

Up to increasing the height and number of rotations even more, we can fit multiple `pullback annuli' inside the padding layer, see \Cref{fig:proofstratpadding}.
This allows us to change the monodromy by any product of Dehn twists along winding number zero curves. We record this as the following theorem, which may be of independent interest.

\begin{thm} \label{thm:paddingaddwinding0twists}
Let $(\alpha_-,\xi_+)$ be a strongly adapted bicontact package on a closed oriented 3-manifold $M$. Let $S$ be a Birkhoff section to $R_-$ whose boundary orbits are all negative, and let $f$ be the monodromy of $S$.
Then for any finite (possibly repeating) sequence of curves $(c_1,...,c_m)$ on $S$ with $\wind_{\xi_+}(S,c_i) = 0$ for each $i=1,...,m$,
and any sequence of integers $(n_1,...,n_m)$,
there exists a bicontact package $(\widehat{\alpha}_-,\widehat{\xi}_+)$, obtained from $(\alpha_-,\xi_+)$ by bicontact padding surgery, then homotopy, then Goodman surgery, for which there is a Birkhoff section $\widehat{S}$ for the Reeb flow $\widehat{R}_-$, such that the monodromy of $\widehat{S}$ is $\tau_{c_1}^{n_1} ... \tau_{c_m}^{n_m} f$.
\end{thm}

\subsection{Toy model: Geodesic flow} \label{subsec:introgeodflow}

We give more intuition for the bicontact padding surgery by describing it on the model example of the geodesic flow: 
Let $\Sigma$ be a closed hyperbolic surface and let $g^t$ be the geodesic flow on $T^1\Sigma$.

One can construct a strongly adapted bicontact package $(\alpha_-,\xi_+)$ on $T^1 \Sigma$ where the Reeb flow $R_-$ is exactly the Seifert fibering and $\xi_+$ is spanned by the directions of the Seifert fiber and the geodesic flow. See \cite[Section 3.1]{Sal25} for details.
A Birkhoff section to $R_-$ can be constructed by taking a section to the Seifert fibering in the complement of finitely many points, which is in turn equivalent to a nonsingular vector field on $\Sigma$ in the complement of finitely many points.

Notice that for such a Birkhoff section $S$, $\xi_+$ completes exactly one rotation before $S$ bumps into itself under $R_-$.

Now on $\intr(S) \times I \subset T^1\Sigma$, we modify the bicontact package by adding a padding layer where $\xi_+$ completes $k-1$ rotations around $R_-$. This gives a bicontact package on a $3$-manifold with torus boundary components, and one can glue to each boundary torus the neighborhood of an orbit of $R_-$ where $\xi_+$ is also modified so that it completes $k-1$ additional rotations. This operation yields a new bicontact package on a (generally) new manifold.

It turns out that if $S$ was chosen to have a single boundary component and if $k$ divides $\chi(\Sigma)$, then one can show that this operation yields the same bicontact package as that of the lifted geodesic flow to the $k$-fold cover of $T^1\Sigma$. We also refer to \cite{PinSal26} where this construction features prominently in the case when $\Sigma$ is a surface with boundary.

The key simplification in this toy model is that the bicontact package is particularly nice near the boundary orbits. We call such nice boundary orbits \textbf{standard}. Since bicontact packages are at least as rigid as Reeb flows near closed orbits, one cannot expect that a general bicontact package will always be isomorphic to that of a geodesic flow near the boundary orbits of a Birkhoff section. 
In order to be able to define bicontact padding surgery in full generality, we use \Cref{prop:standardizeorbit}, which states that up to homotoping the bicontact packages, one can always assume that each boundary orbit is standard. %

\subsection{Reducing the monodromy}

It remains to overcome issue (2). We first observe that the monodromy $f$ preserves the homotopy class of $\ell$, in the sense that there is a continuous family of line fields interpolating between $f_* \ell$ and $\ell$, obtained by flowing $\intr(S)$ backwards along $R_-$ and taking its intersection with $\xi_+$ in the process.

We now have a theorem of Calderon--Salter in \cite{CS23}, which morally says that if $S$ has genus $\geq 5$, then twists along winding number zero curves generate a large subgroup in the group of diffeomorphisms that preserve the homotopy class of $\ell$, see \Cref{thm:factorizationthm}.
In fact, if $S$ has only one boundary component, then these twists generate the whole group, see \Cref{cor:factorizationthmoneboundary}.

Thus the issue can be overcome if we can construct a Birkhoff section $S$ to $R_-$ that has genus $\geq 5$ and one boundary component.
To do so we use the assumption that $\phi$ is negative skew:
A theorem of Marty \cite[Theorem F]{Mar25} states that every skew Anosov flow admits a Birkhoff section with only one boundary component. Inspecting his proof, it is not difficult to see that the genus of the Birkhoff section can be made arbitrarily large.
Meanwhile, an argument of Barthelmé \cite{Bar25} shows that $\phi$ is isotopically equivalent to $R_-$ of a supporting bicontact package, see \Cref{prop:swapanosovR_-}.
This allows us to transfer a Birkhoff section for $\phi$ with the desired properties to one for $R_-$.

\subsection{Outline of proof} \label{subsec:introoutline}

Having explained how to overcome each of the issues in the failed approach, let us summarize this discussion with an outline of the proof of \Cref{thm:flowreducetogeneralizedht} in a logical order.

Let $\phi$ be a negative skew Anosov flow. We construct a Birkhoff section for $\phi$ that has genus $\geq 5$ and only one boundary component.
We then choose a supporting strongly adapted bicontact package $(\alpha_-,\xi_+)$ where the Reeb flow $R_-$ is isotopically equivalent to $\phi$.
In particular, $R_-$ also has a Birkhoff section $S$ that has genus $\geq 5$ and only one boundary component.

We then apply \Cref{thm:paddingaddwinding0twists} and \Cref{cor:factorizationthmoneboundary} to transform the bicontact package via bicontact padding surgery and Goodman surgery into one that has a Birkhoff section with trivial monodromy.
This implies that the resulting Anosov flow lives on a Seifert fibered space. Thus by a result of Barbot \cite{Bar96}, the flow must be a finite cover of a geodesic flow.

\subsection{AI disclosure} \label{subsec:aidisclosure} 

An LLM proofreading review of an earlier draft of this article detected an error in our construction for standardizing boundary orbits, which is a preliminary step we need in order to define our padding surgery in full generality. 
The model (ChatGPT-6 Astra, using a ``plus'' subscription) was then asked for a fix\footnote{The prompt and the model's output are available on the authors' websites, or at this link \url{https://drive.google.com/drive/folders/1N_-YBmhJxkiQsguL5m_2lHJB72LLyFGS?usp=sharing}.} and came up with a different approach which is now \Cref{prop:standardizeorbit}.
This LLM-generated proof was digested then rewrote and presented in this article in a self-contained way and with some clarifying details, albeit keeping the key ideas.

We made an effort to trace the source of the ideas in the proof:
An (AI-aided) literature search reveals that some ideas appear to be taken from the study of the \textit{inverse problem of the calculus of variations}, on which a standard reference seems to be \cite{Kru97}, although the exact argument seems to be new. 
For this reason, we have attributed \Cref{prop:standardizeorbit} to the model used.
However, the authors are insufficiently acquainted with this area of mathematics to feel confident about this judgment.
We therefore invite comments on this point from anyone more familiar with the literature.

Aside from this argument, all other parts of the paper were conceived by the authors.

\subsection*{Acknowledgments}

The authors thank Pierre Dehornoy and Ana Rechtman for helpful discussions.
Some progress was made on this project at the 2024 St.~Louis Topology Conference and the 2025 `Low dimensional topology and Floer theory' workshop in Montr\'eal. We thank the organizers of these events for creating stimulating research environments. 
TB was partially supported by the NSERC Discovery (RGPIN-2024-04412) and Alliance International programs (ALLRP 598447 - 24).
CCT was supported by a CRM postdoctoral fellowship at CIRGET.

\subsection*{Conventions}

In this paper, we will deal with different objects with various notions of signs and orientations. As such, we take the time to specify our conventions.

On manifolds:
\begin{itemize}
    \item Unless otherwise stated, all manifolds in this paper will be oriented, smooth, compact, and connected.
    \item A \textbf{curve} is an a priori unoriented, embedded closed dimension 1 submanifold. (In particular, for us, all curves on surfaces are automatically simple.)
    An \textbf{arc} is an a priori unoriented, properly embedded non-closed dimension 1 submanifold.
    \item Unless stated otherwise, we will write $S^1 = \mathbb{R}/2\pi\mathbb{Z}$, $T^2 = \mathbb{R}^2/(2\pi\mathbb{Z})^2$, and $D^2$ for the unit disc.
\end{itemize}

On orientations:
\begin{itemize}
    \item If $P$ and $Q$ are two oriented submanifolds of $M$ that intersect transversely, then each component of $P \cap Q$ is also an oriented submanifold, where the orientation is defined by the `intersection last' convention.
    That is, if $(x_1,...,x_n)$ is an oriented local coordinate system on $M$ 
    where $(x_1,...,x_{n-q},x_{2n-p-q+1},...,x_n)$ is an oriented local coordinate system on $P$ and $(x_{n-q+1},...,x_{2n-p-q},x_{2n-p-q+1},...,x_n)$ is an oriented local coordinate system on $Q$, 
    then $(x_{2n-p-q+1},...,x_n)$ is an oriented local coordinate system on $P \cap Q$.
    \item If $N$ is a codimension-1 submanifold of $M$, then a coorientation on $N$ induces an orientation on $N$, and vice versa, by specifying that
    for every oriented 1-manifold $c$ transverse to $N$, the sign of each point in $N \cap c$ is consistent when determined using either the orientation or coorientation of $N$.
    \item Unless otherwise stated, each boundary component of an oriented manifold $M$ will be given its \textbf{boundary orientation},
    defined as the orientation induced by the inward coorientation.
    Equivalently, the boundary orientation is defined by the `inward normal last' convention.
\end{itemize}
These conventions capture many of the usual intuitive conventions in dimensions $\leq 3$. 
For example: 
\begin{itemize}[label=-]
    \item The boundary of the unit disc in $\mathbb{R}^2$ is oriented counterclockwise. 
    \item The intersection of two cooriented planes in $\mathbb{R}^3$ is oriented by taking the right-hand rule with their normals.
\end{itemize}

On flows:
\begin{itemize}
    \item A \textbf{flow} on a manifold $M$ is a smooth function $\phi\colon \mathbb{R} \times M \to M$, for which we write $\phi^t(x) = \phi(t,x)$, satisfying $\phi^0(x) = x$ and $\phi^s(\phi^t(x)) = \phi^{s+t}(x)$ for every $x \in M, s, t \in \mathbb{R}$. 
    \item The \textbf{generating vector field} $X$ of $\phi$ is defined by $X|_x = \frac{d}{dt}|_{t=0} \phi^t(x)$ at each point $x \in M$.
    \item Conversely, every smooth vector field $X$ generates a flow. 
    When there is no risk of confusion, we will save on assigning the flow a separate name and instead denote it also by $X$.
    \item Unless otherwise stated, all the flows in this paper have no stationary points, or equivalently their generating vector fields are nonvanishing.  
\end{itemize}

On Dehn twists:
\begin{itemize}
    \item Let $c$ be a curve on a surface $S$. Pick a collar neighborhood $\nu$ of $c$, along with an oriented coordinate system $(r,\theta) \in [0,1] \times \mathbb{R}/\mathbb{Z}$. 
    A \textbf{Dehn twist} of degree $n \in \mathbb{Z}$ along $c$ is a homeomorphism of $S$ supported in $\nu$ of the form $(r,\theta) \mapsto (r,\theta+h(r))$ where $h = 0$ near $r = 0$ and $h = n$ near $r = 1$.
    \item A \textbf{positive Dehn twist} along $c$ is a Dehn twist of degree $1$ along $c$, and will be denoted by $\tau_c$.
\end{itemize}

\begin{figure}
    \centering
    \begingroup%
  \makeatletter%
  \providecommand\color[2][]{%
    \errmessage{(Inkscape) Color is used for the text in Inkscape, but the package 'color.sty' is not loaded}%
    \renewcommand\color[2][]{}%
  }%
  \providecommand\transparent[1]{%
    \errmessage{(Inkscape) Transparency is used (non-zero) for the text in Inkscape, but the package 'transparent.sty' is not loaded}%
    \renewcommand\transparent[1]{}%
  }%
  \providecommand\rotatebox[2]{#2}%
  \newcommand*\fsize{\dimexpr\f@size pt\relax}%
  \newcommand*\lineheight[1]{\fontsize{\fsize}{#1\fsize}\selectfont}%
  \ifx\svgwidth\undefined%
    \setlength{\unitlength}{168.72991967bp}%
    \ifx\svgscale\undefined%
      \relax%
    \else%
      \setlength{\unitlength}{\unitlength * \real{\svgscale}}%
    \fi%
  \else%
    \setlength{\unitlength}{\svgwidth}%
  \fi%
  \global\let\svgwidth\undefined%
  \global\let\svgscale\undefined%
  \makeatother%
  \begin{picture}(1,0.44265855)%
    \lineheight{1}%
    \setlength\tabcolsep{0pt}%
    \put(0,0){\includegraphics[width=\unitlength,page=1]{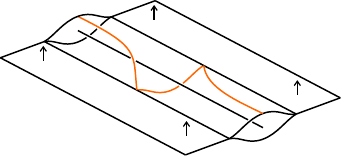}}%
    \put(0.78651805,0.00562626){\color[rgb]{0,0,0}\makebox(0,0)[lt]{\lineheight{1.25}\smash{\begin{tabular}[t]{l}$c$\end{tabular}}}}%
    \put(0.96119243,0.08982802){\color[rgb]{0,0,0}\makebox(0,0)[lt]{\lineheight{1.25}\smash{\begin{tabular}[t]{l}$S$\end{tabular}}}}%
  \end{picture}%
\endgroup%

    \caption{A $+1$-surgery on a curve $c$ on a fiber surface $S$ changes the monodromy from $f$ to $\tau_c f$.}
    \label{fig:mappingtorusdehnsurgery}
\end{figure}

On slopes and Dehn surgery:
\begin{itemize}
    \item Let $c$ be a curve in a 3-manifold $M$. Pick a tubular neighborhood $\nu$ of $c$, along with an oriented cylindrical coordinate system $(r,\theta,z)$.
    A \textbf{slope} on the torus $\partial \nu$ is a projective equivalence class in $H_1(\partial \nu; \mathbb{R})$.
    The \textbf{meridian} is the slope that vanishes in $H_1(\nu;\mathbb{R})$.
    \item Orienting $\partial \nu$ as the boundary of $M \backslash \nu$, there is an induced orientation on $H_1(\partial \nu; \mathbb{R})$, thus an induced clockwise orientation on the set of slopes $PH_1(\partial \nu; \mathbb{R}) \cong S^1$.
    Given two slopes $s_1,s_2$, we write $(s_1,s_2)$ for the open interval containing all slopes $s$ for which $(s_1,s,s_2)$ is ordered positively. 
    We will also write $s_1<s<s_2$ for $s \in (s_1,s_2)$.
    \item A slope is \textbf{rational} if it is a projective equivalence class in $H_1(\partial \nu; \mathbb{Q})$.
    Equivalently, the rational slopes are exactly the isotopy classes of essential closed curves on $\partial \nu$.
    The intersection number between two rational slopes is the geometric intersection number between the corresponding curves.
    \item A \textbf{framing} on $c$ is a rational slope on $\partial \nu$ that has intersection number $1$ with the meridian. 
    Given a framing, the set of slopes can be canonically identified with $\mathbb{R}P^1 = \mathbb{R} \cup \{\infty\}$, where the framing corresponds to $0$ and the meridian corresponds to $\infty$.
    The set of rational slopes is then identified with $\mathbb{Q}P^1 = \mathbb{Q} \cup \{\infty\}$. We denote elements of the latter by fractions $\frac{p}{q}$, where $p,q \in \mathbb{Z}$ satisfy $\gcd(p,q) = 1$.
    \item With a framing specified, the \textbf{$\frac{p}{q}$-Dehn surgery} of $M$ along $c$ is obtained by cutting and regluing $\nu$ to $M \backslash \nu$ in a way so that the meridian on $\nu$ is sent to a curve of slope $\frac{p}{q}$ on $M \backslash \nu$.
\end{itemize}
In particular, if $M = S \times [0,1] / (x,1) \sim (f(x),0)$ is the mapping torus of a surface homeomorphism $f\colon S \to S$, and $c$ is a curve on the fiber surface $S \times \{0\}$, then $\frac{1}{n}$-surgery of $M$ along $c$ gives the mapping torus of $\tau_c^n f$, where we take the natural framing of $c$ induced by the fiber surface $S \times \{0\}$. %
See \Cref{fig:mappingtorusdehnsurgery}.

\section{Preliminaries} \label{sec:prelim}

\subsection{Anosov flows} \label{subsec:anosovflow}

Throughout this paper, $M$ will be a closed oriented smooth 3-manifold.
A flow $\phi$ on $M$ generated by a nonvanishing vector field $X$ is \textbf{Anosov} if there exists a $d\phi^t$-invariant splitting $TM = E^s \oplus \langle X \rangle \oplus E^u$ such that under some Riemannian metric on $M$, there exists $C, r > 0$ such that $||d\phi^t(v^s)|| \leq Ce^{-rt} ||v^s||$ for every $t > 0$ and $v^s \in E^s$, and $||d\phi^t(v^u)|| \leq Ce^{rt} ||v^u||$ for every $t < 0$ and $v^u \in E^u$. %

It can be shown that $E^{s/u}$ integrate uniquely to 1-dimensional foliations on $M$, which we call the \textbf{strong stable/unstable} foliations respectively, while $E^{s/u} \oplus \langle X \rangle$ integrate uniquely to 2-dimensional foliations on $M$, which we call the \textbf{weak stable/unstable} foliations respectively.
We say that $\phi$ has \textbf{orientable foliations} if $E^s$ and $E^u$ are orientable as line bundles on $M$.
Unless otherwise stated, all Anosov flows in this paper will have orientable foliations.

\begin{figure}
    \centering
    \begingroup%
  \makeatletter%
  \providecommand\color[2][]{%
    \errmessage{(Inkscape) Color is used for the text in Inkscape, but the package 'color.sty' is not loaded}%
    \renewcommand\color[2][]{}%
  }%
  \providecommand\transparent[1]{%
    \errmessage{(Inkscape) Transparency is used (non-zero) for the text in Inkscape, but the package 'transparent.sty' is not loaded}%
    \renewcommand\transparent[1]{}%
  }%
  \providecommand\rotatebox[2]{#2}%
  \newcommand*\fsize{\dimexpr\f@size pt\relax}%
  \newcommand*\lineheight[1]{\fontsize{\fsize}{#1\fsize}\selectfont}%
  \ifx\svgwidth\undefined%
    \setlength{\unitlength}{84.10125552bp}%
    \ifx\svgscale\undefined%
      \relax%
    \else%
      \setlength{\unitlength}{\unitlength * \real{\svgscale}}%
    \fi%
  \else%
    \setlength{\unitlength}{\svgwidth}%
  \fi%
  \global\let\svgwidth\undefined%
  \global\let\svgscale\undefined%
  \makeatother%
  \begin{picture}(1,0.9194332)%
    \lineheight{1}%
    \setlength\tabcolsep{0pt}%
    \put(0,0){\includegraphics[width=\unitlength,page=1]{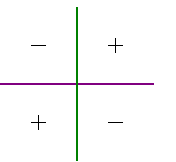}}%
    \put(0.932605,0.3903809){\color[rgb]{0.50196078,0,0.50196078}\makebox(0,0)[lt]{\lineheight{1.25}\smash{\begin{tabular}[t]{l}$E^u$\end{tabular}}}}%
    \put(0.49942132,0.88333393){\color[rgb]{0,0.50196078,0}\makebox(0,0)[lt]{\lineheight{1.25}\smash{\begin{tabular}[t]{l}$E^s$\end{tabular}}}}%
  \end{picture}%
\endgroup%

    \caption{The positive and negative quadrants in an Anosov flow. In this figure, the generating vector field is pointing out of the page, the stable direction is vertical and the unstable direction is horizontal.}
    \label{fig:posnegquadrants}
\end{figure}

At every point $p \in M$, the tangent space $T_p M$ has four complementary regions in the complement of $E^s \oplus \langle X \rangle |_p$ and $E^u \oplus \langle X \rangle |_p$.
We call these the \textbf{quadrants} at $p$.
We say that a quadrant $Q$ is \textbf{positive} if moving in a right-handed manner around $X|_p$, $Q$ is a region that starts with $E^u \oplus \langle X \rangle |_p$ and ends with $E^s \oplus \langle X \rangle |_p$.
Otherwise we say that $Q$ is \textbf{negative}.
See \Cref{fig:posnegquadrants}, where we draw the generating vector field $X$ as pointing out of the page, $E^s$ as the vertical direction and $E^u$ as the horizontal direction within the plane $TM/\langle X \rangle$. Under this convention, the positive/negative quadrants are the ones that contain vectors with positive/negative slope, respectively.

Two flows are \textbf{orbit equivalent} if there is a homeomorphism $h$ sending the flow lines of one to those of the other, preserving their orientations but not necessarily their parametrizations.
Furthermore, they are \textbf{isotopically equivalent} if they are defined on the same 3-manifold and $h$ is isotopic to identity.

We will consider two Anosov flows as equivalent if they are isotopically equivalent.
Under this convention, Anosov flows are structurally stable in the following sense, see \cite[Theorem 1]{Ano67} or \cite[Corollary 5.4.7]{FH19}.

\begin{thm} \label{thm:anosovflowstructuralstability}
Let $\phi$ be an Anosov flow with generating vector field $X$. There exists $\epsilon > 0$ such that every flow $\widetilde{\phi}$ generated by a vector field $\widetilde{X}$ that is $\epsilon$-close to $X$ in the $C^1$-topology is Anosov and isotopically equivalent to $\phi$.
\end{thm}

From \Cref{thm:anosovflowstructuralstability} we can deduce the following.

\begin{thm} \label{thm:anosovflow1parameterfamily}
Let $(\phi_s)_{s\in[0,1]}$ be a family of Anosov flows generated by a $C^1$-continuous family of vector fields, then $\phi_s$ are isotopically equivalent.
\end{thm}

Let $\gamma$ be a closed orbit of an Anosov flow $\phi$. 
Under our assumption that $\phi$ has orientable foliations, $E^s$ (or equivalently $E^u$) induces a framing of $\gamma$ in $M$, sometimes referred to as the \textbf{degeneracy slope} at $\gamma$. Unless otherwise stated, we will always put this framing on $\gamma$ in this context.

\begin{constr}[Goodman surgery \cite{Goo83}] \label{constr:Goodmansurgery}
Let $\gamma$ be a closed orbit of an Anosov flow $\phi$ on $M$.
For every $n \in \mathbb{Z}$, the manifold $M_{\frac{1}{n}}(\gamma)$ obtained by Dehn filling $M \backslash \nu(\gamma)$ along the slope $\frac{1}{n}$ admits an Anosov flow.
More specifically, this flow is generated by a vector field obtained by modifying $X$ in a small tubular neighborhood of $\gamma$; the exact modification is unimportant for this paper.
We refer to this construction as \textbf{$\frac{1}{n}$-Goodman surgery} along $\gamma$.

More generally, given a disjoint collection of closed orbits $\{\gamma_1,...,\gamma_k\}$ and integers $n_1,...,n_k \in \mathbb{Z}$, we can perform $\frac{1}{n_i}$-Goodman surgery along each $\gamma_i$ separately to get an Anosov flow on the surgered 3-manifold $M_{\frac{1}{n_1},...,\frac{1}{n_k}}(\gamma_1,...,\gamma_k)$.
\end{constr}

We conclude this subsection by recalling the class of skew Anosov flows.

Let $\phi$ be an Anosov flow on a closed 3-manifold $M$.
It can be shown that the quotient of the universal cover $\widetilde{M}$ under the lifted flow lines of $\phi$ is a topological plane $\mathcal{O}$, called the \textbf{orbit space} of $\phi$.
Moreover, the lifted leaves of the weak stable/unstable foliation of $\phi$ induce a 1-dimensional foliation of $\mathcal{O}$, which we refer to as the \textbf{stable/unstable foliation} on $\mathcal{O}$ respectively.
See \cite[Section 1.3]{BarMan25} for details.

We say that $\phi$ is a \textbf{positive skew} Anosov flow if there is an orientation-preserving identification $\mathcal{O} \cong \{(x,y) \in \mathbb{R}^2 \mid -y < x < -y+1\}$ with the stable/unstable foliations identified with the foliations by $\{x=x_0\}$-lines/$\{y=y_0\}$-lines, respectively.
We say that $\phi$ is a \textbf{negative skew} Anosov flow if $\mathcal{O} \cong \{(x,y) \in \mathbb{R}^2 \mid y < x < y+1\}$ with the stable/unstable foliations identified with the foliations by $\{x=x_0\}$-lines/$\{y=y_0\}$-lines, respectively.

\subsection{Contact structures} \label{subsec:contactstructure}

A \textbf{positive contact form} on $M$ is a 1-form $\alpha_+$ satisfying $\alpha_+ \wedge d\alpha_+ > 0$.
The \textbf{Reeb vector field} of $\alpha_+$ is defined to be the unique vector field $R_+$ satisfying $\iota_{R_+} \alpha_+ = 1$ and $\iota_{R_+} d\alpha_+ = 0$.
The \textbf{Reeb flow} of $\alpha_+$ is the flow generated by $R_+$.
Observe that this flow preserves $\alpha_+$ since $\mathcal{L}_{R_+} \alpha_+ = d \iota_{R_+} \alpha_+ + \iota_{R_+} d \alpha_+ = 0$.

Similarly, a \textbf{negative contact form} on $M$ is a 1-form $\alpha_-$ satisfying $\alpha_- \wedge d\alpha_- < 0$.
The \textbf{Reeb vector field} of $\alpha_-$ is defined to be the unique vector field $R_-$ satisfying $\iota_{R_-} \alpha_- = 1$ and $\iota_{R_-} d\alpha_- = 0$.
The \textbf{Reeb flow} of $\alpha_-$ is the flow generated by $R_-$, and it preserves $\alpha_-$.

A \textbf{positive/negative contact structure} is a cooriented plane field that can be expressed as the kernel of a positive/negative contact form, respectively.
The following proposition gives a convenient way to think about contact structures.

\begin{prop}[{\cite[Proposition 1.1.5]{ET98}}] \label{prop:contactstructurerudder}
Suppose $\xi$ is a cooriented plane field on $M$. Pick a nonvanishing vector field $V$ that is tangent to $\xi$ at every $p \in M$.
Then $\xi$ is a positive/negative contact structure if and only if at each $p \in M$, the family of planes $\left( V^{-t}_*(\xi|_{V^t(p)}) \right)_{t \in \mathbb{R}}$ rotates in a left/right-handed way around $V|_p$ as $t$ increases, respectively.
See \Cref{fig:contactstructurerudder}.
\end{prop}

In this subsection, we will collect some terminology and facts about contact structures. 
Some good references for this material are \cite{Hon04, ET}. %
While these references only deal with positive contact structures, statements can be carried over to negative contact structures by reversing the orientation of the 3-manifold.

Let $\xi$ be a contact structure on $M$.
An embedded 1-manifold is \textbf{Legendrian} if it is tangent to $\xi$ everywhere.
More generally, an embedded graph is \textbf{Legendrian} if its edges are tangent to $\xi$ everywhere.
An embedded 1-manifold is \textbf{transverse} if it is transverse to $\xi$ everywhere.
If a transverse 1-manifold is oriented, we can further specify whether it is \textbf{positive} or \textbf{negative}.

Let $c$ be a Legendrian curve. The contact structure $\xi$ induces a framing of $c$ in $M$, sometimes referred to as the \textbf{Legendrian framing} at $c$. Unless otherwise stated, we will always put this framing on such a curve $c$.

We have the following standard neighborhood theorem for Legendrian curves. %

\begin{prop}[{\cite[Theorem 1.2.4]{ET}}] \label{lemma:legendriancurvestdnbd}
Let $\xi$ be a contact structure, and let $c$ be a Legendrian curve.
Then there exists a tubular neighborhood $\nu \cong [-\epsilon,\epsilon]_x \times [-\epsilon,\epsilon]_y \times S^1_z$ of $c$ in which $c = \{(0,0,z) \mid z \in S^1\}$ and $\xi = \ker (dx \pm ydz)$, where the sign is positive/negative if $\xi$ is positive/negative, respectively.
\end{prop}

We move on from curves to surfaces.
Let $S$ be an embedded surface, possibly with boundary. The \textbf{characteristic foliation} on $S$, which we denote by $S_\xi$, is the singular oriented 1-dimensional foliation on $S$ tangent to the singular oriented line field $\xi \cap TS$. %

A vector field $V$ is \textbf{contact} if the flow it generates preserves $\xi$.
We say that $S$ is \textbf{convex} if there exists a contact vector field $V$ which is transverse to $S$.
In this case, we define the \textbf{dividing set} on $S$ to be the set $\Gamma = \{p \in S \mid V|_p \in \xi|_p\}$, and we define $S_{+/-}$ to be the subset of $S$ where $V$ is positively/negatively transverse to $\xi$, respectively.
It can be shown that $\Gamma$ is an embedded 1-manifold that separates $S$ into $S_+$ and $S_-$.

It can also be shown that the isotopy class of $\Gamma$ is independent of the transverse contact vector field $V$. However, this fact is less important to us since we will keep a choice of $V$ fixed in our applications.

Finally, it can be shown that $\Gamma$ divides the characteristic foliation $S_\xi$, in the following sense.

\begin{defn}[Divide] \label{defn:divide}
Let $S$ be a surface. Let $\Gamma$ be an embedded 1-manifold that separates $S$ into two subsurfaces, $S_+$ and $S_-$. Let $\mathcal{F}$ be a singular oriented 1-dimensional foliation on $S$.
We say that $\Gamma$ \textbf{divides} $\mathcal{F}$ if $\Gamma$ is transverse to $\mathcal{F}$ and
there exists an area form $\omega$ on $S$ and a vector field $Z$ directing $\mathcal{F}$ such that the divergence $\div_\omega Z$, defined by the equation $d \iota_Z \omega = (\div_\omega Z) \omega$,
is positive/negative on $S_{+/-}$ respectively.
\end{defn}

We say that an embedded graph $G$ on $S$ is \textbf{non-isolating} if 
it is transverse to $\Gamma$, in the sense that its vertices are disjoint from $\Gamma$ and its edges are transverse to $\Gamma$, and every component of $S \backslash G$ meets $\Gamma$ or a boundary component of $S$. %

The case that will matter for this paper is when $V$ is the Reeb vector field $R$ of a contact form $\alpha$ for $\xi$.
In this case, we have $\Gamma = \varnothing$.

For the rest of this subsection, we will review two results on convex surfaces and their characteristic foliations.
These are variations of what are commonly known as the Giroux flexibility theorem and the Legendrian realization principle. 
Since we could not find versions in the literature that match our setting, we take the time to provide some details of the proofs.
We also take the opportunity to state the results in greater generality than we need in this paper, in anticipation of future applications; see \Cref{subsec:relationwithbbmy} for more discussion.

Let $S$ be a convex surface with transverse contact vector field $V$. Then one can construct a thin collar neighborhood $V^{[-\epsilon,\epsilon]}(S) \cong S \times [-\epsilon,\epsilon]_t$ around $S$, under which the contact structure $\xi$ is invariant in the $t$-direction.
In other words, we can write $\xi = \ker(\beta + f dt)$, where $\beta$ is a 1-form and $f$ is a function on $S$.

Considering $S \times [-\epsilon,\epsilon]$ as a subset of an abstract copy of $S \times \mathbb{R}$, we can extend the contact structure $\xi = \ker(\beta+fdt)$ from $S \times [-\epsilon,\epsilon]$ to $S \times \mathbb{R}$ by $t$-invariance.
Here we emphasize that $S \times \mathbb{R}$ is \emph{not} considered as a submanifold of $M$.

\begin{thm} \label{thm:flexibility}
Let $S$ be a surface, possibly with boundary.
Consider $S \times \mathbb{R}_t$ with a $t$-invariant contact structure $\xi = \ker(\beta + f dt)$, where $\beta$ is a 1-form and $f$ is a function on $S$. 
Since $S \times \{0\}$ is transverse to the contact vector field $\frac{\partial}{\partial t}$, it is convex with dividing set $\Gamma = \{f=0\}$.

Then for every singular oriented foliation $\mathcal{F}$ on $S$ that is divided by $\Gamma$, there is an embedding $i\colon  S \hookrightarrow S \times \mathbb{R}_t$ transverse to the $t$-direction and isotopic to the inclusion $i_0 \colon S = S \times \{0\} \subset S \times \mathbb{R}$, such that $i(S)_\xi = \mathcal{F}$.

Moreover, if $\mathcal{F} = (S \times \{0\})_\xi$ on a subsurface $W$ that is disjoint from $\Gamma$ and the singularities of $\mathcal{F}$, and for which $\mathcal{F}$ is transverse to $\partial W$ and pointing into $W \cap S_{+/-}$ and out of $W \cap S_{-/+}$, depending on whether $\xi$ is positive/negative respectively, then we can arrange for $i$ to coincide with $i_0$ on $W$.
\end{thm}
\begin{proof}
The case when $W = \varnothing$ is contained in the proof of \cite[Theorem 8.3]{Hon04}. We briefly recap the argument so that we can indicate the necessary modifications when $W \neq \varnothing$. For concreteness suppose $\xi$ is a positive contact structure.

By assumption, there exists an area form $\omega$ on $S$ and a vector field $Z$ directing $\mathcal{F}$ such that $\div_\omega Z$ is positive/negative on $S_{+/-}$ respectively.
From this data, one can construct a function $f_1$ on $S$ that is positive/negative on $S_{+/-}$ respectively, such that $\xi_1 = \ker(\iota_Z \omega + f_1 dt)$ is a $t$-invariant contact structure on $S \times \mathbb{R}$. 

More specifically, one can follow the proof of \cite[Theorem 3.2.9]{ET}: 
Take a collar neighborhood $\nu \cong \Gamma \times [-\delta,\delta]_s$ of $\Gamma$, where $Z = \frac{\partial}{\partial s}$ and $\Gamma \times \{\mp \delta\} \subset S_\pm$. One can then set $f_1 = K_\pm$ on $S_\pm \backslash \nu$, for any constants $K_- < 0 < K_+$, and interpolate in $\nu$ by producing a smoothly varying solution of the inequality $f_1 \div_\omega Z - \frac{\partial f_1}{\partial s} > 0$ along each arc $x \times [-\delta,\delta]$ in $\nu$.
That this is possible is because $K_+ > K_-$, so we can take $\frac{\partial f_1}{\partial s} \leq 0$, with strict inequality away from the endpoints, while $f_1 \div_\omega Z \geq 0$, with strict inequality away from $\Gamma = \{s=0\}$. Or see \cite[Theorem 3.2.9]{ET} for an explicit solution.

Then the characteristic foliation of $S \times \{0\}$ with respect to $\xi_1$ is $\mathcal{F}$.
The strategy is to find a $t$-invariant isotopy that takes $\xi_1$ to $\xi$, then take the image of $S \times \{0\}$ under this isotopy.

Up to such an isotopy, we can first assume that $\xi_1 = \xi$ in a collar neighborhood of $\Gamma$. 
Afterwards, we can choose contact forms $\alpha = \beta \pm dt$ and $\alpha_1 = \beta_1 \pm dt$ for $\xi$ and $\xi_1$ respectively on $S_\pm$ away from a smaller collar neighborhood of $\Gamma$.
Then $(\alpha_s = (1-s) \beta + s \beta_1 \pm dt)_{s \in [0,1]}$ is a 1-parameter family of contact forms.

Finally, we use Moser's trick: For each $s \in [0,1]$, there is a $t$-invariant vector field $X_s$ satisfying $\mathcal{L}_{X_s} \alpha_s = \frac{d}{ds} \alpha_s$. We integrate this vector field to find the isotopy. %
Since $X_s$ is $t$-invariant, the same is true for the isotopy. 
This completes the proof when $W = \varnothing$.

Now if $\mathcal{F} = (S \times \{0\})_\xi$ on a subsurface $W$, then up to rescaling, we can assume that $\beta = \iota_Z \omega$ on $W$. 
We use the fact that $W$ is disjoint from the singularities of $\mathcal{F}$ here to ensure that the rescaling factor is smooth.
We can then choose a function $f_1$ on $S$ such that 
\begin{itemize}
    \item $f_1$ agrees with $f$ on $W$, 
    \item $f_1$ is positive/negative on $S_{+/-}$ respectively, and
    \item $\xi_1 = \ker(\iota_Z \omega + f_1 dt)$ is a $t$-invariant contact structure on $S \times \mathbb{R}$.
\end{itemize}

More specifically, we take a collar neighborhood $\nu \cong (\Gamma \cup \partial W) \times [-\delta,\delta]_s$ of $\Gamma \cup \partial W$ where $Z = \frac{\partial}{\partial s}$, $\Gamma \times \{\mp \delta\} \subset S_\pm$, and $\partial (W \cap S_\pm) \times \{\pm \delta\} \subset W \cap S_\pm$.
We use the facts that $W$ is disjoint from $\Gamma$ and that $\mathcal{F}$ is transverse to $\partial W$ here.
One can then set $f_1 = f$ on $W$ and $f_1 = K_\pm$ on $S_\pm \backslash (W \cup \nu)$, and interpolate in $\nu$ by producing a smoothly varying solution of the inequality $f_1 \div_\omega Z - \frac{\partial f_1}{\partial s} > 0$ along each arc $x \times [-\delta,\delta]$ in $\nu$.
One can use the same argument for this being possible as in the $W = \varnothing$ case provided that $|K_\pm| > \max_{W \cap S_\pm} |f|$.
We use the fact that $\mathcal{F}$ is pointing into $W \cap S_+$ and out of $W \cap S_-$ here.

We then repeat the argument above in order to find a $t$-invariant isotopy that takes $\xi_1$ to $\xi$, then take the image of $S \times \{0\}$ under this isotopy.
The key point here is that since $\xi = \xi_1$ on $W \times \mathbb{R}$, the contact forms $\alpha = \beta \pm dt$ and $\alpha_1 = \beta_1 \pm dt$ satisfy $\beta = \beta_1$ on $W$.
Thus the 1-parameter family $\alpha_s$ stays constant on $W \times \mathbb{R}$, so the same is true for the isotopy. %
\end{proof}

\begin{thm} \label{thm:legendrianrealization}
Let $S$ be a surface, possibly with boundary.
Consider $S \times \mathbb{R}_t$ with a $t$-invariant contact structure $\xi = \ker(\beta + f dt)$, where $\beta$ is a 1-form and $f$ is a function on $S$. 
Since $S \times \{0\}$ is transverse to the contact vector field $\frac{\partial}{\partial t}$, it is convex with dividing set $\Gamma = \{f=0\}$.

Let $G$ be a non-isolating graph on $S$. 
Then there is an embedding $i \colon S \hookrightarrow S \times \mathbb{R}_t$ transverse to the $t$-direction and isotopic to the inclusion $i_0 \colon S = S \times \{0\} \subset S \times \mathbb{R}$, such that $i(G)$ is Legendrian.

Moreover, if we further assume that every boundary component $b$ of $S \times \{0\}$ is:
\begin{itemize}
    \item disjoint from $\Gamma$;
    \item positively transverse if $b$ lies in $S_{+/-}$ and  negatively transverse if $b$ lies in $S_{-/+}$, depending on whether $\xi$ is positive/negative respectively; and
    \item $G \times \{0\}$ is already Legendrian near $\partial S \times \{0\}$,
\end{itemize}
then we can arrange for $i$ to coincide with $i_0$ in a collar neighborhood of $\partial S$.
\end{thm}
\begin{proof}
This is essentially contained in the proof of \cite[Theorem 3.4.5]{ET}. We briefly recap the argument for completeness. 
For concreteness suppose $\xi$ is a positive contact structure.

The strategy is to define a foliation $\mathcal{F}$ on $S$ that contains $G$ and which is divided by $\Gamma$, then apply \Cref{thm:flexibility}.
In the case when the boundary components of $S \times \{0\}$ that lie in $S_{+/-}$ are positively/negatively transverse, the original characteristic foliation $(S \times \{0\})_\xi$ is a nonsingular foliation by non-separating arcs in a collar neighborhood of such boundary components, with the arcs pointing into/out of the collar neighborhood. 
The foliation $\mathcal{F}$ will be defined to agree with $(S \times \{0\})_\xi$ in this neighborhood, so that, if $\partial S \cap \Gamma = \varnothing$, we can ensure the embedding in \Cref{thm:flexibility} is the inclusion $p \mapsto (p,0)$ in this neighborhood.

Up to subdividing the edges of $G$, we can assume that every component of $G \cap S_\pm$ contains at least one vertex.
We can then classify each vertex as positive or negative depending on whether it lies in $S_+$ or $S_-$, and declare each positive/negative vertex as a positive/negative elliptic singularity of $\mathcal{F}$ respectively.
We then declare the midpoint between any two adjacent singularities of the same sign in $G$ to be a hyperbolic singularity, and declare the remaining complementary intervals in $G$ to be nonsingular leaves of $\mathcal{F}$.

As explained above, in the setting of the `moreover' statement, we further declare $\mathcal{F}$ to agree with $(S \times \{0\})_\xi$ in a collar neighborhood of each boundary component, so that in particular $\mathcal{F}$ is pointing into/out of collar neighborhoods of the boundary components that lie in $S_{+/-}$, respectively.

Finally, we extend $\mathcal{F}$ into each complementary region of $G \cup \Gamma$ as a gradient vector field, as in \cite[Theorem 3.4.5]{ET}.
\end{proof}

We emphasize that the embeddings provided by \Cref{thm:flexibility} and \Cref{thm:legendrianrealization} map into $S \times \mathbb{R}$ instead of the neighborhood $S \times [-\epsilon,\epsilon]$ of $S$ in $M$.
In other words, these embeddings `may not fit' inside $M$.
We will `make room' for these embeddings to fit by adding padding in \Cref{subsec:factoroutwindzerocurve}.

\begin{rmk} \label{rmk:legendrianrealization}
The reader who is more familiar with contact geometry may be confused here:
It is true that there is a contact embedding of $S \times \mathbb{R}$ into $S \times [-\epsilon,\epsilon]$, i.e., an embedding that preserves the contact structures. However, in general this embedding distorts the $t$-direction.
Since it is important for us to keep track of the \emph{given} contact vector field $V$, such an embedding is unhelpful to us.
\end{rmk}

\subsection{Bicontact packages} \label{subsec:bicontactpackage}

For the purposes of this paper, we define a \textbf{bicontact structure} to be a pair $(\xi_-,\xi_+)$ consisting of a negative contact structure $\xi_-$ and a positive contact structure $\xi_+$ that are transverse to each other at every point.
We say that $(\xi_-,\xi_+)$ \textbf{supports} a flow generated by a nonvanishing vector field $X$ if 
$X$ is contained in and points in the positive direction of $\xi_- \cap \xi_+$ at every point.

Building on this notion, we define a \textbf{bicontact package} to be a pair $(\alpha_-,\xi_+)$ consisting of a negative contact \emph{form} $\alpha_-$ and a positive contact \emph{structure} $\xi_+$, such that $\xi^+$ and $\ker \alpha^-$ are transverse.
We say that $(\alpha_-,\xi_+)$ \textbf{supports} a flow $\phi$ if its underlying bicontact structure $(\xi_- = \ker \alpha_-,\xi_+)$ supports $\phi$.
We say that $(\alpha_-,\xi_+)$ is \textbf{strongly adapted} if the Reeb vector field $R_-$ of $\alpha_-$ is tangent to $\xi_+$ everywhere.

The significance of this definition lies in the following result of Hozoori.

\begin{thm}[Hozoori {\cite[Theorem A]{Hoz25}}] \label{thm:hozoori}
A flow $\phi$ on a closed oriented 3-manifold is Anosov with orientable foliations if and only if it is supported by a strongly adapted bicontact package $(\alpha_-,\xi_+)$.
\end{thm}

In the context of \Cref{thm:hozoori}, one can deduce from \Cref{prop:contactstructurerudder} that the contact structures $\xi_- = \ker \alpha_-$ and $\xi_+$ must stay in the positive and negative quadrants at every point. 
In particular, a curve $c$ is an orbit of $\phi$ if and only if it is Legendrian with respect to both $\xi_-$ and $\xi_+$, and in this case the degeneracy slope agrees with both of the Legendrian framings at $c$.

We have the following stability results from \Cref{thm:anosovflowstructuralstability}.

\begin{prop} \label{prop:negcontactformstructuralstability}
Let $(\alpha_-,\xi_+)$ be a strongly adapted bicontact package supporting an Anosov flow $\phi$. There exists $\epsilon > 0$ such that for every 1-form $\widetilde{\alpha}_-$ that is $\epsilon$-close to $\alpha_-$ in the $C^2$-topology, there exists $\widetilde{\xi}_+$ such that $(\widetilde{\alpha}_-,\widetilde{\xi}_+)$ is a strongly adapted bicontact package and supports an Anosov flow $\widetilde{\phi}$ isotopically equivalent to $\phi$.
\end{prop}
\begin{proof}
Choose a positive contact form $\alpha_+$ for $\xi_+$.
Any $\widetilde{\alpha}_-$ that is $C^2$-close enough to $\alpha_-$ is a negative contact form whose Reeb vector field $\widetilde{R}_-$ is $C^1$-close to $R_-$. We can thus define a form $\widetilde{\alpha}_+ = \alpha_+ - \alpha_+(\widetilde{R}_-) \widetilde{\alpha}_- = \alpha_+ - \alpha_+(\widetilde{R}_- - R_-) \widetilde{\alpha}_-$ which is $C^1$-close to $\alpha_+$, thus a positive contact form.
By construction, $(\widetilde{\alpha}_-,\widetilde{\xi}_+ = \ker \widetilde{\alpha}_+)$ is a strongly adapted bicontact package, thus by \Cref{thm:hozoori} supports an Anosov flow $\widetilde{\phi}$.

Since $\widetilde{\alpha}_{\pm}$ are $C^1$-close to $\alpha_{\pm}$, the generating vector field of $\widetilde{\phi}$ is $C^1$-close to that of $\phi$. Thus by \Cref{thm:anosovflowstructuralstability}, $\widetilde{\phi}$ is isotopically equivalent to $\phi$.
\end{proof}

\begin{defn}[Homotopy of bicontact packages] \label{defn:bicontactpackagehomotopy}
We say that two strongly adapted bicontact packages $((\alpha_-)_0,(\xi_+)_0)$ and $((\alpha_-)_1,(\xi_+)_1)$ are \textbf{homotopic} if they can be extended into a $C^2$-continuous family of strongly adapted bicontact packages $((\alpha_-)_s,(\xi_+)_s)_{s \in [0,1]}$.
\end{defn}

\begin{prop} \label{prop:bicontactpackage1parameterfamily}

Let $((\alpha_-)_0,(\xi_+)_0)$ and $((\alpha_-)_1,(\xi_+)_1)$ be two strongly adapted bicontact packages that are homotopic. Then the Anosov flows that they support are isotopically equivalent.
\end{prop}
\begin{proof}
This follows immediately from \Cref{thm:hozoori} and \Cref{thm:anosovflow1parameterfamily}.
\end{proof}

\begin{prop} \label{prop:poscontactstructureflowisotopy} %
Let $(\alpha_-,\xi_+)$ be a strongly adapted bicontact package. Let $h$ be a diffeomorphism obtained via an isotopy along the flow lines of $R_-$, i.e., $h(p) = R_-^{t(p)}(p)$ for some $t\colon M \to \mathbb{R}$, and let $\widetilde{\xi}_+ = h^*(\xi_+)$. Then $(\alpha_-,\widetilde{\xi}_+)$ is a strongly adapted bicontact package that is homotopic to $(\alpha_-,\xi_+)$. 
\end{prop}
\begin{proof}
For $s \in [0,1]$, we define $h_s(p) = R_-^{st(p)}(p)$, and $(\xi_+)_s = h_s^*(\xi_+)$.
Then $(\alpha_-,(\xi_+)_s)$ is a smooth family of strongly adapted bicontact packages. 
\end{proof}

We also have the following surgical operation, due to Salmoiraghi \cite{Sal25}.

\begin{constr}[Legendrian-transverse surgery] \label{constr:legendriantransversesurgery}
Let $(\alpha_-,\xi_+)$ be a strongly adapted bicontact package. 
We say that a curve $c$ is \textbf{$\xi_-$-Legendrian-$\xi_+$-transverse} if it is Legendrian with respect to $\xi_- = \ker \alpha_-$ and transverse with respect to $\xi_+$.
In this case, recall that $c$ naturally carries a Legendrian framing from $\xi_-$.
For every $n \in \mathbb{Z}_+$, the manifold $M_{\frac{1}{n}}(c)$ admits a strongly adapted bicontact package, obtained by modifying $(\alpha_-,\xi_+)$ in an arbitrarily small tubular neighborhood of $c$.

Symmetrically, a curve $c$ is \textbf{$\xi_+$-Legendrian-$\xi_-$-transverse} if it is Legendrian with respect to $\xi_+$ and transverse with respect to $\xi_-$.
In this case, for every $n \in \mathbb{Z}_-$, the manifold $M_{\frac{1}{n}}(c)$ admits a strongly adapted bicontact package obtained by modifying $(\alpha_-,\xi_+)$ in an arbitrarily small tubular neighborhood of $c$.

We refer to this construction as \textbf{$\frac{1}{n}$-Legendrian-transverse surgery} along $c$.

Goodman surgery (\Cref{constr:Goodmansurgery}) is a specific case of Legendrian-transverse surgery in the following sense:
Given a closed orbit $\gamma$ of an Anosov flow supported by a strongly adapted bicontact package $(\alpha_-,\xi_+)$, one can find a $\xi_-$-Legendrian-$\xi_+$-transverse curve $c$ arbitrarily close to and parallel to $\gamma$, then for every $n \in \mathbb{Z}_+$, $\frac{1}{n}$-Legendrian-transverse surgery along $c$ is $\frac{1}{n}$-Goodman surgery along $\gamma$.
Similarly, one can find a $\xi_+$-Legendrian-$\xi_-$-transverse curve $c$ arbitrarily close to and parallel to $\gamma$, then for every $n \in \mathbb{Z}_-$, $\frac{1}{n}$-Legendrian-transverse surgery along $c$ is $\frac{1}{n}$-Goodman surgery along $\gamma$.
We refer to \cite[Section 7]{Sal25} for details. %
\end{constr}

 Now that Legendrian-transverse surgery is defined, we explain how it can be used to turn any Anosov flow into a transitive one: %

\begin{prop}\label{prop:legendrian_transverse_gives_transitive}
    Let $\phi$ be an Anosov flow on $M$ with orientable foliations. Then $\phi$ can be transformed into a transitive Anosov flow by Legendrian-transverse surgery.
\end{prop}

\begin{proof}
    Let $(\alpha_-,\xi_+)$ be a strongly adapted bicontact package for $\phi$. As transitive flows are $C^1$-generic among Reeb flows (see \cite[Theorem C.4]{CDHR24}), we may assume that $R_-$ is transitive. By Pugh's $C^1$-closing lemma (or Irie's smooth closing lemma for Reeb flows in $3$-manifolds \cite[Lemma 3.1]{Iri15}), we can further assume that for any fixed $\epsilon>0$, there exists an $\epsilon$-dense periodic orbit of $R_-$.%

In particular, we may assume that $R_-$ admits a periodic orbit $\alpha$ that intersects every embedded incompressible torus in $M$. Now let $\psi$ be the flow on $N$ obtained by doing a large enough Legendrian-transverse surgery on $\alpha$. If $\psi$ is not transitive, then $N$ must admit at least one torus that is transverse to the flow $\psi$ (see \cite{Bru93}). %

Let $T$ be such a torus. Since we chose a large enough slope for the surgery, $T$ does not intersect the surgery curve $\alpha$. In particular, $T$ must be compressible in $M$. Moreover, since Legendrian-transverse surgery only changes the vector field locally around the surgery curve $\alpha$, $T$ must be transverse also to the flow $\phi$. But a torus transverse to an Anosov flow is automatically incompressible (see e.g., \cite[Lemma 6.4.2]{BarMan25}), so $T$ must be incompressible in $M$, a contradiction.
\end{proof}

Note that thanks to \Cref{prop:legendrian_transverse_gives_transitive}, \Cref{thm:LT_and_bicontact_gives_suspension} will follow from \Cref{cor:flowreducetosusp}.

\section{Bicontact padding surgery} \label{sec:paddingsurgery}

\subsection{Birkhoff section} \label{subsec:birkhoffsection}

Let $\phi$ be a flow on a closed oriented 3-manifold $M$ generated by a vector field $X$. 

\begin{defn} \label{defn:birkhoffsection}
A \textbf{Birkhoff section} to $\phi$ is a compact cooriented immersed surface $S \looparrowright M$ satisfying:
\begin{enumerate}
    \item the interior of $S$ is embedded and positively transverse to $X$,
    \item each boundary component of $S$ is tangent to $X$, equivalently, each boundary component of $S$ lies along a closed orbit of $\phi$, and
    \item every forward and backward trajectory meets $S$, i.e., for every $p \in M$, there exist $T_- < 0 < T_+$ such that $\phi^{T_-}(p), \phi^{T_+}(p) \in S$. %
\end{enumerate}
\end{defn}

Suppose $b$ is a boundary component of $S$ lying along a closed orbit $\gamma$ of $\phi$. We say that $b$ is a \textbf{positive} boundary component if its orientation as a boundary component of $S$ agrees with the orientation pulled back from $\gamma$ as a closed orbit of $\phi$. In this case we also say that $\gamma$ is a \textbf{positive} boundary orbit. 
If the orientations disagree, we say that $b$ is a \textbf{negative} boundary component and $\gamma$ is a \textbf{negative} boundary orbit.

By (3), the drilled manifold $M \backslash \partial S$ is a mapping torus $\intr(S) \times [0,1]/(x,1) \sim (f(x),0)$. 
We refer to $f$ as the \textbf{monodromy} of $S$. This is a well-defined element of the mapping class group of $\intr(S)$.

The following construction is helpful for analyzing Birkhoff sections near their boundary orbits.

\begin{constr}[Blow-up] \label{constr:blowup}
Suppose $\gamma$ is a closed orbit of $\phi$. The \textbf{blow-up} of $M$ along $\gamma$ is the 3-manifold with boundary $M_\gamma$ obtained by replacing $\gamma$ with the quotient bundle $((TM|_\gamma/T\gamma) \backslash \{0\})/\mathbb{R}_+$. The \textbf{blow-up} of $\phi$ is the flow $\phi_\gamma$ on $M_\gamma$ that agrees with $\phi$ away from $\partial M_\gamma = ((TM|_\gamma/T\gamma) \backslash \{0\})/\mathbb{R}_+$ and equals the linearization of $\phi$ along $\partial M_\gamma$.

More concretely, pick a tubular neighborhood $\nu$ of $\gamma$ and pick coordinates $\nu \cong D^2 \times S^1_z$ where $\gamma = \{0\} \times S^1$.
Then the blown-up manifold $M_\gamma$ is obtained by removing $\gamma$ from $M$ and gluing in $[0,1]_r \times S^1_\theta \times S^1_z$ along $\nu \backslash \gamma \cong (0,1]_r \times S^1_\theta \times S^1_z$.
The blown-up flow $\phi_\gamma$ agrees with $\phi$ away from the boundary torus $\{r=0\}$ and is defined by $\phi^t_\gamma(0,\theta,z) = (0,\theta(t),z+t)$, where $d\phi^t|_{(0,0,z)}([\cos \theta \frac{\partial}{\partial x} + \sin \theta \frac{\partial}{\partial y}]) = [\cos \theta(t) \frac{\partial}{\partial x} + \sin \theta(t) \frac{\partial}{\partial y}]$, on $\{r=0\}$.

In particular, note that the boundary torus is foliated by flow lines of $\phi_\gamma$.
There is a natural fibration $\partial M_\gamma = ((TM|_\gamma/T\gamma) \backslash \{0\})/\mathbb{R}_+ \to \gamma$ with the fibers being the meridional curves $\{r=0, z=z_0\}$. The flow lines of $\phi_\gamma$ on $\partial M_\gamma$ are transverse to these fibers.
This allows us to talk about the \textbf{slope} $f$ of the foliation by flow lines on $\partial M_\gamma$:
If we identify the slopes with $\mathbb{R} \cup \{\infty\}$ so that $\frac{\partial}{\partial z}$ is $0$ and $\frac{\partial}{\partial \theta}$ is $\infty$, then $f$ is $\frac{1}{2\pi}$ times the translation number of a lift of the first return map on a fiber to $\mathbb{R}$. %
In particular, the slope $f$ is a rational number $\frac{p}{q}$ if and only if there is a closed orbit of $\phi_\gamma$ on $\partial M_\gamma$ of slope $\frac{p}{q}$.
\end{constr}

Suppose $S$ is a Birkhoff section to $\phi$.
If $\gamma$ is a boundary orbit of $S$, then there is an induced immersion of $S$ into $M_\gamma$: The immersion agrees with that into $M$ away from $\partial M_\gamma$, and is induced by the tangent plane of $S$ along any boundary component that lies along $\gamma$. To avoid confusion, we will often denote this copy of $S$ by $S_\gamma$.
Since $S$ is immersed into $M$, the boundary components of $S_\gamma$ that lie on $\partial M_\gamma$ must be transverse to fibers of the fibration $\partial M_\gamma \to \gamma$. 
In particular, the slope $\sigma$ of $S_\gamma \cap \partial M_\gamma$ cannot coincide with the meridian $\mu$.

Following the terminology of \cite{CDHR24}, we say that $S$ is \textbf{$\partial$-strong} at $\gamma$ if $S_\gamma \cap \partial M_\gamma$ is embedded and a section of $\phi_\gamma$ on $\partial M_\gamma$.
In this case, the slope $\sigma$ of $S_\gamma \cap \partial M_\gamma$ cannot coincide with the slope $f$ of the foliation by flow lines of $\phi_\gamma$.
Moreover, $\gamma$ is a positive boundary orbit if and only if $\sigma \in (f,\mu)$.

Conversely, given a rational slope $\sigma \neq f, \mu$, one can construct a section $c$ to $\phi_\gamma$ on $\partial M_\gamma$ of slope $\sigma$ that is transverse to the fibers of $\partial M_\gamma \to \gamma$. If $S_\gamma$ is a section to $\phi_\gamma$ on $M_\gamma$ extending $c$, then we can take its image under the collapsing map $M_\gamma \to M$ to get a Birkhoff section to $\phi$. %

We say that a Birkhoff section is \textbf{$\partial$-strong} if it is $\partial$-strong near all of its boundary orbits.
Unless otherwise stated, all Birkhoff sections in this paper will be $\partial$-strong.

It is instructive to analyze how a Birkhoff section $S$ can fail to be $\partial$-strong.
Suppose we blow up along a boundary orbit $\gamma$.
Since $S$ is embedded in its interior, the same is true for $S_\gamma$, so the only way that these boundary components can fail to be embedded is that they touch, but not cross, tangentially. 
Moreover, since $S$ is positively transverse to $\phi$ in its interior, the same is true for $S_\gamma$, so the boundary components are either positively transverse or tangent to $\phi_\gamma$ at each point on $\partial M_\gamma$.

\begin{figure}
    \centering
    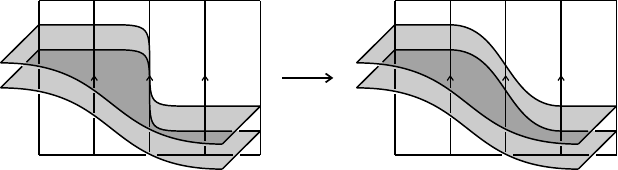
    \caption{Most of the time, it is possible to perturb a Birkhoff section along flow lines to become $\partial$-strong.}
    \label{fig:perturbforboundarystrong}
\end{figure}

If a boundary component is positively transverse to $\phi_\gamma$ at some point, then we can perturb $S_\gamma$ to make the whole boundary component positively transverse to $\phi_\gamma$, while maintaining transversality to the fibration $\partial M_\gamma \to \gamma$. %
If each boundary component on $\partial M_\gamma$ is positively transverse to $\phi_\gamma$ at some point, then we can do these perturbations in a way so that the boundary components become embedded and disjoint from each other.
See \Cref{fig:perturbforboundarystrong}.

We can then take the image of the perturbed $S_\gamma$ under the collapsing map $M_\gamma \to M$. The result is a copy of $S$ in $M$ that is perturbed along flow lines and now $\partial$-strong at $\gamma$.
This discussion shows that being $\partial$-strong is a rather mild hypothesis, for example we can deduce the following proposition.

\begin{prop} \label{prop:anosovboundarystrong}
Let $S$ be a Birkhoff section to an Anosov flow $\phi$. Up to isotoping $S$ along the flow lines of $\phi$, we can arrange for $S$ to be $\partial$-strong.
\end{prop}
\begin{proof}
If $\phi$ is Anosov, then on each blown-up boundary torus $\partial M_\gamma$, the blown-up flow $\phi_\gamma$ contains a finite number of attracting closed orbits and an equal number of repelling closed orbits. These have degeneracy slope and are arranged alternatively.

Meanwhile, the boundary components of $S_\gamma$ on $\partial M_\gamma$ must have a slope different from the degeneracy slope, otherwise $S$ would miss some (infinite) flow line that limits onto $\gamma$ in forward or backward time, contradicting the fact that $S$ is a Birkhoff section.
So each boundary component is positively transverse to $\phi_\gamma$ at some point, and we can perturb as discussed above.
\end{proof}

\subsection{Standard Reeb orbit} \label{subsec:stdorbit}

Consider a thickened torus $\mathbf{T} = T^2 \times I_r$, where $I$ is an interval. 
We say that a 1-form $\alpha$ on $\mathbf{T}$ is \textbf{standard} if under an affine coordinate system $T^2 \cong S^1_x \times S^1_y$, we have $\alpha = A(r) dx + B(r) dy$ for some functions $A,B$ depending only on $r$. %
More invariantly, one can think of $\alpha$ as the parametrized path $\widetilde{\alpha}(r) = (B(r), -A(r))$ on the plane $\widetilde{T}^2 = \mathbb{R}^2$. 
A short computation shows that $\alpha$ is negative contact if and only if $A(r)B'(r) - A'(r)B(r) > 0$, i.e., if $\widetilde{\alpha}$ winds around the origin counterclockwise as $r$ increases. Compare with \Cref{prop:contactstructurerudder}.
In this case, we write $\alpha = \alpha_-$ and compute that its Reeb vector field is $R_- = \frac{1}{A(r)B'(r) - A'(r)B(r)} (B'(r) \frac{\partial}{\partial x} - A'(r) \frac{\partial}{\partial y})$.
In particular, each level torus $T^2 \times \{r_0\}$ is foliated by orbits of $R_-$.

In this setting, there are two families of slopes at play:
\begin{itemize}
    \item The slope $\chi(r)$ of the characteristic foliation of each $T^2 \times \{r\}$ with respect to $\ker \alpha_-$. 
    (Here $\chi$ stands for `characteristic'.)
    If we think of $\alpha_-$ as its corresponding parametrized path $\widetilde{\alpha_-}$, this is also the slope of the line from the origin to $\widetilde{\alpha_-}(r)$. As $r$ increases, this slope has to decrease strictly.
    (Recall that the set of slopes is oriented \emph{clockwise}.)
    \item The slope $f(r)$ determined by the orbits of the Reeb flow $R_-$ on each $T^2 \times \{r\}$. 
    (Here $f$ stands for `flow'.)
    If we think of $\alpha_-$ as its corresponding parametrized path $\widetilde{\alpha_-}$, then $f(r)$ is also the slope of 
    the derivative $\widetilde{\alpha_-}'(r)$. 
    As $r$ increases, this slope may increase or decrease, but can never coincide with $\chi(r)$.
\end{itemize}

\begin{figure}
    \centering
    \begingroup%
  \makeatletter%
  \providecommand\color[2][]{%
    \errmessage{(Inkscape) Color is used for the text in Inkscape, but the package 'color.sty' is not loaded}%
    \renewcommand\color[2][]{}%
  }%
  \providecommand\transparent[1]{%
    \errmessage{(Inkscape) Transparency is used (non-zero) for the text in Inkscape, but the package 'transparent.sty' is not loaded}%
    \renewcommand\transparent[1]{}%
  }%
  \providecommand\rotatebox[2]{#2}%
  \newcommand*\fsize{\dimexpr\f@size pt\relax}%
  \newcommand*\lineheight[1]{\fontsize{\fsize}{#1\fsize}\selectfont}%
  \ifx\svgwidth\undefined%
    \setlength{\unitlength}{165.8134929bp}%
    \ifx\svgscale\undefined%
      \relax%
    \else%
      \setlength{\unitlength}{\unitlength * \real{\svgscale}}%
    \fi%
  \else%
    \setlength{\unitlength}{\svgwidth}%
  \fi%
  \global\let\svgwidth\undefined%
  \global\let\svgscale\undefined%
  \makeatother%
  \begin{picture}(1,0.95480298)%
    \lineheight{1}%
    \setlength\tabcolsep{0pt}%
    \put(0,0){\includegraphics[width=\unitlength,page=1]{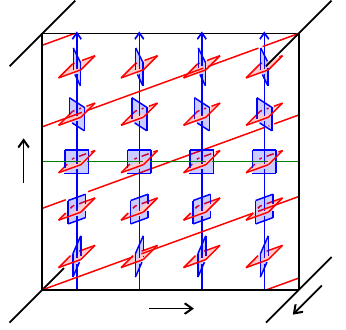}}%
    \put(0.91239407,0.05114539){\color[rgb]{0,0,0}\makebox(0,0)[lt]{\lineheight{1.25}\smash{\begin{tabular}[t]{l}$r$\end{tabular}}}}%
    \put(0.46352332,0.00209893){\color[rgb]{0,0,0}\makebox(0,0)[lt]{\lineheight{1.25}\smash{\begin{tabular}[t]{l}$x$\end{tabular}}}}%
    \put(-0.00224352,0.47401281){\color[rgb]{0,0,0}\makebox(0,0)[lt]{\lineheight{1.25}\smash{\begin{tabular}[t]{l}$y$\end{tabular}}}}%
    \put(0.37583356,0.90456933){\color[rgb]{0,0,1}\makebox(0,0)[lt]{\lineheight{1.25}\smash{\begin{tabular}[t]{l}$f$\end{tabular}}}}%
    \put(0.89810222,0.4678689){\color[rgb]{0,0.50196078,0}\makebox(0,0)[lt]{\lineheight{1.25}\smash{\begin{tabular}[t]{l}$\ell$\end{tabular}}}}%
    \put(0.89810222,0.61887206){\color[rgb]{1,0,0}\makebox(0,0)[lt]{\lineheight{1.25}\smash{\begin{tabular}[t]{l}$\chi$\end{tabular}}}}%
  \end{picture}%
\endgroup%

    \caption{A standard bicontact package on $\mathbf{T}$.}
    \label{fig:stdtorusshell}
\end{figure}
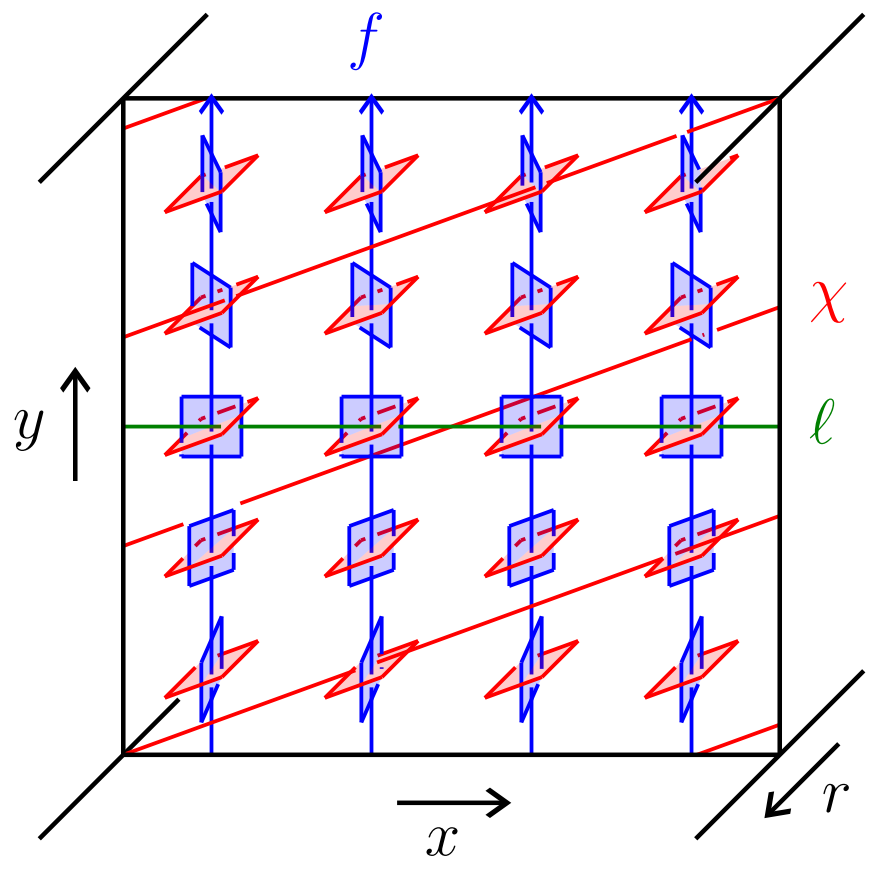

See \Cref{fig:stdtorusshell} for an example, ignoring $\xi_+$ in the figure for now.

We say that $\alpha_-$ is \textbf{rational} if the flow slope $f(r)$ is a constant rational slope.
In this case, we can choose coordinates $T^2 \cong S^1_x \times S^1_y$, so that $f(r)$ is the $y$-direction. This forces $\alpha_- = A(r) dx + B dy$ with $B$ being constant, and $R_- = \frac{1}{B} \frac{\partial}{\partial y}$.

We say that a strongly adapted bicontact package $(\alpha_-,\xi_+)$ on $\mathbf{T}$ is \textbf{standard} if the negative contact form $\alpha_-$ is standard.
In this setting, we add to our list of slopes:
\begin{itemize}
    \item For each $r$, consider the locus of positive tangency between $T^2 \times \{r\}$ and $\xi_+$. By \Cref{prop:contactstructurerudder}, this set is a multicurve.
    We denote the isotopy class of this multicurve by $l$ and refer to it as the \textbf{Legendrian class}. 
    (Here $l$ stands for `Legendrian'.)
    As the notation suggests, note that $l$ does \emph{not} depend on $r$.
    Also, $f(r)$ never coincides with $l$. 
    In fact, we can naturally coorient $l$ using the direction of the Reeb flow, thus the Legendrian class can be thought of as a homology class, justifying its terminology.
    \item For certain applications, we will also keep track of a slope $\sigma$, which will be the slope induced by a Birkhoff section to the Reeb flow generated by $R_-$.
    (Here $\sigma$ stands for `section'.)
    As such, $f(r)$ will also never coincide with $\sigma$.
\end{itemize}
See \Cref{fig:stdtorusshell} for an example.

Now suppose further that $\alpha_-$ is rational.
As mentioned above, in this case we can choose coordinates $T^2 \cong S^1_x \times S^1_y$, under which $\alpha_- = A(r) dx + B dy$ and $R_- = \frac{1}{B} \frac{\partial}{\partial y}$.
We say that $\xi_+$ is \textbf{linear} if $\xi_+ = \ker(\cos(px+qy)dr - \sin(px+qy)dx)$ for some $p \in \mathbb{Z}, q \in \mathbb{Z}_+$.
Notice that $\xi_+$ being linear is \emph{not} a requirement for $(\alpha_-,\xi_+)$ to be standard. This is to allow for sufficient flexibility in applications.

We have the following lemmas for manipulating standard strongly adapted bicontact packages on thickened tori.

\begin{lemma}[Interpolating the positive contact structure] \label{lemma:fillstdorbitinterpolatepos}
Let $T^2 \times [r_0,r_3]$ be a thickened torus with a standard rational negative contact form $\alpha_-$. Choose coordinates $T^2 \cong S^1_x \times S^1_y$ under which $\alpha_- = A(r) dx + B dy$ and $R_- = \frac{1}{B} \frac{\partial}{\partial y}$.
Suppose $\alpha_-$ is part of a strongly adapted bicontact package $(\alpha_-,\xi_+)$ on $T^2 \times [r_0,r_1]$ and $T^2 \times [r_2,r_3]$. Then we can write $\xi_+ = \ker(\cos \vphi(x,y,r) dr - \sin \vphi(x,y,r) dx)$ where $\vphi\colon \mathbb{R}^2 \times ([r_0,r_1] \cup [r_2,r_3]) \to \mathbb{R}$.

If $\vphi(2\pi,y,r_1)-\vphi(0,y,r_1) = \vphi(2\pi,y,r_2)-\vphi(0,y,r_2)$ and $\vphi(x,2\pi,r_1)-\vphi(x,0,r_1) = \vphi(x,2\pi,r_2)-\vphi(x,0,r_2)$ then we can extend $\xi_+$ over $T^2 \times [r_1,r_2]$ so that $(\alpha_-,\xi_+)$ is a strongly adapted bicontact package on $T^2 \times [r_0,r_3]$.
\end{lemma}
\begin{proof}
By definition, $\vphi$ is a lift of a function $\vphi\colon \mathbb{R}^2/(2\pi\mathbb{Z})^2 \times ([r_0,r_1] \cup [r_2,r_3]) \to \mathbb{R}/2\pi\mathbb{Z}$, thus $\frac{1}{2\pi}(\vphi(x+2\pi,y,r)-\vphi(x,y,r))$ and $\frac{1}{2\pi}(\vphi(x,y+2\pi,r)-\vphi(x,y,r))$ are integers for every $(x,y) \in \mathbb{R}^2$ and every $r \in [r_0,r_1] \cup [r_2,r_3]$.
Moreover, since $\xi_+$ is positive contact, $\vphi$ is increasing in $y$ by \Cref{prop:contactstructurerudder}.
If we have the equalities $\vphi(2\pi,y,r_1)-\vphi(0,y,r_1) = \vphi(2\pi,y,r_2)-\vphi(0,y,r_2)$ and $\vphi(x,2\pi,r_1)-\vphi(x,0,r_1) = \vphi(x,2\pi,r_2)-\vphi(x,0,r_2)$, then we can extend $\vphi$ over $[r_1,r_2]$ so that $\vphi$ remains increasing in $y$, and so that $\frac{1}{2\pi}(\vphi(x+2\pi,y,r)-\vphi(x,y,r))$ and $\frac{1}{2\pi}(\vphi(x,y+2\pi,r)-\vphi(x,y,r))$ are constant integers over $T^2 \times [r_0,r_3]$.
We can thus set $\xi_+ = \ker(\cos \vphi(x,y,r) dr - \sin \vphi(x,y,r) dx)$ on $T^2 \times [r_0,r_3]$.
\end{proof}

\begin{lemma}[Linearizing the positive contact structure] \label{lemma:fillstdorbitstdpos}
Let $T^2 \times [r_0,r_3]$ be a thickened torus with a standard rational negative contact form $\alpha_-$. 
Suppose $\alpha_-$ is part of a strongly adapted bicontact package $(\alpha_-,\xi_+)$ on $T^2 \times [r_2,r_3]$. Then we can extend $\xi_+$ over $T^2 \times [r_0,r_2]$ so that $(\alpha_-,\xi_+)$ is a strongly adapted bicontact package on $T^2 \times [r_0,r_3]$ and $\xi_+$ is linear on $T^2 \times [r_0,r_1]$.
\end{lemma}
\begin{proof}
Choose coordinates $T^2 \cong S^1_x \times S^1_y$ under which $\alpha_- = A(r) dx + B dy$ and $R_- = \frac{1}{B} \frac{\partial}{\partial y}$, and write $\xi_+ = \ker(\cos \vphi(x,y,r) dr - \sin \vphi(x,y,r) dx)$ where $\vphi\colon \mathbb{R}^2 \times [r_2,r_3] \to \mathbb{R}$.

By \Cref{prop:contactstructurerudder}, we have $\frac{\partial \vphi(x,y,r)}{\partial y} > 0$, hence $\vphi(x,2\pi,r) - \vphi(x,0,r) = 2\pi q$ for some positive integer $q$.
We also have $\vphi(2\pi,y,r) - \vphi(0,y,r) = 2\pi p$ for some integer $p$.

We then apply \Cref{lemma:fillstdorbitinterpolatepos} to the positive contact structure $\ker(\cos(px+qy)dr - \sin(px+qy)dx)$ on $T^2 \times [r_0,r_1]$ and the given $\xi_+$ on $T^2 \times [r_2,r_3]$.
\end{proof} %

\begin{lemma}[Rotating the characteristic and flow slopes] \label{lemma:fillstdorbitrotateneg}
Let $T^2 \times [r_2,r_3]$ be a thickened torus with a standard strongly adapted bicontact package $(\alpha_-,\xi_+)$. 
Suppose $\alpha_-$ is rational and $\xi_+$ is linear.
Then, given any standard negative contact form $\alpha^\circ_-$ on $T^2 \times [r_0,r_1]$ such that 
$f(r_2) < f(r_1) < \chi(r_2) < \chi(r_1) < l$ and $f(r) \neq l$ for $r \in [r_0,r_1]$,
we can extend $(\alpha_-,\xi_+)$ over $T^2 \times [r_0,r_2]$ so that $(\alpha_-,\xi_+)$ remains a standard strongly adapted bicontact package, $\alpha_-$ is a rational multiple of $\alpha_-^\circ$ on $T^2 \times [r_0,r_1]$, and $f(r)$ is decreasing for $r \in [r_1,r_2]$.
\end{lemma}
\begin{proof}

Recall the parametrized path $\widetilde{\alpha}$ associated to a standard 1-form $\alpha$. 
We currently have two such paths: $\widetilde{\alpha_-}$, which is defined on $[r_2,r_3]$, and $\widetilde{\alpha_-^\circ}$, which is defined on $[r_0,r_1]$. 
We extend $\widetilde{\alpha_-}$ over $[r_0,r_2]$ so that 
\begin{itemize}
    \item $\widetilde{\alpha_-}$ remains a path that winds around the origin counterclockwise,
    \item $\widetilde{\alpha_-}$ is a rescaling of $\widetilde{\alpha_-^\circ}$ on $[r_0,r_1]$, and
    \item the slope of the derivative $\widetilde{\alpha_-}'$ is decreasing on $[r_1,r_2]$.
\end{itemize}
Here we used the hypotheses that 
$f(r_2) < f(r_1) < \chi(r_2) < \chi(r_1) < l$.
See \Cref{fig:fillstdorbitrotatenegslopes} for an example.

\begin{figure}
    \centering
    \fontsize{10pt}{10pt}\selectfont
    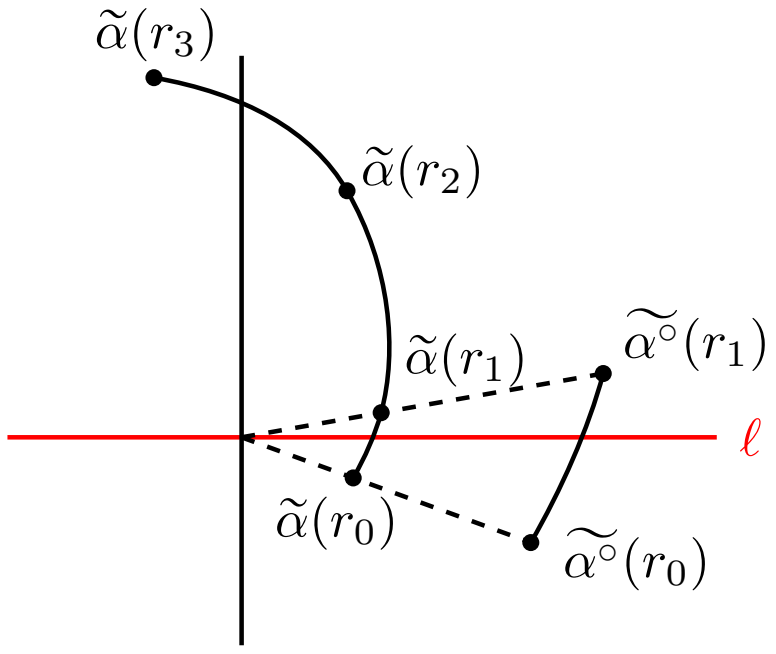
    \caption{Extending the path $\widetilde{\alpha_-}$ in \Cref{lemma:fillstdorbitrotateneg}.}
    \label{fig:fillstdorbitrotatenegslopes}
\end{figure}

It remains to extend $\xi_+$. 
By assumption of linearity, we can choose coordinates $T^2 \cong S^1_x \times S^1_y$ under which $\xi_+ = \ker(\cos(ky) dr - p\sin(ky) dx - q\sin(ky)dy)$ on $T^2 \times [r_2,r_3]$, for some $p,q,k \in \mathbb{Z}$ where $p,k > 0$.
Under such coordinates, the Legendrian slope $l$ is given by $\frac{\partial}{\partial x}$.
We write $\alpha_- = A(r) dx + B(r) dy$, so that $R_- = \frac{1}{A(r)B'(r) - A'(r)B(r)} (B'(r) \frac{\partial}{\partial x} - A'(r) \frac{\partial}{\partial y})$.
Note that $A'(r) \neq 0$ for $r \in [r_2,r_3]$ since the flow slope $f(r)$ never coincides with $l$ where the bicontact package is defined.
Since $f(r)$ is decreasing for $r \in [r_1,r_2]$ and $f(r_2) < f(r_1) < l$, the same is true there.
Finally, the same is also true on $[r_0,r_1]$ by hypothesis.

We set $C(r) = -\frac{B'(r)}{A'(r)}$, so that the flow slope $f(r)$ is in the direction of $C(r)\frac{\partial}{\partial x} + \frac{\partial}{\partial y}$.
In particular, $C'(r) \leq 0$ on $[r_1,r_2]$, and by assumption of rationality, $C(r)$ is constant on $[r_2,r_3]$.
Choose a function $\rho(r)$ such that $\rho(r) \geq 1 + \frac{p}{k} \max\{0,-C'(r)\}$ for $r \in [r_0,r_2]$ and $\rho(r) = 1$ for $r \in [r_2,r_3]$. %

We can now set $\alpha_+ = \rho(r) \cos(ky) dr - p\sin(ky) dx + p C(r) \sin(ky) dy$ and $\xi_+ = \ker \alpha_+$ on $T^2 \times [r_0,r_2]$.
A routine computation shows that 
$$\alpha_+ \wedge d\alpha_+ = p \left(k + \left( k (\rho(r)-1) + pC'(r) \sin^2(ky) \right) \right) dx \wedge dy \wedge dr.$$ 
Since $\rho(r) \geq 1 - \frac{p}{k} C'(r)$, this shows that $\xi_+$ remains a positive contact structure.
\end{proof}

Now consider a solid torus $\mathbf{D} = D^2 \times S^1_z$.

\begin{defn}[Standard Reeb orbit] \label{defn:stdorbit}
We define a \textbf{standard} negative contact form on $\mathbf{D}$ to be
\begin{align*}
\alpha_- &= -r^2 d\theta + (L+ \frac{p}{q} r^2)dz \text{ in polar coordinates $(r,\theta)$ for $D^2$} \\
&= -xdy+ydx+(L+ \frac{p}{q}(x^2+y^2))dz \text{ in Euclidean coordinates $(x,y)$ for $D^2$}
\end{align*}
for some $L \in \mathbb{Q}_+$, $\frac{p}{q} \in \mathbb{Q}$.
Its Reeb vector field is 
\begin{align*}
R_- &= \frac{1}{L} (\frac{p}{q} \frac{\partial}{\partial \theta} + \frac{\partial}{\partial z}) \\
&= \frac{1}{L} (- \frac{p}{q} y\frac{\partial}{\partial x} + \frac{p}{q} x\frac{\partial}{\partial y} + \frac{\partial}{\partial z}).
\end{align*}
In particular, the $z$-axis $\{0\} \times S^1$ is a closed orbit of $R_-$.  
From \Cref{constr:blowup}, one can check that the same expression holds for the blow-up of $R_-$ along the $z$-axis, thus $R_-$ is tangent to each concentric torus $\{r=r_0\}$, and the foliations of the tori by flow lines have the same slope $\frac{p}{q}$. %

Consider the plane field
$$\xi_+ = \ker(\cos(nz)dx - \sin(nz)dy + \frac{p}{q}(y\cos(nz) + x\sin(nz))dz)$$ 
in Euclidean coordinates $(x,y,z)$, where $n \in \mathbb{Z}$.
Note that $R_- \subset \xi_+$, and $\xi_+$ is a positive contact structure if and only if $n+\frac{p}{q} > 0$.
In this case, we refer to $(\alpha_-,\xi_+)$ as a \textbf{standard} strongly adapted bicontact package on $\mathbf{D}$ with flow slope $\frac{p}{q}$ and Legendrian slope $-n$.
We refer to the $z$-axis as a \textbf{standard} Reeb orbit.

More generally, given a strongly adapted bicontact package $(\alpha_-,\xi_+)$ and a closed orbit $\gamma$ of the Reeb flow $R_-$, we say that $\gamma$ is a \textbf{standard} Reeb orbit if there is an isomorphism of bicontact packages between a neighborhood of $\gamma$ and a neighborhood of the $z$-axis under a standard strongly adapted bicontact package on $\mathbf{D}$.
In particular, given any standard strongly adapted bicontact package $(\alpha_-,\xi_+)$ on $\mathbf{D}$ and any $\lambda \in \mathbb{Q}_+$, the $z$-axis remains a standard Reeb orbit under the rescaled bicontact package $(\lambda \alpha_-,\xi_+)$.
\end{defn}

Note that the restriction of a standard strongly adapted bicontact package on $\mathbf{D}$ to $\mathbf{D} \backslash (\{0\} \times S^1)$ is also a standard strongly adapted bicontact package. Moreover, $\alpha_-$ is rational, but $\xi_+$ is \emph{not} linear.
Indeed, if we convert from Euclidean to cylindrical coordinates, we have
$$\xi_+ = \ker(\cos(\theta+nz) dr - r \sin(\theta+nz) d(\theta-\frac{p}{q}z)).$$ %
One can make $\xi_+$ linear by reparametrizing the $r$-coordinate on $\mathbf{D} \backslash (\{0\} \times S^1)$.

\begin{prop} \label{prop:fillstdorbit}
Let $T^2 \times [r_0,r_1]$ be a thickened torus with a standard strongly adapted bicontact package $(\alpha_-,\xi_+)$, where $\alpha_-$ is rational.
Suppose $f(r_1) < \chi(r_1) < l$.
Then for every slope $m$ lying in $(\chi(r_1),l)$ and with intersection number $i(m,l) = 1$, we can extend $T^2 \times [r_1-\epsilon,r_1]$ along with $(\alpha_-,\xi_+)$ into a solid torus $\mathbf{D}$ by topologically filling along the slope $m$, such that $\{0\} \times S^1$ is a standard Reeb orbit.

Moreover, we can perform the extension so that given any multicurve $c$ transverse to $R_-$ on $T^2 \times \{r_1\}$, we can extend it into a local Birkhoff section to $R_-$ within $\mathbf{D}$ which has $\{0\} \times S^1$ as its boundary orbit.
\end{prop}
\begin{proof}
For the duration of this proof, we denote $\mathbf{T} = T^2 \times [r_0,r_1]$ and $\mathbf{D} = D^2 \times S^1$.
We denote the given bicontact package on $\mathbf{T}$ by $(\alpha_-^{\mathbf{T}},\xi_+^{\mathbf{T}})$ and denote its characteristic and flow slopes by $\chi^{\mathbf{T}}(r)$ and $f^{\mathbf{T}}(r)$.

Pick a slope $\zeta$ with $i(\zeta,m) = 1$. 
Pick coordinates $T^2 \cong S^1_\theta \times S^1_z$ such that 
$m$ is $\frac{\partial}{\partial \theta}$ and $\zeta$ is $\frac{\partial}{\partial z}$.
By identifying these with the coordinates $\theta$ and $z$ on $\mathbf{D}$, this allows us to compare slopes between $\mathbf{T}$ and $\mathbf{D}$.
Furthermore, we can identify these slopes with $\mathbb{R}P^1$ by specifying that $m$ corresponds to $\infty$ and $\zeta$ corresponds to $0$. 
Then because $i(m,l)=1$, the Legendrian slope $l$ corresponds to some integer $-n$.

Pick a rational slope $f_0$ in $(f(r_1),\chi(r_1))$, and suppose it corresponds to $\frac{p}{q} \in \mathbb{Q}$. Then since $f_0 < \chi(r_1) < m < l$, we have $\frac{p}{q} > -n$.
We consider the standard bicontact package on $\mathbf{D}$ with flow slope $\frac{p}{q}$ and Legendrian slope $-n$, denoting it by $(\alpha_-^{\mathbf{D}},\xi_+^{\mathbf{D}})$ and denoting its characteristic and flow slopes by $\chi^{\mathbf{D}}(r)$ and $f^{\mathbf{D}}(r)$.

Since $f(r_1) < f_0 < \chi(r_1) < m < l$,
we can take $\epsilon$ small enough so that $f_{\mathbf{T}}(r_1-3\epsilon) < f_{\mathbf{D}}(4\epsilon) < \chi_{\mathbf{T}}(r_1-3\epsilon) < \chi_{\mathbf{D}}(4\epsilon) < l$.

We apply \Cref{lemma:fillstdorbitstdpos} on the restriction of $\alpha_-^{\mathbf{T}}$ to $T^2 \times [r_1-3\epsilon,r_1] \subset \mathbf{T}$ and the restriction of $\xi_+^{\mathbf{T}}$ to $T^2 \times [r_1-\epsilon,r_1]$ to get a new positive contact structure $\widehat{\xi}_+^{\mathbf{T}}$ on $T^2 \times [r_1-3\epsilon,r_1]$ so that $(\alpha_-^{\mathbf{T}},\widehat{\xi}_+^{\mathbf{T}})$ is a strongly adapted bicontact package with $\widehat{\xi}_+^{\mathbf{T}} = \xi_+^{\mathbf{T}}$ on $T^2 \times [r_1-\epsilon,r_1]$ and $\widehat{\xi}_+^{\mathbf{T}}$ linear on $T^2 \times [r_1-3\epsilon,r_1-2\epsilon]$.

We then apply \Cref{lemma:fillstdorbitrotateneg} to the restriction of $(\alpha_-^{\mathbf{T}},\widehat{\xi}_+^{\mathbf{T}})$ to $T^2 \times [r_1-3\epsilon,r_1-2\epsilon]$ and the restriction of $\alpha_-^{\mathbf{D}}$ to $T^2 \times [3\epsilon,4\epsilon]$ to obtain a bicontact package $(\alpha_-,\xi_+)$ on $T^2 \times [3\epsilon,r_1-2\epsilon]$ so that $(\alpha_-,\xi_+) = (\alpha_-^{\mathbf{T}},\widehat{\xi}_+^{\mathbf{T}})$ on $T^2 \times [r_1-3\epsilon,r_1-2\epsilon]$ and $\alpha_- = \lambda \alpha_-^{\mathbf{D}}$ on $T^2 \times [3\epsilon,4\epsilon]$ for some $\lambda \in \mathbb{Q}_+$.

We then apply \Cref{lemma:fillstdorbitinterpolatepos} to the restriction of $\lambda \alpha_-^{\mathbf{D}}$ to $T^2 \times [\epsilon,4\epsilon]$ and the restrictions of $\xi_+$ to $T^2 \times [3\epsilon,4\epsilon]$ and $\xi_+^{\mathbf{D}}$ to $T^2 \times [\epsilon,2\epsilon]$ in order to extend $\xi_+$ to $T^2 \times [\epsilon,4\epsilon]$ in a way so that $(\lambda \alpha_-^{\mathbf{D}},\xi_+)$ is a strongly adapted bicontact package and $\xi_+ = \xi_+^{\mathbf{D}}$ on $T^2 \times [\epsilon,2\epsilon]$.
Here the choice of the Legendrian slope of $\xi_+^{\mathbf{D}}$ as $-n$ comes into play when verifying the equalities in the hypothesis of \Cref{lemma:fillstdorbitinterpolatepos}.

For the `moreover' statement, note that given any slope $\sigma \neq f(r_1)$, we can pick $f_0$ so that $f_0 \in (f(r_1),\sigma)$ as well. In this case, the flow slope $f(r)$ of $(\alpha_-,\xi_+)$ 
\begin{itemize}
    \item coincides with $f^{\mathbf{T}}(r)$ for $r > r_1 - 3\epsilon$, hence is constant at $f(r_1)$ there since $\alpha_-^{\mathbf{T}}$ is rational,
    \item coincides with $f^{\mathbf{D}}(r)$ for $r < 4 \epsilon$, hence is constant at $f_0 \in (f(r_1),\sigma)$ there, and
    \item is decreasing for $r \in [4\epsilon, r_1 - 3\epsilon]$ by \Cref{lemma:fillstdorbitrotateneg},
\end{itemize}
so $f(r)$ never rotates beyond the given slope $\sigma$.

Now given a multicurve $c$ transverse to $R_-$ on $T^2 \times \{r_1\}$, the slope $\sigma$ of $c$ will be different from $f(r_1)$.
The paragraph above ensures that we can choose a smoothly varying family of multicurves $c(r)$ of slope $\sigma$ on the level tori $T^2 \times \{r\}$ of the blow-up $T^2 \times [0,r_1]$ that is transverse to the flow lines.
Moreover, we can choose $c(r_1)$ to be the given multicurve $c$, and arrange for $c(0)$ also to be transverse to the fibering by meridional curves on the innermost torus $\{r=0\}$, or agree with a meridional fiber if $\sigma$ happens to be the meridian.
Then the surface obtained by taking the union of the isotoped curves on each torus is transverse to the blown-up flow, and collapses down to give a local ($\partial$-strong) Birkhoff section of the original flow $R_-$. %
\end{proof}

\subsection{Padding with standard boundary orbits} \label{subsec:paddingsurgery}

In this subsection, we present the bicontact padding surgery operation.

\begin{constr}[Bicontact padding surgery] \label{constr:paddingsurgery}
Let $(\alpha_-,\xi_+)$ be a strongly adapted bicontact package on a closed oriented 3-manifold $M$. 
Suppose the Reeb flow $R_-$ of $\alpha_-$ admits a Birkhoff section $S$.
Furthermore, suppose that every boundary orbit of $S$ is negative and is a standard Reeb orbit.

By definition of a standard Reeb orbit, there exists a tubular neighborhood $\mathbf{D}$ of each boundary orbit $\gamma$ on which 
the restriction of $(\alpha_-,\xi_+)$ has the form $\alpha_- = -r^2 d \theta + (L+ \frac{p}{q} r^2) dz$ and $\xi_+ = \ker(\cos(nz)dx - \sin(nz) dy+\frac{p}{q}(y\cos(nz)+x\sin(nz))dz)$, under some Euclidean coordinate system $(x,y,z)$ or cylindrical coordinate system $(r,\theta,z)$.
In particular, the restriction of $\alpha_-$ to $\mathbf{D} \backslash \gamma$ is standard and rational, and under the framing induced by the $\theta$ coordinate, we have $f = \frac{p}{q} < m = \infty < l = -n$.
Here, we remark that the constants $L, \frac{p}{q}, n$ all depend on $\gamma$, but we suppress this from our notation to avoid clutter.

Since every boundary orbit is negative, the slope $\sigma$ induced by $S$ lies in the interval $(f,m)$.
Up to an isotopy along flow lines, we can assume that $S$ intersects each level torus $\{r=r_0\}$ in evenly spaced straight lines, with respect to the Euclidean structure induced by the coordinates $\theta$ and $z$.

We pick $\epsilon > 0$ small enough so that $f < \sigma < \chi(\epsilon) < l$, and so that $\{r = \epsilon\}$ is a torus in $\mathbf{D}$.
We remove the flow-invariant open solid torus $\nu = \{r < \frac{\epsilon}{2}\}$ around each boundary orbit, 
and denote $M^\circ = M \backslash \nu$ and $\mathbf{T} = \mathbf{D} \backslash \nu$.
We then make a linear coordinate change on $\mathbf{T}$, replacing $(\theta,z) \in S^1 \times S^1$ by $(u,v) \in \mathbb{R}/\mathbb{Z} \times \mathbb{R}/\mathbb{Z}$ so that $S \cap \mathbf{T} = \{v = \frac{k}{K} \mid k = 0,...,K-1\}$, where $K$ is the number of components of $S \cap \mathbf{T}$.
Under these coordinates, we can write $\alpha_- = A(r) du + B(r) dv$ for some functions $A,B$, and $R_- = a\frac{\partial}{\partial u} + b\frac{\partial}{\partial v}$ for some constants $a,b$.
Here, since $L$ is rational, $a,b$ are rational multiples of $\frac{1}{2\pi}$. Also, up to reversing $u$ and $v$, we can assume $b > 0$. %

Meanwhile, let $S^\circ = S \backslash \nu$. 
For small enough $\eta > 0$, $P = R_-^{[-\eta,\eta]}(S^\circ)$ is a subset of $M^\circ$ diffeomorphic to a product $S^\circ \times [-\eta,\eta]_t$.
On this product, we will use the coordinates $(x,t)$ where $x \in S^\circ$ and $t \in [-\eta,\eta]$. 
Observe that we can restrict $(r,u)$ to a coordinate system in the collar neighborhood $S^\circ \cap \mathbf{T}$ of $\partial S^\circ$ in $S^\circ$, and we have $u(x,t) = u(x) + at$ for $x \in S^\circ \cap \mathbf{T}$.
Since the Reeb flow preserves $\alpha_-$, we can write $\alpha_- = \beta + dt$ on $P$, where $\beta$ is a 1-form on $S^\circ$.
Observe that $P$ intersects $\mathbf{T}$ in the subsets $\{v \in [\frac{k}{K} - b\eta, \frac{k}{K} + b\eta]\}$.
See \Cref{fig:paddingsurgery} left for an illustration.

\begin{figure}
    \centering
    \fontsize{6pt}{6pt}\selectfont
    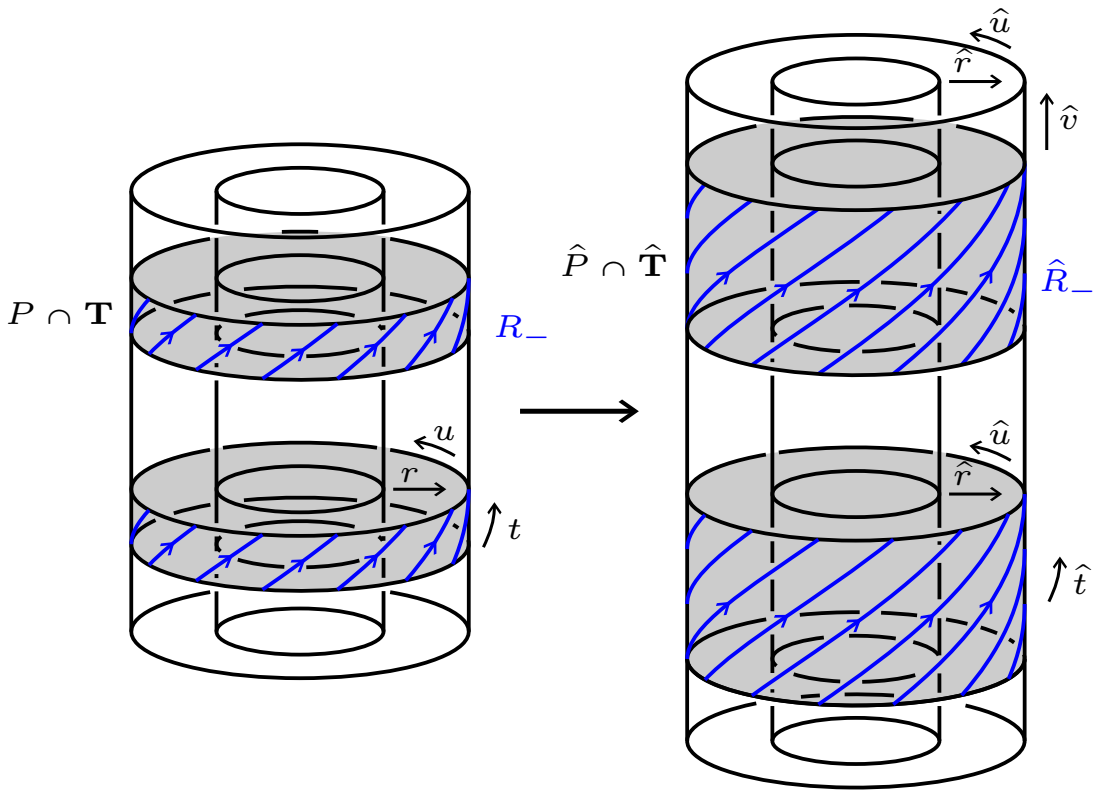
    \caption{Bicontact padding surgery near boundary orbits.}
    \label{fig:paddingsurgery}
\end{figure}

Now let $\widehat{P} = S^\circ \times [-H-\eta,H+\eta]_{\widehat{t}}$,
where $H > 0$ is to be chosen later.
We define $\widehat{\alpha}_- = \beta + d\widehat{t}$, so that its Reeb vector field is $\widehat{R}_- = \frac{\partial}{\partial \widehat{t}}$, and define a positive contact structure $\widehat{\xi}_+$ by specifying that:
\begin{itemize}
    \item After identifying $S^\circ \times [-H-\eta,-H-\frac{1}{2}\eta] \subset \widehat{P}$ with $S^\circ \times [-\eta,-\frac{1}{2}\eta] \subset P$ via $\widehat{t} = t - H$ and identifying $S^\circ \times [H+\frac{1}{2}\eta,H+\eta] \subset \widehat{P}$ with $S^\circ \times [\frac{1}{2}\eta,\eta] \subset P$ via $\widehat{t} = t + H$, we have $\widehat{\xi}_+ = \xi_+$ on these subsets.
    \item For each $x \in S^\circ$, $\widehat{\xi}_+$ completes $N$ more left-handed full rotations along $x \times [-H-\eta,H+\eta]$ than $\xi_+$ does along $x \times [-\eta,\eta]$. 
    More precisely, the family of plane fields $((\widehat{R}_-^{\widehat{t}})^*(\widehat{\xi}_+|_{(x,\widehat{t})}))_{\widehat{t} \in [-H-\eta,H+\eta]}$ rotates in a left-handed way around $\frac{\partial}{\partial \widehat{t}}$ and completes $N$ more full rotations as $\widehat{t}$ increases than $((R_-^t)^*(\xi_+|_{(x,t)}))_{t \in [-\eta,\eta]}$ does as $t$ increases.
    Here, $N$ is a positive integer to be chosen later.
\end{itemize}

We now cut $P$ away from $M^\circ$ and glue in $\widehat{P}$ by identifying $S^\circ \times \{\pm(H+\eta)\} \subset \widehat{P}$ with $S^\circ \times \{\pm\eta\} \subset P$ via $\widehat{t} = t \pm H$.
Since the bicontact packages on $\widehat{P}$ and $P$ agree near the surfaces of gluing, the glued manifold $\widehat{M}^\circ$ inherits a bicontact package $(\widehat{\alpha}_-,\widehat{\xi}_+)$.
Observe that $\widehat{M}^\circ$ is naturally homeomorphic to $M^\circ$, but for our discussion, it will be helpful to regard them as separate manifolds and decorate objects on the former with a hat.
For example, we let $\widehat{\mathbf{T}}$ be the set $(\mathbf{T} \backslash P) \cup ((S \cap \mathbf{T}) \times [-H-\eta,H+\eta])$.

Recall that we have the coordinates $(r,u)$ on $S^\circ \cap \mathbf{T}$.
We define $\widehat{r}(x,\widehat{t}) = r(x)$ and define $\widehat{u}(x,\widehat{t}) = u(x)+a\widehat{t}$ for $(x,\widehat{t}) \in (S^\circ \cap \mathbf{T}) \times [-H-\eta,H+\eta]$.
The functions $r$ and $\widehat{r}$ match up near the surfaces of gluing, that is, if $t \in [-\eta,-\frac{1}{2}\eta]$ then $r(x,t) = \widehat{r}(x,t-H)$ and if $t \in [\frac{1}{2}\eta,\eta]$ then $r(x,t) = \widehat{r}(x,t+H)$. 
We can thus glue them up into a function on $\widehat{\mathbf{T}}$, which we also denote by $\widehat{r}$.

However, the same is not true in general for $u$ and $\widehat{u}$. Instead, if $t \in [-\eta,-\frac{1}{2}\eta]$ then $u(x,t) = \widehat{u}(x,t-H) + aH$ and if $t \in [\frac{1}{2}\eta,\eta]$ then $u(x,t) = \widehat{u}(x,t+H)-aH$. 
In view of this, we pick $H$ so that $aH \in \mathbb{Z}$, i.e. we pick $H$ to be an integer multiple of $\frac{1}{a}$ if $a \neq 0$. Then the functions $u$ and $\widehat{u}$ can now be glued up into a function on $\widehat{\mathbf{T}}$, which we also denote by $\widehat{u}$.
Here we use the fact that although $a$ depends on the boundary orbit $\gamma$, it is always a rational multiple of $\frac{1}{2\pi}$, so that we can pick a common integral multiple of $\frac{1}{a}$ over all boundary orbits.

Similarly, we define a function $\widehat{v}\colon \widehat{\mathbf{T}} \to \mathbb{R}/(2HbK+1)\mathbb{Z}$, acting as the analogue of the $v$ coordinate on $\mathbf{T}$, by specifying that
\begin{itemize}
    \item $\widehat{v}(r,u,v) = v + (2k+1)Hb$ if $(r,u,v) \in \mathbf{T} \backslash P$ and $\frac{k}{K} + b\eta \leq v \leq \frac{k+1}{K} - b\eta$, and
    \item $\widehat{v}(x,\widehat{t}) = \widehat{v}(x,-H-\eta) + b(\widehat{t}+H+\eta)$ if $(x,\widehat{t}) \in (S \cap \mathbf{T}) \times [-H-\eta,H+\eta]$.
\end{itemize}
Then $(\widehat{r},\widehat{u},\widehat{v})$ is a coordinate system on $\widehat{\mathbf{T}}$.
See \Cref{fig:paddingsurgery} right.

Observe that on $\widehat{\mathbf{T}}$, we have $\widehat{\alpha}_- = A(\widehat{r})d\widehat{u} + B(\widehat{r})d\widehat{v}$, where $A,B$ are the same functions as those for $\alpha_-$, and its Reeb vector field $\widehat{R}_- = a \frac{\partial}{\partial \widehat{u}} + b \frac{\partial}{\partial \widehat{v}}$, where $a,b$ are the same constants as those for $R_-$.
In particular $\widehat{\alpha}_-$ is standard, hence so is $(\widehat{\alpha}_-,\widehat{\xi}_+)$, and $\widehat{\alpha}_-$ is rational, because if $a \frac{\partial}{\partial u} + b \frac{\partial}{\partial v}$ is rational on $(\mathbb{R}/\mathbb{Z})_u \times (\mathbb{R}/\mathbb{Z})_v$, then $a \frac{\partial}{\partial \widehat{u}} + b \frac{\partial}{\partial \widehat{v}}$ is rational on $(\mathbb{R}/\mathbb{Z})_{\widehat{u}} \times (\mathbb{R}/(2HbK+1)\mathbb{Z})_{\widehat{v}}$ since $H$ is an integer multiple of $\frac{1}{a}$.

Note that the thickened torus obtained by collapsing all flow lines of $R_-$ in $P \cap \mathbf{T}$ can be canonically identified with that obtained by collapsing all flow lines of $\widehat{R}_-$ in $\widehat{P} \cap \widehat{\mathbf{T}}$. Through this, we can compare slopes and homology classes on $\mathbf{T}$ and $\widehat{\mathbf{T}}$.
Under this convention, the new Legendrian class $\widehat{l}$ is the sum of the old one $l$ and $N$ copies of the multicurve $\partial S^\circ$ (where we coorient all curves using $R_-$).
Since the old bicontact package was a restriction of one on a solid torus, the intersection number between $l$ and the meridian is $1$, so $l$ is primitive.
Thus there exists arbitrarily large values of $N$ for which $\widehat{l}$ is also primitive, for example when $N$ is a multiple of the intersection number $i(l,\partial S^\circ)$.
We fix such a choice of $N$. 

Recall that we have $f < \sigma < \chi(\epsilon) < l$.
Up to increasing $H$ (while maintaining the property that it is an integer multiple of $\frac{1}{a}$), the slope $\widehat{\chi}(\epsilon)$ can be made arbitrarily close to $\sigma$, so we have $\widehat{f} = f < \sigma < \widehat{\chi}(\epsilon) < \widehat{l}$.
We can then pick a slope $\widehat{m}$ lying in $(\widehat{\chi}(\epsilon),\widehat{l})$ and with intersection number $i(\widehat{m},\widehat{l})=1$.

Finally, we apply \Cref{prop:fillstdorbit} and extend the new bicontact structure onto the closed 3-manifold $\widehat{M}$ obtained by Dehn filling $\widehat{M}^\circ \cong M^\circ$ along the slope $\widehat{m}$.
In particular, note that $\widehat{M}$ is a Dehn surgery of $M$ along $\partial S$. 

We also apply the `moreover' statement of \Cref{prop:fillstdorbit} to complete $S^\circ \times \{0\} \subset \widehat{P} \subset \widehat{M}^\circ$ into a Birkhoff section $\widehat{S}$ of $\widehat{R}_-$ on $\widehat{M}$, where every boundary orbit is a standard Reeb orbit.
Since $\widehat{M}^\circ \cong M^\circ$, the monodromy of $\widehat{S}$ is the same as that of $S$. %
Since $\sigma \in (\widehat{f},\widehat{m})$, every boundary orbit of the new Birkhoff section $\widehat{S}$ is negative.

We refer to this construction as \textbf{bicontact padding surgery} of $(\alpha_-,\xi_+)$ along $S$.
We refer to the number $H$ as the \textbf{height} of the surgery and the number $N$ as the \textbf{degree} of the surgery.
We will also refer to the subset $\widehat{P}$ as the \textbf{padding region}.
\end{constr}

\begin{rmk}[Surgery slopes of bicontact padding surgery] \label{rmk:paddingsurgeryslope}
We work out more explicitly the surgery slopes that are permissible in bicontact padding surgery.
Let us perform computations at a fixed boundary orbit. 
Using the same notation as in \Cref{constr:paddingsurgery}, the flow slope is $n+\frac{p}{q}$ under the Legendrian framing, which is the framing we will use throughout this computation.
Suppose the boundary slope of $S$ is $\frac{c}{d}$, where $c, d > 0$.
Then the new Legendrian class $\widehat{l}$ has slope 
$\frac{Nc}{Nd-1}$,
where $N$ has to be such that
$\gcd(Nc,Nd-1) = 1$.
Then the surgery slope $\widehat{m}$ can be 
$\frac{r+kNc}{s+k(Nd-1)}$,
where $Ncs-(Nd-1)r = 1$, and $k$ is sufficiently large such that $\widehat{m} \in (\widehat{\chi}(\epsilon),\widehat{l})$.

In particular, the set of permissible surgery slopes is a `2-parameter family'.
\end{rmk}

\begin{rmk}[Comparison with Legendrian-transverse surgery] \label{rmk:paddingcomparewithltsurgery}
Each boundary orbit of a Birkhoff section to $R_-$ is also a $\xi_+$-Legendrian $\xi_-$-transverse curve in the sense of \Cref{constr:legendriantransversesurgery}, thus in the setting of \Cref{constr:paddingsurgery}, we can also do Legendrian-transverse surgery along $\partial S$.
While both bicontact padding surgery along $S$ and Legendrian-transverse surgery along $\partial S$ end up defining strongly adapted bicontact packages on 3-manifolds obtained by Dehn surgery on $M$, we emphasize that they are completely different.

Cosmetically, bicontact padding surgery involves modification over a neighborhood of a 2-dimensional subset, whereas Legendrian-transverse surgery involves modification over a neighborhood of a 1-dimensional subset. Thus one should expect bicontact padding surgery to have a more dramatic effect on the dynamics of the supported Anosov flow.

More fundamentally, the sets of surgery slopes are also different.
For Legendrian-transverse surgery, the surgery slopes are $\frac{1}{n}$ where $n \in \mathbb{Z}_-$, which is in particular a `1-parameter family'.
For bicontact padding surgery, the set of surgery slopes is described in \Cref{rmk:paddingsurgeryslope}, and is in particular a `2-parameter family' in $(N,k)$.

However, we do conjecture that the parameter $k$ corresponds to subsequent Legendrian surgery, in the sense that first performing bicontact padding surgery with parameters $(N,k)$, then performing $-\frac{1}{j}$-Legendrian-transverse surgery along the cores of the filling tori produces a bicontact package equivalent to the one obtained by performing bicontact padding surgery with parameters $(N,k+j)$.
Here, equivalent means that the bicontact packages support isotopically equivalent Anosov flows.

We also point out that the setting of bicontact padding surgery is much more restrictive than that of Legendrian-transverse surgery. 
For example, a necessary condition for a collection of orbits to be boundary orbits of a Birkhoff section is for some linear combination of them to be null-homologous, and in general this is very far from sufficient.
\end{rmk}

The next proposition records the properties of the padded bicontact package that will matter most in our application (\Cref{thm:paddingaddwinding0twists}).

\begin{prop} \label{prop:paddingproperties}
Let $(\alpha_-,\xi_+)$ be a strongly adapted bicontact package on a closed oriented 3-manifold $M$. Let $S$ be a Birkhoff section to $R_-$ whose boundary orbits are all negative standard Reeb orbits.
Let $\beta = \alpha_-|_S$ be the 1-form on $S$ obtained by restricting $\alpha_-$.
Let $(\widehat{\alpha}_-,\widehat{\xi}_+)$ be a bicontact package obtained by bicontact padding surgery on $(\alpha_-,\xi_+)$ along $S$.
Then:
\begin{enumerate}
    \item Up to increasing the height $H$ of the surgery, we can embed an arbitrarily large collar neighborhood $S^\circ \times [-H,H]_{\widehat{t}}$ of $S^\circ$ in $\widehat{M}$, so that $\widehat{\alpha}_-$ pulls back to $\beta + d\widehat{t}$ on $S^\circ \times [-H,H]$.
    \item Up to increasing the degree $N$ of the surgery, the plane fields $((\widehat{R}_-^{\widehat{t}})^*(\widehat{\xi}_+|_{\widehat{R}_-^{\widehat{t}}(x)}))_{\widehat{t} \in [-H,H]}$ complete an arbitrarily large number of left-handed rotations for every $x \in S^\circ$. \qed
\end{enumerate}
\end{prop}

\begin{rmk} \label{rmk:partialBirkhoffsection}
\Cref{constr:paddingsurgery} works even when $S$ is just a \textbf{partial Birkhoff section}, i.e., a Birkhoff section but not necessarily meeting every forward and backward trajectory, as long as its boundary orbits are still negative and standard. 
\end{rmk}

\section{Modifying the monodromy} \label{sec:modifymonodromy}

In this section, we explain how to modify the monodromy of a Birkhoff section using bicontact padding surgery.
In \Cref{subsec:windingnumber}, we set up the terminology of winding numbers. 
In \Cref{subsec:factoroutwindzerocurve}, we show that we can add Dehn twists along winding number zero curves.
In \Cref{subsec:factorizationthm}, we state a theorem of Calderon--Salter that says that such winding number zero twists offer a large amount of freedom in modifying the monodromy.

\subsection{Winding numbers} \label{subsec:windingnumber}

For the duration of this subsection, we fix the following setting:
Let $S$ be a compact surface and let $\ell$ be a nonsingular oriented line field on $\intr(S)$.

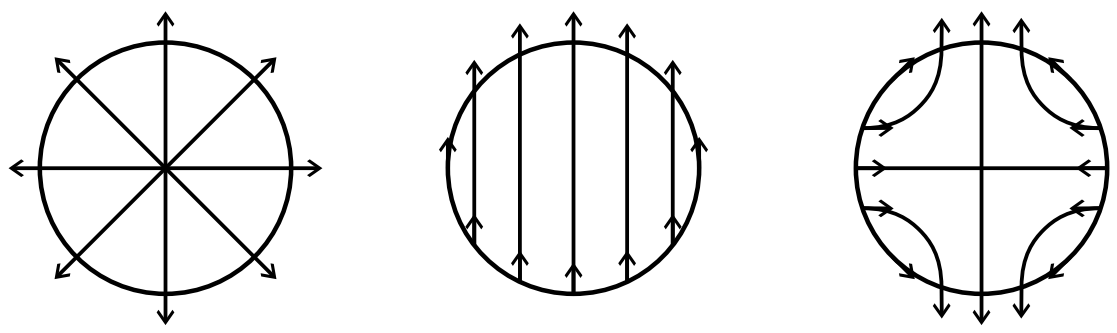
\begin{figure}
    \centering
    \fontsize{10pt}{10pt}\selectfont
    \begingroup%
  \makeatletter%
  \providecommand\color[2][]{%
    \errmessage{(Inkscape) Color is used for the text in Inkscape, but the package 'color.sty' is not loaded}%
    \renewcommand\color[2][]{}%
  }%
  \providecommand\transparent[1]{%
    \errmessage{(Inkscape) Transparency is used (non-zero) for the text in Inkscape, but the package 'transparent.sty' is not loaded}%
    \renewcommand\transparent[1]{}%
  }%
  \providecommand\rotatebox[2]{#2}%
  \newcommand*\fsize{\dimexpr\f@size pt\relax}%
  \newcommand*\lineheight[1]{\fontsize{\fsize}{#1\fsize}\selectfont}%
  \ifx\svgwidth\undefined%
    \setlength{\unitlength}{211.97486457bp}%
    \ifx\svgscale\undefined%
      \relax%
    \else%
      \setlength{\unitlength}{\unitlength * \real{\svgscale}}%
    \fi%
  \else%
    \setlength{\unitlength}{\svgwidth}%
  \fi%
  \global\let\svgwidth\undefined%
  \global\let\svgscale\undefined%
  \makeatother%
  \begin{picture}(1,0.36810226)%
    \lineheight{1}%
    \setlength\tabcolsep{0pt}%
    \put(0,0){\includegraphics[width=\unitlength,page=1]{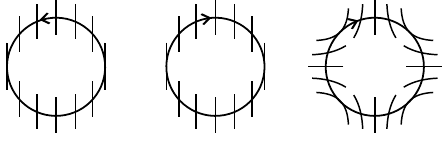}}%
    \put(-0.00175494,0.0060537){\color[rgb]{0,0,0}\makebox(0,0)[lt]{\lineheight{1.25}\smash{\begin{tabular}[t]{l}$\wind_\ell(c)=1$\end{tabular}}}}%
    \put(0.34498318,0.00605268){\color[rgb]{0,0,0}\makebox(0,0)[lt]{\lineheight{1.25}\smash{\begin{tabular}[t]{l}$\wind_\ell(c)=-1$\end{tabular}}}}%
    \put(0.69879763,0.00447222){\color[rgb]{0,0,0}\makebox(0,0)[lt]{\lineheight{1.25}\smash{\begin{tabular}[t]{l}$\wind_\ell(c)=-2$\end{tabular}}}}%
  \end{picture}%
\endgroup%

    \caption{Some examples of winding numbers.}
    \label{fig:windeg}
\end{figure}

\begin{defn}[Winding number] \label{defn:windingnumber}
The \textbf{winding number} of an oriented, possibly immersed, curve $c$ on $\intr(S)$ with respect to $\ell$, denoted by $\wind_\ell(c)$, is the signed number of rotations completed by the derivative $c'$ along $c$ with respect to $\ell$.
We illustrate some examples in \Cref{fig:windeg}.

Note that $\wind_\ell(c)$ is invariant under regular homotopy of $c$, i.e., homotopy through smoothly immersed curves.
Using this fact, we can also define the \textbf{winding number} of a boundary component $b$ of $S$ to be the winding number of a slight push-off of $b$ into $\intr(S)$.

We define the \textbf{winding number} of an oriented multicurve to be the sum of winding numbers of its components.
\end{defn}

\begin{figure}
    \centering
    \fontsize{10pt}{10pt}\selectfont
    \begingroup%
  \makeatletter%
  \providecommand\color[2][]{%
    \errmessage{(Inkscape) Color is used for the text in Inkscape, but the package 'color.sty' is not loaded}%
    \renewcommand\color[2][]{}%
  }%
  \providecommand\transparent[1]{%
    \errmessage{(Inkscape) Transparency is used (non-zero) for the text in Inkscape, but the package 'transparent.sty' is not loaded}%
    \renewcommand\transparent[1]{}%
  }%
  \providecommand\rotatebox[2]{#2}%
  \newcommand*\fsize{\dimexpr\f@size pt\relax}%
  \newcommand*\lineheight[1]{\fontsize{\fsize}{#1\fsize}\selectfont}%
  \ifx\svgwidth\undefined%
    \setlength{\unitlength}{206.40468808bp}%
    \ifx\svgscale\undefined%
      \relax%
    \else%
      \setlength{\unitlength}{\unitlength * \real{\svgscale}}%
    \fi%
  \else%
    \setlength{\unitlength}{\svgwidth}%
  \fi%
  \global\let\svgwidth\undefined%
  \global\let\svgscale\undefined%
  \makeatother%
  \begin{picture}(1,0.35287944)%
    \lineheight{1}%
    \setlength\tabcolsep{0pt}%
    \put(0,0){\includegraphics[width=\unitlength,page=1]{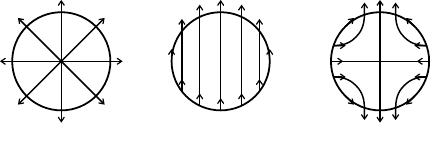}}%
    \put(0.00807646,0.00918582){\color[rgb]{0,0,0}\makebox(0,0)[lt]{\lineheight{1.25}\smash{\begin{tabular}[t]{l}$\wind_\ell(c)=0$\end{tabular}}}}%
    \put(0.35690618,0.00918582){\color[rgb]{0,0,0}\makebox(0,0)[lt]{\lineheight{1.25}\smash{\begin{tabular}[t]{l}$\wind_\ell(c)=-1$\end{tabular}}}}%
    \put(0.72753724,0.00918582){\color[rgb]{0,0,0}\makebox(0,0)[lt]{\lineheight{1.25}\smash{\begin{tabular}[t]{l}$\wind_\ell(c)=-2$\end{tabular}}}}%
  \end{picture}%
\endgroup%

    \caption{We can extend $\ell$ into a vector field with a zero of index $\wind_\ell(c)+1$ in the disc bounded by each boundary component $c$.}
    \label{fig:windboundaryextend}
\end{figure}

\begin{prop} \label{prop:windpoincarehopf}
Let $S$ be a compact surface and let $\ell$ be a nonsingular orientable line field on $\intr(S)$. For any subsurface $W$ of $S$, we have $\chi(W) = \wind_\ell(\partial W)$.
\end{prop}
\begin{proof}
Let $\overline{W}$ be the closed surface obtained by gluing a disc along each boundary component of $W$. We can extend $\ell$ into a vector field on $\overline{W}$ with a zero of index $\wind(c)+1$ in the disc bounded by each boundary component $c$ (where $c$ is oriented as a boundary component of $W$), see \Cref{fig:windboundaryextend}.
Now applying the Poincaré--Hopf theorem, we have 
$$\chi(W) + |\partial W| = \chi(\overline{W}) = \sum_{c \subset \partial W} (\wind_\ell(c)+1) = \wind_\ell(\partial W)+|\partial W|$$
which implies the proposition.
\end{proof}

\begin{defn}[Framed mapping class group] \label{defn:framedmappingclass}
We say that two nonsingular oriented line fields $\ell_0$ and $\ell_1$ on $\intr(S)$ are \textbf{homotopic} if they can be extended into a continuous family of nonsingular oriented line fields $(\ell_s)_{s \in [0,1]}$ on $\intr(S)$.

We define the \textbf{framed mapping class group} $\PMod(\intr(S))[\ell]$ to be the group of orientation-preserving diffeomorphisms $f$ on $\intr(S)$ that act trivially on the set of punctures of $\intr(S)$ and for which $f_* \ell$ is homotopic to $\ell$, modulo isotopy.
\end{defn}

\subsection{Twisting along winding number zero curves} \label{subsec:factoroutwindzerocurve}

Suppose $(\alpha_-,\xi_+)$ is a strongly adapted bicontact package. Let $\phi$ be the Reeb flow generated by the Reeb vector field $R_-$ of $\alpha_-$.
Let $S$ be a Birkhoff section for $\phi$. 
Then since $R_- \subset \xi_+$ and $R_-$ is positively transverse to $\intr(S)$, we have an induced nonsingular oriented line field $\ell = \xi_+ \cap TS$ on $\intr(S)$.

In this context, we will write $\wind_{\xi_+}(S,c) = \wind_\ell(c)$ to emphasize that these numbers arise from $\xi_+$, and to contrast with an (a priori) different notion of winding number to be defined in \Cref{subsec:birkhoffsectionanosovflow}.

If $S$ has only one boundary component lying along each boundary orbit, then $f$ preserves the punctures of $\intr(S)$.
In this case, the monodromy $f$ of $S$ lies in $\PMod(\intr(S))[\ell]$: A homotopy from $\ell$ to $f_* \ell$ is obtained by flowing $\intr(S)$ backwards along $R_-$ and taking its intersection with $\xi_+$.

We now deduce the following results from \Cref{prop:paddingproperties}.

\begin{prop} \label{prop:paddingrealizeanosovorbit}
Let $(\alpha_-,\xi_+)$ be a strongly adapted bicontact package on a closed oriented 3-manifold $M$. Let $S$ be a Birkhoff section to $R_-$ whose boundary orbits are all negative standard Reeb orbits.
Then for any finite (possibly repeating) sequence of curves $(c_1,...,c_m)$ on $S$ with $\wind_{\xi_+}(S,c_i) = 0$ for each $i=1,...,m$,
there exists a bicontact package $(\widehat{\alpha}_-,\widehat{\xi}_+)$ on a closed oriented 3-manifold $\widehat{M}$, obtained from $(\alpha_-,\xi_+)$ by bicontact padding surgery, then homotopy, such that:
\begin{itemize}
    \item There is a Birkhoff section $\widehat{S}$ for the Reeb flow $\widehat{R}_-$ such that the monodromies of $S$ and $\widehat{S}$ are equal.
    \item There exists a topological embedding $S \times [0,m+1] \hookrightarrow \widehat{M}$ such that for every $x \in S$, $\{x\} \times [0,m+1]$ lie along flow-lines of $\widehat{R}_-$. %
    \item For each $i=1,...,m$, there exists a curve $\widehat{c}_i$ on $S \times \{i\}$ homotopic to $c_i \times \{i\}$, and which is Legendrian with respect to both $\widehat{\xi}_+$ and $\widehat{\xi}_- = \ker \widehat{\alpha}_-$, thus is a closed orbit of the Anosov flow supported by $(\widehat{\alpha}_-,\widehat{\xi}_+)$. Moreover, the degeneracy slope of $\widehat{c}_i$ is the framing induced by the surface $S \times \{i\}$.
\end{itemize}
\end{prop}
\begin{proof}
For simplicity, we first explain the proof when $m=1$.
We first perform bicontact $N$-padding surgery along $S$ to get a padded bicontact package $(\widehat{\alpha}_-,\widehat{\xi}_+)$ on $\widehat{M}$.
In the notation of \Cref{constr:paddingsurgery}, up to shrinking the solid torus $\nu$ that we cut out, we can assume that $c_1$ lies on the truncated surface $S^\circ$.
There are then two steps:
\begin{enumerate}
    \item We isotope $S^\circ$ so that $c_1$ is Legendrian with respect to $\widehat{\xi}_- = \ker \widehat{\alpha}_-$.
    \item We isotope $\widehat{\xi}_+$ along flow lines of $\widehat{R}_-$ so that $c_1$ is Legendrian with respect to $\widehat{\xi}_+$ as well. 
\end{enumerate}

For step 1, let $\beta = \alpha_-|_S$ and consider $S^\circ \times \mathbb{R}_t$ with the $t$-invariant contact form $\alpha_- = \beta + dt$.
Note that by \Cref{prop:windpoincarehopf}, a winding number zero curve is automatically non-isolating.
Thus by \Cref{thm:legendrianrealization}, there is an embedding $i\colon S^\circ \hookrightarrow S^\circ \times [-H,H] \subset S^\circ \times \mathbb{R}$, transverse to $\frac{\partial}{\partial t}$, isotopic to $S^\circ \times \{0\}$, and fixed in a collar neighborhood of the boundary, such that $i(c_1)$ is Legendrian with respect to $\xi_- = \ker \alpha_-$.
The assumptions that all boundary orbits of $S$ are negative and that $c_1$ is non-isolating are used here.
By \Cref{prop:paddingproperties}(1), up to increasing the height, there is a contact embedding of $S^\circ \times [-H,H]$ into a collar neighborhood $U$ of $S^\circ$ in $\widehat{M}$.
We can thus isotope $S^\circ$ in $\widehat{M}$ to the embedded copy of $i(S^\circ)$ in $U$. Once we do this, $c_1$ is Legendrian with respect to $\widehat{\xi}_-$.

Note that since the line fields induced by $\widehat{\xi}_+$ and $\xi_+$ on $S^\circ$ are homotopic, we have $\wind_{\widehat{\xi}_+}(S^\circ,c_1) = \wind_{\xi_+}(S,c_1) = 0$.

We move on to step 2.
After isotoping $S^\circ$, $U$ can be written as $\{\widehat{R}_-^t(x) \mid x \in S^\circ, t \in [H_-(x),H_+(x)]\}$.
For each $x \in c_1$, consider the planes $((\widehat{R}_-^{\widehat{t}})^*(\widehat{\xi}_+|_{\widehat{R}_-^{\widehat{t}}(x)}))_{\widehat{t} \in [H_-(x),H_+(x)]}$. 
By \Cref{prop:paddingproperties}(2), up to increasing the degree $N$, this family completes an arbitrarily large number of rotations.
Together with the condition that $\wind_{\widehat{\xi}_+}(S^\circ,c_1) = 0$, this ensures that we can find a function $\widehat{t}_1\colon c_1 \to \mathbb{R}$ such that for every $x \in c_1$, the plane $(\widehat{R}_-^{\widehat{t}_1(x)})^*(\widehat{\xi}_+|_{\widehat{R}_-^{\widehat{t}_1(x)}(x)})$ contains the line $Tc_1|_x$. 
Here, more specifically, the condition that $\wind_{\xi_+}(c_1) = 0$ ensures that $\widehat{t}_1$ is continuous, while increasing $N$ ensures that $\widehat{t}_1(x) \in [H_-(x),H_+(x)]$.

We can now isotope $\widehat{\xi}_+$ along flow lines of $\widehat{R}_-$, pulling back $\widehat{\xi}_+|_{\widehat{R}_-^{\widehat{t}_1(x)}(x)}$ to $x$ for every $x \in c_1$, so that $c_1$ becomes Legendrian with respect to $\widehat{\xi}_+$. %
See \Cref{fig:homotopeforlegendrian}.
By \Cref{prop:poscontactstructureflowisotopy}, after isotopy along flow lines, $(\widehat{\alpha}_-,\widehat{\xi}_+)$ is still a strongly adapted bicontact package.

Since $\widehat{\xi}_+$ is transverse to $TS^\circ$, the degeneracy slope on $c_1$, which is equal to its Legendrian framing with respect to $\widehat{\xi}_+$, must coincide with the framing induced by $S^\circ$. 
This concludes the proof when $m=1$.

More generally, if we are given $m$ curves $c_1,...,c_m$, we divide the large collar neighborhood provided by \Cref{prop:paddingproperties}(1) into $m$ sub-neighborhoods, then run the above strategy for each $c_i$ individually on each sub-neighborhood.
The fact that the curves $\widehat{c}_i$ in the statement of the proposition are only homotopic to $c_i \times \{i\}$ and cannot be made to coincide with them exactly is because of the different isotopies of $S^\circ$ in each sub-neighborhood.
\end{proof}

The next result is \Cref{thm:paddingaddwinding0twists}, but with the added assumption that the boundary orbits are standard. The extension to \Cref{thm:paddingaddwinding0twists} will be done in \Cref{subsec:standardizeorbit} after explaining \Cref{prop:standardizeorbit}. %
\begin{thm} \label{thm:paddingaddwinding0twistsstdboundary}
Let $(\alpha_-,\xi_+)$ be a strongly adapted bicontact package on a closed oriented 3-manifold $M$. Let $S$ be a Birkhoff section to $R_-$ whose boundary orbits are all negative standard Reeb orbits, and let $f$ be the monodromy of $S$.
Then for any finite (possibly repeating) sequence of curves $(c_1,...,c_m)$ on $S$ with $\wind_{\xi_+}(S,c_i) = 0$ for each $i=1,...,m$,
and any sequence of integers $(n_1,...,n_m)$,
there exists a bicontact package $(\widehat{\alpha}_-,\widehat{\xi}_+)$, obtained from $(\alpha_-,\xi_+)$ by bicontact padding surgery, then homotopy, then Goodman surgery, for which there is a Birkhoff section $\widehat{S}$ for the Reeb flow $\widehat{R}_-$, such that the monodromy of $\widehat{S}$ is $\tau_{c_1}^{n_1} ... \tau_{c_m}^{n_m} f$.
\end{thm}
\begin{proof}
We take the bicontact package furnished by \Cref{prop:paddingrealizeanosovorbit} and perform $\frac{1}{n_i}$-Goodman surgery along the Anosov orbits $\widehat{c}_i$.

More specifically, we perform these Goodman surgeries within the framework of Legendrian-transverse surgery, as explained in \Cref{constr:legendriantransversesurgery}. In particular, the bicontact package is only modified in an arbitrarily small neighborhood of the Anosov orbits $\widehat{c}_i$. Thus the Birkhoff section $\widehat{S}$ provided by \Cref{prop:paddingrealizeanosovorbit} is untouched and remains a Birkhoff section after these surgeries. %
\end{proof}

\subsection{A factorization theorem} \label{subsec:factorizationthm}

The goal of this subsection is to review the results of \cite{CS23}, which we will use in the proof of our main theorem to reduce the monodromy to identity.

\begin{defn}[Auxiliary curve system, {\cite[Definitions 6.8 and 6.11]{CS23}}]
Let $S$ be a compact orientable surface and $\ell$ be a nonsingular oriented line field on $\intr(S)$. 

Let $b$ be a boundary component of $S$. An \textbf{auxiliary curve} for $b$ is an oriented curve $c$ in $\intr(S)$ with $\wind_\ell(c) = \begin{cases} \text{$1$ or $-1$} & \text{if $\wind_\ell(b)$ is odd} \\ \text{$2$ or $-2$} & \text{if $\wind_\ell(b)$ is even} \end{cases}$ and which separates $b$ from the other boundary components of $S$.

Let $b_1$ and $b_2$ be two boundary components of $S$. An \textbf{auxiliary curve} for $\{b_1,b_2\}$ is an oriented curve $c$ in $\intr(S)$ that separates $b_1$ and $b_2$ from the other boundary components of $S$.

An \textbf{auxiliary curve system} for $S$ is a collection of auxiliary curves
\begin{itemize}
    \item for each boundary component of $S$, and 
    \item for each pair $\{b_1,b_2\}$ of boundary components of $S$ such that $\wind_\ell(b_1)$ and $\wind_\ell(b_2)$ are both even.
\end{itemize}
\end{defn}

\begin{thm}[{Calderon--Salter \cite[Corollary 6.12]{CS23}}] \label{thm:factorizationthm}
Let $S$ be a compact orientable surface and $\ell$ be a nonsingular oriented line field on $\intr(S)$ satisfying $\wind_\ell(b) < 0$ for every boundary component $b$ of $S$.
If $S$ has genus $\geq 5$, then there exists a collection $\mathcal{C}$ of winding number zero curves such that for any auxiliary curve system $\mathcal{A}$, given any $f \in \PMod(\intr(S))[\ell]$, there exist curves $c_i \in \mathcal{C}$ and $d_j \in \mathcal{A}$ such that $f = \prod \tau_{c_i}^{n_i} \prod \tau_{d_j}^{m_j}$ for some integers $n_i, m_j$.
\end{thm}
\begin{proof}
\cite[Corollary 6.12]{CS23} states that, under the assumption that $\wind_\ell(b) < -1$ for every boundary component $b$, there exists a collection $\mathcal{C}$ of winding number zero curves such that for any auxiliary curve system $\mathcal{A}$, $\PMod(\intr(S))[\ell]$ is generated by Dehn twists along $\mathcal{C} \cup \mathcal{A}$. We briefly indicate why their proof
\begin{itemize}
    \item carries over when $\wind_\ell(b) = -1$ for some boundary components $b$, and
    \item shows that one can put all the twists along $\mathcal{C}$ in front of the twists along $\mathcal{A}$.
\end{itemize}

We adopt their notation for the duration of this proof, including writing the framed mapping class group as $\PMod^n_g[\overline{\phi}]$.
They consider a group $\Mod^*_{g,n}[\phi]$ for which there is a surjection $p_* \colon \Mod^*_{g,n}[\phi] \to \PMod^n_g[\overline{\phi}]$.
In \cite[Proposition 6.7]{CS23}, they show that there is a short exact sequence
$$1 \to \Mod_{g,n}[\phi] \to \Mod^*_{g,n}[\phi] \overset{D}{\to} PR' \to 1.$$

\cite[Theorem B]{CS23} states that $\Mod_{g,n}[\phi]$ is generated by Dehn twists along any collection $\mathcal{C}$ of winding number zero curves where
\begin{enumerate}
    \item $S$ deformation retracts onto the union of curves in $\mathcal{C}$, 
    \item the intersection graph $\Lambda_{\mathcal{C}}$ of $\mathcal{C}$ is a tree, and
    \item $\Lambda_{\mathcal{C}}$ contains the E6 Dynkin diagram.
\end{enumerate}
\cite[Figure 15]{CS23} shows the existence of such a collection, provided that $\wind_\ell(b) < -1$ for every boundary component $b$.
This result generalizes to the case when $\wind_\ell(b) = -1$ for some boundary components $b$: 
First, one can cap off such boundary components and extend the nonsingular line field to the capped-off surface $\overline{S}$, as in \Cref{fig:windboundaryextend} middle. 
Then one takes a collection $\overline{\mathcal{C}}$ of winding number zero curves satisfying properties (1)--(3) on $\overline{S}$ as in \cite[Figure 15]{CS23}. 

Up to an isotopy we can assume that the elements of $\overline{\mathcal{C}}$ are disjoint from the capping-off discs, so that $\overline{\mathcal{C}}$ still determines a collection of winding number zero curves on $S$. Moreover, $\overline{\mathcal{C}}$ continues to satisfy properties (2) and (3) on $S$.
To arrange for property (1) on $S$, we add in a curve $c$ for every capped-off boundary component $b$, so that $c$ and $b$ cobound a pair of pants with a winding number zero curve that is already in the collection. By \Cref{prop:windpoincarehopf}, $c$ has winding number zero, and one can check that properties (2) and (3) are preserved by adding $c$.

Meanwhile, the proof of \cite[Proposition 6.7]{CS23} shows that for any auxiliary curve system $\mathcal{A}$, there exist certain elements in $\Mod^*_{g,n}[\phi]$, written as $A_k, B_{i,j}$, whose images under $D$ generate $PR'$, and whose images in $\PMod^n_g[\overline{\phi}]$ are exactly Dehn twists along the given auxiliary curves.

Thus every element of $\Mod^*_{g,n}[\phi]$ can be written as a product of twists along $\mathcal{C}$ in front of a word in $A_k, B_{i,j}$. Taking the image under $p_*$ gives the required statement.
\end{proof}

\begin{cor} \label{cor:factorizationthmoneboundary}
Let $S$ be a compact orientable surface with one boundary component $b$ and $\ell$ be a nonsingular oriented line field on $\intr(S)$ satisfying $\wind_\ell(b) < 0$. 
If $S$ has genus $\geq 5$ then there exists a collection $\mathcal{C}$ of winding number zero curves such that the Dehn twists along the curves in $\mathcal{C}$ generate $\PMod(\intr(S))[\ell]$.
\end{cor}
\begin{proof}
By \Cref{prop:windpoincarehopf}, if $S$ just has one boundary component $b$, then $b$ has odd winding number. Since there are no other boundary components, we can simply pick a nullhomotopic curve as an auxiliary curve for $b$, of which Dehn twists along are of course isotopic to identity.
\end{proof}

For the purposes of the discussion in \Cref{subsec:variation}, we will also prove the following results when $S$ has more than one boundary component.

\begin{lemma} \label{lemma:disjointauxiliarycurvescriterion}
Let $S$ be a compact orientable surface with boundary and $\ell$ be a nonsingular oriented line field on $\intr(S)$ satisfying $\wind_\ell(b) < 0$ for every boundary component $b$ of $S$.
If $S$ has at most two boundary components with even winding number, then there exists an auxiliary curve system $\mathcal{A}$ for $S$ consisting of disjoint curves.
Moreover, the complementary regions of $\mathcal{A}$ in $S$ can each be assumed to be either an annulus or have genus $\geq 1$.

Conversely, if $S$ has more than two boundary components with even winding number, then there does not exist an auxiliary curve system for $S$ consisting of disjoint curves.
\end{lemma}
\begin{proof}
For each boundary component $b$, \Cref{prop:windpoincarehopf} implies that any curve $c_b$ that cobounds a genus $g_b = \begin{cases} \frac{-\wind_{\ell}(b)-1}{2} & \text{if $\wind_{\ell}(b)$ is odd} \\ \frac{-\wind_{\ell}(b)-2}{2} & \text{if $\wind_{\ell}(b)$ is even} \end{cases}$ subsurface $W_b$ with $b$ is an auxiliary curve for $b$.

Meanwhile, using \Cref{prop:windpoincarehopf} again, we compute
\begin{align*}
\sum_b g_b &= \sum_b \frac{-\wind_{\ell}(b)-1}{2} - \frac{1}{2} ~\text{\# ($\partial$ components with even winding number)} \\
&= \frac{- \chi(S) - \text{\# $\partial$ components}}{2} - \frac{1}{2} ~\text{\# ($\partial$ components with even winding number)} \\
&= \genus(S) - 1 - \frac{1}{2} ~\text{\# ($\partial$ components with even winding number)} < \genus(S).
\end{align*}
So we can pick $W_b$, thus $c_b$, to be disjoint from each other. See the red curves in \Cref{fig:disjointauxiliarycurvescriterion}.

\begin{figure}
    \centering
    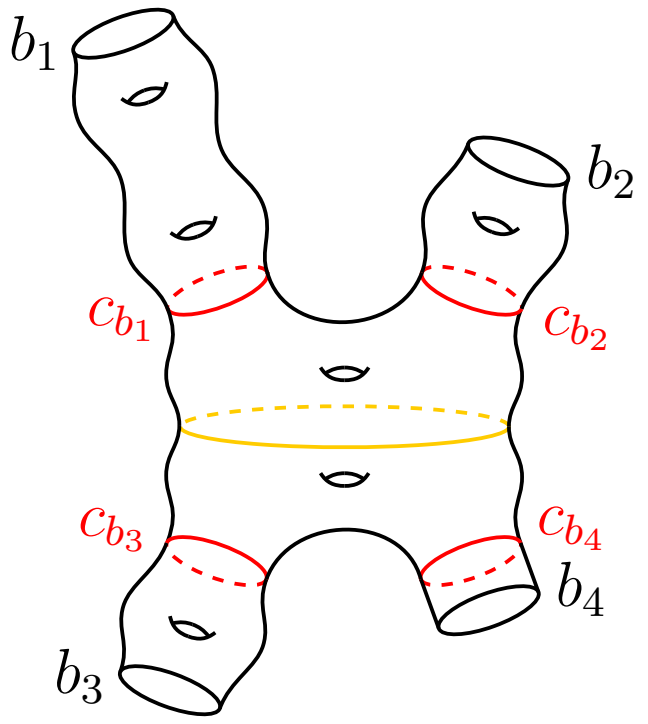
    \caption{The auxiliary curve system in \Cref{lemma:disjointauxiliarycurvescriterion}.}
    \label{fig:disjointauxiliarycurvescriterion}
\end{figure}

From \Cref{prop:windpoincarehopf}, one can deduce that the number of boundary components with even winding number must be even.
If there are no such boundary components, then we are done. 
If there are two such boundary components $b_1$ and $b_2$, then by the computation above, $S \backslash \bigcup_b W_b$ has genus $2$, so we can pick an auxiliary curve $c$ for $\{b_1,b_2\}$ to be a curve in $S \backslash \bigcup_b W_b$ that cobounds a genus one subsurface with $c_{b_1}$ and $c_{b_2}$. See the yellow curve in \Cref{fig:disjointauxiliarycurvescriterion}.

If there are three boundary components $b_1, b_2, b_3$ with even winding number, then any auxiliary curve system has to contain auxiliary curves for $\{b_1,b_2\}$ and for $\{b_2,b_3\}$, and these must intersect each other.
\end{proof}

\begin{cor} \label{cor:factorizationthmtwoevenboundary}
Let $S$ be a compact orientable surface and $\ell$ be a nonsingular oriented line field on $\intr(S)$ satisfying $\wind_\ell(b) < 0$ for every boundary component $b$ of $S$.
If $S$ has genus $\geq 5$ and has at most two boundary components with even winding number, then there exists a collection $\mathcal{C}$ of winding number zero curves and a collection $\mathcal{A}$ of disjoint curves on $S$ where the complementary regions of $\mathcal{A}$ in $S$ each have genus $\geq 1$, such that given any $f \in \PMod(\intr(S))[\ell]$, there exist curves $c_i \in \mathcal{C}$ and $d_j \in \mathcal{A}$ such that $f = \prod \tau_{c_i}^{n_i} \prod \tau_{d_j}^{m_j}$ for some integers $n_i, m_j$.
\end{cor}
\begin{proof}
Apply \Cref{thm:factorizationthm} to the auxiliary curve system in \Cref{lemma:disjointauxiliarycurvescriterion}, and omit $\tau_c$ for every auxiliary curve $c$ parallel to a boundary component.
\end{proof}

\section{Standardizing boundary orbits} \label{subsec:standardizeorbit}

In this section, we prove the following proposition.

\begin{prop}[ChatGPT-6 Astra] \label{prop:standardizeorbit} 
Let $(\alpha_-,\xi_+)$ be a strongly adapted bicontact package, and let $\gamma$ be a closed orbit of the Reeb flow $R_-$. Then $(\alpha_-,\xi_+)$ is homotopic, through strongly adapted bicontact packages where $\gamma$ is a closed orbit of the Reeb flow, to a strongly adapted bicontact package $(\widehat{\alpha}_-,\widehat{\xi}_+)$ where $\gamma$ is a standard Reeb orbit of the Reeb flow $\widehat{R}_-$.

In fact, the positive contact structure can be chosen to be constant throughout the homotopy, and the homotopy can be chosen to be supported in an arbitrarily small tubular neighborhood $\nu$ of $\gamma$.
Moreover, if $\gamma$ is a negative boundary orbit of a $\partial$-strong Birkhoff section $S$ to $R_-$, then $S \cap (M \backslash \nu)$ can be extended into a $\partial$-strong Birkhoff section $\widehat{S}$ to $\widehat{R}_-$. In particular, there is a natural identification between $S$ and $\widehat{S}$ which preserves their monodromies and induced line fields from $\xi_+$ and $\widehat{\xi}_+$, up to homotopy.
\end{prop}

As discussed in \Cref{subsec:aidisclosure}, \Cref{prop:standardizeorbit} and its proof were explained to us by ChatGPT-6 Astra.

\subsection{Formula for blowing up}

For the purpose of computations in this section, we need an explicit formula for the blow-up operation (\Cref{constr:blowup}).

Let $\gamma$ be a closed orbit of a flow $\phi$ on $M$.
More concretely, pick a tubular neighborhood $\nu$ of $\gamma$ and pick coordinates $\nu \cong D^2 \times S^1_z$ where $\gamma = \{0\} \times S^1$. In Euclidean coordinates $(x,y)$ on $D^2$, we can write 
\begin{equation} \label{eq:vfbeforeblowup}
X = A \frac{\partial}{\partial x} + B \frac{\partial}{\partial y} + C \frac{\partial}{\partial z}
\end{equation}
for some smooth functions $A, B, C$. 
We refer to $A, B$ as the \textbf{transversal coefficients} and $C$ as the \textbf{tangential coefficient}.
Since $\gamma$ is a closed orbit, we have $A(0,0,z) = B(0,0,z) = 0$, so we can write $A = A_1 x + A_2 y$ and $B = B_1 x + B_2 y$ where $A_i, B_i$ are smooth functions on $D^2 \times S^1$.
More specifically, we can take
\begin{equation} \label{eq:blowupcoefficients}
\begin{split}
A_1(x,y,z) &=
\begin{cases}
\frac{A(x,0,z)-A(0,0,z)}{x} & \text{if $x \neq 0$} \\
\frac{\partial A}{\partial x}(0,0,z) & \text{if $x = 0$}
\end{cases} \\
A_2(x,y,z) &= 
\begin{cases}
\frac{A(x,y,z)-A(x,0,z)}{y} & \text{if $y \neq 0$} \\
\frac{\partial A}{\partial y}(x,0,z) & \text{if $y = 0$}
\end{cases} \\
B_1(x,y,z) &=
\begin{cases}
\frac{B(x,0,z)-B(0,0,z)}{x} & \text{if $x \neq 0$} \\
\frac{\partial B}{\partial x}(0,0,z) & \text{if $x = 0$}
\end{cases} \\
B_2(x,y,z) &= 
\begin{cases}
\frac{B(x,y,z)-B(x,0,z)}{y} & \text{if $y \neq 0$} \\
\frac{\partial B}{\partial y}(x,0,z) & \text{if $y = 0$}
\end{cases}
\end{split}.
\end{equation}
Converting to polar coordinates $(r,\theta)$ on $D^2$, we can then write 
\begin{align} \label{eq:normalblowup}
X =~& A \left( \cos \theta \frac{\partial}{\partial r} - \frac{1}{r} \sin \theta \frac{\partial}{\partial \theta} \right) + B \left( \sin \theta \frac{\partial}{\partial r} + \frac{1}{r} \cos \theta \frac{\partial}{\partial \theta} \right) + C \frac{\partial}{\partial z} \nonumber \\
\begin{split}
=~& r(A_1 \cos^2 \theta + A_2 \sin \theta \cos \theta + B_1 \cos \theta \sin \theta + B_2 \sin^2 \theta) \frac{\partial}{\partial r} \\
~& + (- A_1 \cos \theta \sin \theta - A_2 \sin^2 \theta + B_1 \cos^2 \theta + B_2 \sin \theta \cos \theta) \frac{\partial}{\partial \theta} + C \frac{\partial}{\partial z}
\end{split},
\end{align} %
at least on $(D^2 \backslash \{0\}) \times S^1_z = (0,1]_r \times S^1_\theta \times S^1_z$.

As in \Cref{constr:blowup}, the blown-up manifold $M_\gamma$ is obtained by removing $\gamma$ from $M$ and gluing in $[0,1]_r \times S^1_\theta \times S^1_z$ along $\nu \backslash \gamma \cong (0,1]_r \times S^1_\theta \times S^1_z$.
We define the \textbf{blown-up vector field} $X_\gamma$ by \Cref{eq:normalblowup},
where by smoothness of $A_i, B_i$ on $D^2 \times S^1$, we can extend the coefficients smoothly from $(0,1]_r \times S^1_\theta \times S^1_z$ to $[0,1]_r \times S^1_\theta \times S^1_z$.

We record an observation coming out of \Cref{eq:normalblowup}.

\begin{prop} \label{prop:C1closeimpliesblowupC0close}
Let $\phi$ and $\widetilde{\phi}$ be two flows on $M$, generated by vector fields $X$ and $\widetilde{X}$ respectively. Suppose $\gamma$ is a closed orbit of both $\phi$ and $\widetilde{\phi}$. 
If the transversal coefficients of $X$ and $\widetilde{X}$ are $C^1$-close, and their tangential coefficients are $C^0$-close, then the blown-up vector fields $X_\gamma$ and $\widetilde{X}_\gamma$ are $C^0$-close.
\end{prop}
\begin{proof}
Using the notation of \Cref{eq:vfbeforeblowup} and \Cref{eq:blowupcoefficients}, the coefficients $A$ and $\widetilde{A}$, and $B$ and $\widetilde{B}$ being $C^1$-close implies that the functions $A_i$ and $\widetilde{A}_i$, and $B_i$ and $\widetilde{B}_i$ are $C^0$-close. Together with the fact that $C$ and $\widetilde{C}$ are $C^0$-close, this in turn implies that the blown-up vector fields $X_\gamma$ and $\widetilde{X}_\gamma$ as defined by \Cref{eq:normalblowup} are $C^0$-close.
\end{proof}

\subsection{A correspondence between contact forms and functions}

The idea of the proof of \Cref{prop:standardizeorbit} is the following:
By \Cref{lemma:legendriancurvestdnbd}, there is a tubular neighborhood of $\gamma$ with coordinates $(x,y,z) \in [-\epsilon,\epsilon] \times [-\epsilon,\epsilon] \times S^1$ under which $\gamma = \{(0,0,z) \mid z \in S^1\}$ and $\xi_+ = \ker(dx+ydz)$. 
Up to reversing $x$ and $z$, we assume that $\gamma$ is oriented in the direction of increasing $z$. %

\Cref{lemma:negcontactformcorrfunction} below essentially says that the negative contact forms $\alpha_-$ on such a neighborhood whose Reeb vector fields are contained in $\xi_+$ can be treated as pairs of functions $(L,F)$ that satisfy two auxiliary inequalities.
Moreover, \Cref{rem_standardReeb_L_F} gives us a pair of functions $(L_\std, F_\std)$ that determines a standard Reeb orbit.
Thus the strategy is to modify the negative contact form $\alpha_-$ by modifying the corresponding functions $(L,F)$ until they are standard;
we will say more on how we do so in the subsequent subsections.

\begin{lemma} \label{lemma:negcontactformcorrfunction}
Consider $\nu = [-\epsilon,\epsilon]_x \times [-\epsilon,\epsilon]_y  \times S^1_z$ and $\xi_+ = \ker(dx+ydz)$ on $\nu$.
Suppose $\alpha_-$ is a negative contact form on $\nu$ whose Reeb vector field $R_-$ is contained in $\xi_+$. Then, up to shrinking $\epsilon$, there exist functions $L$ and $F$ such that
\begin{equation} \label{eq:negcontactformcorrfunction}
\alpha_- = - \frac{\partial L}{\partial y} dx + \left( L - y\frac{\partial L}{\partial y} \right) dz + dF.
\end{equation}
In this case, the functions $(L,F)$ must satisfy the two inequalities
\begin{equation} \label{eq:negcontactformcorrfunctioncondition1}
\frac{\partial^2 L}{\partial y^2} < 0
\end{equation}
and
\begin{equation} \label{eq:negcontactformcorrfunctioncondition2}
- L \frac{\partial^2 L}{\partial y^2} + W F > 0,
\end{equation}
where
\begin{equation} \label{eq:negcontactformcorrfunctionconditionW}
W := y \frac{\partial^2 L}{\partial y^2} \frac{\partial}{\partial x} + \left( \frac{\partial L}{\partial x} + \frac{\partial^2 L}{\partial y \partial z} - y \frac{\partial^2 L}{\partial x \partial y} \right) \frac{\partial}{\partial y} - \frac{\partial^2 L}{\partial y^2} \frac{\partial}{\partial z},
\end{equation}
and the Reeb vector field is given by
\begin{equation} \label{eq:negcontactformcorrfunctionR-}
R_- = \frac{1}{- L \frac{\partial^2 L}{\partial y^2} + W F} W.
\end{equation}

Conversely, given functions $L$ and $F$ satisfying \Cref{eq:negcontactformcorrfunctioncondition1} and \Cref{eq:negcontactformcorrfunctioncondition2}, the 1-form defined by \Cref{eq:negcontactformcorrfunction} is a negative contact form whose Reeb vector field is contained in $\xi_+$.
\end{lemma}
\begin{proof}
Suppose $\alpha_- = A(x,y,z) dx + B(x,y,z) dy + C(x,y,z) dz$ is a negative contact form on $\nu$ whose Reeb vector field $R_-$ is contained in $\xi_+$.
Define $F(x,y,z) = \int_0^y B(x,s,z) ds$, so that we can write 
\begin{equation} \label{eq:negcontactformcorrfunctionalpha}
\alpha_- - dF = A_1(x,y,z) dx + C_1(x,y,z) dz.
\end{equation}
We then define $L(x,y,z) = C_1(x,y,z) - yA_1(x,y,z)$.

We verify \Cref{eq:negcontactformcorrfunction}:
Taking the exterior derivative of \Cref{eq:negcontactformcorrfunctionalpha}, we have
$$d\alpha_- = - \frac{\partial A_1}{\partial y} dx \wedge dy + \left( -\frac{\partial A_1}{\partial z} + \frac{\partial C_1}{\partial x} \right) dx \wedge dz + \frac{\partial C_1}{\partial y} dy \wedge dz.$$
Since $R_-$ spans the kernel of $d\alpha_-$ and lies in $\xi_+ = \ker(dx+ydz)$, we have
$$0 = d\alpha_- \wedge (dx+ydz) = \left( \frac{\partial C_1}{\partial y} - y\frac{\partial A_1}{\partial y} \right) dx \wedge dy \wedge dz.$$
Thus 
$$\frac{\partial L}{\partial y} = \frac{\partial C_1}{\partial y} - y\frac{\partial A_1}{\partial y} - A_1 = -A_1$$
$$L - y\frac{\partial L}{\partial y} = C_1 - y\frac{\partial C_1}{\partial y} + y^2 \frac{\partial A_1}{\partial y} = C_1$$
which together with \Cref{eq:negcontactformcorrfunctionalpha} verifies \Cref{eq:negcontactformcorrfunction}.

Now starting with \Cref{eq:negcontactformcorrfunction}, we will verify \Cref{eq:negcontactformcorrfunctioncondition1}, \Cref{eq:negcontactformcorrfunctioncondition2}, and \Cref{eq:negcontactformcorrfunctionR-}:
We first take the exterior derivative of \Cref{eq:negcontactformcorrfunction} to get that
\begin{align*}
d\alpha_- &= \frac{\partial^2 L}{\partial y^2} dx \wedge dy + \left( \frac{\partial^2 L}{\partial y \partial z} + \frac{\partial L}{\partial x} - y \frac{\partial^2 L}{\partial x \partial y} \right) dx \wedge dz - y \frac{\partial^2 L}{\partial y^2} dy \wedge dz \\
&= - \iota_W (dx \wedge dy \wedge dz)
\end{align*}
Thus $W$ is a multiple of $R_-$.
More specifically, we can compute 
$$\alpha_-(W) = \left( - \frac{\partial L}{\partial y} dx + \left( L - y\frac{\partial L}{\partial y} \right) dz \right) (W) + dF(W) = - L \frac{\partial^2 L}{\partial y^2} + W F$$
so $R_- = \frac{1}{\alpha_-(W)} W = \frac{1}{- L \frac{\partial^2 L}{\partial y^2} + W F} W$.
This verifies \Cref{eq:negcontactformcorrfunctionR-}.

We also have
\begin{align} \label{eq:negcontactformcorrfunctionconditionalphawedgedalpha}
\alpha_- \wedge d\alpha_- &= - \alpha_-(W) (dx \wedge dy \wedge dz) \nonumber \\
&= - \left( - L \frac{\partial^2 L}{\partial y^2} + W F \right) dx \wedge dy \wedge dz
\end{align}
So \Cref{eq:negcontactformcorrfunctioncondition2} follows from the fact that $\alpha_-$ is a negative contact form.
This in turn implies that $W$ is a positive multiple of $R_-$.

Now, since $\gamma$ is a closed orbit of $R_-$, we have that $R_-$ is a positive multiple of $\frac{\partial}{\partial z}$ along $\gamma$. Thus up to shrinking $\epsilon$, we can assume that the $\frac{\partial}{\partial z}$ coefficient of $R_-$ is positive.
Meanwhile, the $\frac{\partial}{\partial z}$ coefficient of $W$ is $-\frac{\partial^2 L}{\partial y^2}$.
From this we have \Cref{eq:negcontactformcorrfunctioncondition1}.

Conversely, given functions $L$ and $F$ satisfying \Cref{eq:negcontactformcorrfunctioncondition1} and \Cref{eq:negcontactformcorrfunctioncondition2}, the above computation shows that the 1-form defined by \Cref{eq:negcontactformcorrfunction} satisfies \Cref{eq:negcontactformcorrfunctionconditionalphawedgedalpha}, hence is a negative contact form, and satisfies \Cref{eq:negcontactformcorrfunctionR-}.
Since $(dx+ydz)(W) = 0$, $R_-$ is contained in $\xi_+ = \ker(dx+ydz)$.
\end{proof}

\begin{rmk}[The functions $(L_{\mathrm{std}},F_{\mathrm{std}})$ giving a standard Reeb orbit]\label{rem_standardReeb_L_F}
Consider the functions
\begin{align}
L_{\mathrm{std}}(x,y,z) &= \Lambda  + \frac{1}{2} \omega_0^2 x^2 - \frac{1}{2} y^2 \label{eq:stdL} \\
F_{\mathrm{std}}(x,y,z) &= - \frac{1}{2} xy \label{eq:stdF}
\end{align}
where $\Lambda$ and $\omega_0$ are positive rational numbers. 
Then $\alpha_-$ as defined by \Cref{eq:negcontactformcorrfunction} is
$$\alpha_- = \frac{1}{2} (y dx - x dy) + (\Lambda + \frac{1}{2} \omega_0^2 x^2 + \frac{1}{2} y^2) dz.$$
Under the coordinate change
$x = \sqrt{\frac{2}{\omega_0}} u$ and $y = \sqrt{2\omega_0} v$, 
we have
\begin{align*}
\xi_+ &= \ker(du + \omega_0 v dz) \\
\alpha_- &= - u dv + v du + (\Lambda + \omega_0 (u^2 + v^2)) dz
\end{align*}
which describes the bicontact package near a standard Reeb orbit (\Cref{defn:stdorbit}) with flow slope $\omega_0$ and Legendrian slope $0$.
\end{rmk}

\subsection{Initial form of the functions}

The fact that \Cref{eq:stdL} and \Cref{eq:stdF} are both quadratic polynomials suggests one to modify $(L,F)$ by expanding out its Taylor series and modifying the terms one-by-one.
The following lemma does this for the constant and linear terms, and most of the quadratic terms.

\begin{lemma} \label{lemma:negcontactformcorrfunctiontaylorexpansion}
Let $(\alpha_-,\xi_+)$ be a strongly adapted bicontact package, and let $\gamma$ be a closed orbit of the Reeb flow $R_-$.
Then there exists a tubular neighborhood $\nu$ of $\gamma$ with coordinates $(x,y,z) \in [-\epsilon,\epsilon] \times [-\epsilon,\epsilon] \times S^1$ under which $\gamma = \{(0,0,z) \mid z \in S^1 \}$, $\xi_+ = \ker(dx+ydz)$, and $\alpha_-$ can be written as \Cref{eq:negcontactformcorrfunction} for some functions $(L,F)$, such that along $\gamma$,
\begin{enumerate}
    \setcounter{enumi}{-1}
    \item $F = 0$, \label{item:negcontactformcorrfunctiontaylorexpansionF=0}
    \item $L$ is constant, \label{item:negcontactformcorrfunctiontaylorexpansionLconstant}
    \item $\frac{\partial L}{\partial x} = \frac{\partial L}{\partial y} = 0$, \label{item:negcontactformcorrfunctiontaylorexpansionLfirstorder}
    \item $\frac{\partial^2 L}{\partial y^2} = -1$, and \label{item:negcontactformcorrfunctiontaylorexpansionLyy=-1}
    \item $\frac{\partial^2 L}{\partial x \partial y} = 0$. \label{item:negcontactformcorrfunctiontaylorexpansionLxy=0}
\end{enumerate}

Items (\ref{item:negcontactformcorrfunctiontaylorexpansionLconstant})--(\ref{item:negcontactformcorrfunctiontaylorexpansionLxy=0}) can be rephrased in the following way: For each $z$, we expand $L(x,y,z)$ into its Taylor series in $(x,y)$ centered at $(0,0)$ up to degree 2 terms, i.e., we write 
$$L(x,y,z) = \sum_{i+j \leq 2} \frac{1}{i! j!} \frac{\partial^{i+j} L}{\partial x^i \partial y^j}(0,0,z) x^i y^j + E(x,y,z)$$
where $\frac{\partial^{i+j} E}{\partial x^i\partial y^j}(x,y,z) = O((|x|+|y|)^{3-i-j})$, for all $i+j \leq 3$, uniformly in $z$. 
The content of items (\ref{item:negcontactformcorrfunctiontaylorexpansionLconstant})--(\ref{item:negcontactformcorrfunctiontaylorexpansionLxy=0}) is that in fact we have
\begin{equation} \label{eq:negcontactformcorrfunctiontaylorexpansion}
L(x,y,z) = \Lambda + \frac{1}{2} Q(z) x^2 - \frac{1}{2} y^2 + E(x,y,z)
\end{equation}
where $\Lambda$ is a constant and $Q$ is a function in $z$.
\end{lemma}
\begin{proof}
By \Cref{lemma:legendriancurvestdnbd} and \Cref{lemma:negcontactformcorrfunction}, there exists a tubular neighborhood $\nu \cong [-\epsilon,\epsilon]_x \times [-\epsilon,\epsilon]_y \times S^1_z$ under which $\gamma = \{(0,0,z) \mid z \in S^1 \}$, $\xi_+ = \ker(dx+ydz)$, and $\alpha_-$ can be written as \Cref{eq:negcontactformcorrfunction} for some functions $(L,F)$. In fact, the choice of $F$ exhibited there already satisfies (\ref{item:negcontactformcorrfunctiontaylorexpansionF=0}). The goal of this proof is to show that, up to modifying the coordinate system on $\nu$, we can arrange for $L$ to satisfy the remaining items one-by-one. 

\textbf{Arranging for (\ref{item:negcontactformcorrfunctiontaylorexpansionLconstant}).} 
We define $\Lambda = \frac{1}{2\pi} \int_0^{2\pi} L(0,0,s) ds$, and set %
\begin{align*}
x_1(x,y,z) &= x \\
y_1(x,y,z) &= \frac{\Lambda}{L(0,0,z)}y \\
z_1(x,y,z) &= \frac{1}{\Lambda} \int_0^z L(0,0,s) ds.
\end{align*}
Then 
\begin{align*}
dx_1 &= dx \\
dy_1 &= \frac{\Lambda}{L(0,0,z)} dy - \frac{\partial L}{\partial z}(0,0,z) \frac{\Lambda}{L(0,0,z)^2}y dz \\
dz_1 &= \frac{L(0,0,z)}{\Lambda} dz.
\end{align*} %
In particular, we have 
$$\frac{\partial}{\partial y_1} = \frac{L(0,0,z)}{\Lambda} \frac{\partial}{\partial y}.$$
Thus if we set 
$$L_1(x,y,z) = \frac{\Lambda}{L(0,0,z)} L(x,y,z),$$
then we have 
\begin{itemize}
    \item $\xi_+ = \ker(dx+ydz) = \ker(dx_1+y_1dz_1)$, 
    \item $\gamma = \{(0,0,z) \mid z \in S^1\} = \{(0,0,z_1) \mid z_1 \in S^1\}$, 
    \item $\alpha_- = - \frac{\partial L}{\partial y} dx + \left( L - y\frac{\partial L}{\partial y} \right) dz + dF = - \frac{\partial L_1}{\partial y_1} dx_1 + \left( L_1 - y_1\frac{\partial L_1}{\partial y_1} \right) dz_1 + dF$, and 
    \item $L_1(0,0,z) = \Lambda$ is a constant.
\end{itemize}
In other words, we have arranged for (\ref{item:negcontactformcorrfunctiontaylorexpansionLconstant}) while maintaining (\ref{item:negcontactformcorrfunctiontaylorexpansionF=0}).

\textbf{Arranging for (\ref{item:negcontactformcorrfunctiontaylorexpansionLfirstorder}).} 
Since $\gamma$ is a closed orbit of $R_-$, by \Cref{eq:negcontactformcorrfunctionR-}, the vector field $W_1$ defined as in \Cref{eq:negcontactformcorrfunctionconditionW} must be a positive multiple of $\frac{\partial}{\partial z_1}$ along $\gamma = \{(0,0,z_1) \mid z_1 \in S^1\}$.
In other words, we have 
\begin{equation} \label{eq:negcontactformcorrfunctionremovelinearterm}
\frac{\partial L_1}{\partial x_1} = - \frac{\partial^2 L_1}{\partial y_1 \partial z_1}
\end{equation}
along $\gamma$.
With this in mind, we set 
$$L_2(x_1,y_1,z_1) = L_1(x_1,y_1,z_1) + \frac{\partial^2 L_1}{\partial y_1 \partial z_1}(0,0,z_1) x_1 - \frac{\partial L_1}{\partial y_1}(0,0,z_1) y_1$$
and 
$$F_2(x_1,y_1,z_1) = F(x_1,y_1,z_1) - \frac{\partial L_1}{\partial y_1}(0,0,z_1) x_1.$$ 
Then 
\begin{itemize}
    \item $\alpha_- = - \frac{\partial L_1}{\partial y_1} dx_1 + \left( L_1 - y_1\frac{\partial L_1}{\partial y_1} \right) dz_1 + dF = - \frac{\partial L_2}{\partial y_1} dx_1 + \left( L_2 - y_1\frac{\partial L_2}{\partial y_1} \right) dz_1 + dF_2$,
\end{itemize} 
and along $\gamma = \{(0,0,z_1) \mid z_1 \in S^1\}$, we have 
\begin{itemize}
    \item $F_2 = F = 0$,
    \item $L_2 = L_1$ is constant, and 
    \item $\frac{\partial L_2}{\partial x_1} = \frac{\partial L_2}{\partial y_1} = 0$, where we used \Cref{eq:negcontactformcorrfunctionremovelinearterm}.
\end{itemize}
In other words, we now have (\ref{item:negcontactformcorrfunctiontaylorexpansionF=0})--(\ref{item:negcontactformcorrfunctiontaylorexpansionLfirstorder}).

\textbf{Arranging for (\ref{item:negcontactformcorrfunctiontaylorexpansionLyy=-1}).} 
Thanks to \Cref{eq:negcontactformcorrfunctioncondition1}, we also have $\frac{\partial^2 L_2}{\partial y_1^2} <0$, so we can set $a(z_1) := \left( -\frac{\partial^2 L_2}{\partial y_1^2}(0,0,z_1) \right)^{\frac{1}{2}}$, and define
\begin{align*}
x_3(x_1,y_1,z_1) &= a(z_1) x_1 \\
y_3(x_1,y_1,z_1) &= - a'(z_1) x_1 + a(z_1) y_1 \\
z_3(x_1,y_1,z_1) &= z_1. 
\end{align*}
Then
\begin{align*}
dx_3 &= a(z_1) dx_1 + a'(z_1) x_1 dz_1 \\
dy_3 &= - a'(z_1) dx_1 + a(z_1) dy_1 + (- a''(z_1) x_1 + a'(z_1) y_1) dz_1 \\
dz_3 &= dz_1.
\end{align*}
In particular, we have 
$$\frac{\partial}{\partial y_3} = \frac{1}{a(z_1)} \frac{\partial}{\partial y_1}.$$
Thus we have 
\begin{itemize}
    \item $\xi_+ = \ker(dx_1+y_1dz_1) = \ker(a(z_1)(dx_1+y_1dz_1)) = \ker(dx_3+y_3dz_3)$, 
    \item $\gamma = \{(0,0,z_1) \mid z_1 \in S^1\} = \{(0,0,z_3) \mid z_3 \in S^1\}$, 
    \item $\alpha_- = - \frac{\partial L_2}{\partial y_1} dx_1 + \left( L_2 - y_1 \frac{\partial L_2}{\partial y_1} \right) dz_1 + dF_2 = - \frac{\partial L_2}{\partial y_3} dx_3 + \left( L_2 - y_3 \frac{\partial L_2}{\partial y_3} \right) dz_3 + dF_2$.
\end{itemize}
Moreover, along $\gamma$, we have 
\begin{align*}
dx_3 &= a(z_1) dx_1 \\
dy_3 &= - a'(z_1) dx_1 + a(z_1) dy_1 \\
dz_3 &= dz_1
\end{align*}
thus
\begin{align*}
\frac{\partial}{\partial x_3} &=  \frac{1}{a(z_1)} \frac{\partial}{\partial x_1} + \frac{a'(z_1)}{a(z_1)^2} \frac{\partial}{\partial y_1} \\
\frac{\partial}{\partial y_3} &= \frac{1}{a(z_1)} \frac{\partial}{\partial y_1} \\
\frac{\partial}{\partial z_3} &= \frac{\partial}{\partial z_1}.
\end{align*}
From this, we can compute that along $\gamma$, 
\begin{itemize}
    \item $\frac{\partial L_2}{\partial x_3} = \frac{1}{a(z_1)} \frac{\partial L_2}{\partial x_1} + \frac{a'(z_1)}{a(z_1)^2} \frac{\partial L_2}{\partial y_1} = 0$, 
    \item $\frac{\partial L_2}{\partial y_3} = \frac{1}{a(z_1)} \frac{\partial L_2}{\partial y_1} = 0$, and 
    \item $\frac{\partial^2 L_2}{\partial y_3^2} = \frac{1}{a(z_1)^2} \frac{\partial^2 L_2}{\partial y_1^2} = -1$ by definition of $a(z_1)$. 
\end{itemize}
In other words, we now have (\ref{item:negcontactformcorrfunctiontaylorexpansionF=0})--(\ref{item:negcontactformcorrfunctiontaylorexpansionLyy=-1}).

\textbf{Arranging for (\ref{item:negcontactformcorrfunctiontaylorexpansionLxy=0}).} 
We let $b(z_3) = \frac{\partial^2 L_2}{\partial x_3 \partial y_3}(0,0,z_3)$, and set 
$$L_4(x_3,y_3,z_3) = L_2(x_3,y_3,z_3) - b(z_3) x_3 y_3 + \frac{1}{2} b'(z_3) x_3^2,$$
and 
$$F_4(x_3,y_3,z_3) = F_2(x_3,y_3,z_3) - \frac{1}{2} b(z_3) x_3^2.$$ %
Then 
\begin{itemize}
    \item $\alpha_- = - \frac{\partial L_2}{\partial y_3} dx_3 + \left( L_2 - y_3 \frac{\partial L_2}{\partial y_3} \right) dz_3 + dF_2 = - \frac{\partial L_4}{\partial y_3} dx_3 + \left( L_4 - y_3 \frac{\partial L_4}{\partial y_3} \right) dz_3 + dF_4$,
\end{itemize} 
and along $\gamma = \{(0,0,z_3) \mid z_3 \in S^1\}$, we have that
\begin{itemize}
    \item $F_4 = F_2 = 0$,
    \item $L_4 = L_2$ is constant, 
    \item $\frac{\partial L_4}{\partial x_3} = \frac{\partial L_2}{\partial x_3} = 0$, 
    \item $\frac{\partial L_4}{\partial y_3} = \frac{\partial L_2}{\partial y_3} = 0$, 
    \item $\frac{\partial^2 L_4}{\partial y_3^2} = \frac{\partial^2 L_2}{\partial y_3^2} = -1$, and
    \item $\frac{\partial^2 L_4}{\partial x_3 \partial y_3} = \frac{\partial^2 L_2}{\partial x_3 \partial y_3} - b(z_3) = 0$ by definition of $b(z_3)$.
\end{itemize}
In other words, we now have (\ref{item:negcontactformcorrfunctiontaylorexpansionF=0})--(\ref{item:negcontactformcorrfunctiontaylorexpansionLxy=0}).
\end{proof}

\subsection{Modifying the functions}

Picking up from \Cref{lemma:negcontactformcorrfunctiontaylorexpansion}, it remains to modify $(L,F)$ in the following 3 steps:
\begin{enumerate}[label=(\Roman*)]
    \item Eliminating the cubic term $E(x,y,z)$ from $L$. %
    \item Interpolating the quadratic coefficient $Q(z)$ in $L$ to a constant $\omega_0^2$, where $\omega_0$ is a positive rational number.
    \item Interpolating $F$ to \Cref{eq:stdF}.
\end{enumerate}

\begin{constr}[Eliminating the cubic term from $L$] \label{constr:removecubicterm}
Consider $\nu = [-\epsilon,\epsilon]_x \times [-\epsilon,\epsilon]_y  \times S^1_z$. Suppose $\alpha_-$ is a negative contact form defined by \Cref{eq:negcontactformcorrfunction} for some $(L,F)$, where $L$ is of the form of \Cref{eq:negcontactformcorrfunctiontaylorexpansion}.

Let $\rho\colon \mathbb{R}_{\geq 0} \to [0,1]$ be a function where $\rho(t) = 0$ for $t$ close to $0$ and $\rho(t) = 1$ for $t \geq 1$.
For $s \in [0,1]$ and $\delta > 0$, we define
$$L_s(x,y,z) = \Lambda + \frac{1}{2} Q(z) x^2 - \frac{1}{2} y^2 + \left( 1- s \left( 1 - \rho \left( \frac{r}{\delta} \right) \right) \right) E(x,y,z)$$
where $r = \sqrt{x^2+y^2}$.
Then $L_0 = L$ and $L_1 = \Lambda + \frac{1}{2} Q(z) x^2 - \frac{1}{2} y^2$ in a smaller tubular neighborhood of $\{0\} \times S^1$.

We claim that for sufficiently small $\delta$, $(L_s,F)$ satisfies \Cref{eq:negcontactformcorrfunctioncondition1} and \Cref{eq:negcontactformcorrfunctioncondition2} for all $s \in [0,1]$.
Indeed, 
\begin{align*}
\left| \frac{\partial L_s}{\partial x} - \frac{\partial L}{\partial x} \right| &= \left| \frac{\partial}{\partial x} \left( s \left( \rho \left( \frac{r}{\delta} \right) - 1 \right) E(x,y,z) \right) \right| \\
&= \left| \frac{s}{\delta}\rho' \left( \frac{r}{\delta} \right) \frac{\partial r}{\partial x} E(x,y,z) + s \left( \rho \left( \frac{r}{\delta} \right) - 1 \right) \frac{\partial E}{\partial x}(x,y,z) \right| \\
&\leq \frac{1}{\delta} \left| \rho' \left( \frac{r}{\delta} \right) \right| |E(x,y,z)| + \left| \rho \left( \frac{r}{\delta} \right) - 1 \right| \left| \frac{\partial E}{\partial x}(x,y,z) \right| \\
&= O(\delta^{-1}) \cdot O(1) \cdot O(\delta^3) + O(1) \cdot O(\delta^2) = O(\delta^2)
\end{align*}
since $E(x,y,z) = O(r^3)$ and $\frac{\partial E}{\partial x}(x,y,z) = O(r^2)$, and $\rho(\frac{r}{\delta}) = 1$ and $\rho'(\frac{r}{\delta}) = 0$ unless $r \leq \delta$.
More generally, similar computations show that 
\begin{equation} \label{eq:removecubictermestimate}
\left| \frac{\partial^{i+j+k} L_s}{\partial x^i \partial y^j \partial z^k} - \frac{\partial^{i+j+k} L}{\partial x^i \partial y^j \partial z^k} \right| = O(\delta^{3-i-j})
\end{equation}
for each $i+j \leq 3$.
In particular, $\left| \frac{\partial^2 L_s}{\partial y^2} - \frac{\partial^2 L}{\partial y^2} \right| = O(\delta)$.
Hence \Cref{eq:negcontactformcorrfunctioncondition1} for $L_s$ follows from those for $L$, for small enough $\delta$.
Also one can check that 
$$\left| \left( -L\frac{\partial^2 L}{\partial y^2}+WF \right) - \left( -L_s\frac{\partial^2 L_s}{\partial y^2}+W_sF \right) \right| \leq \left|L\frac{\partial^2 L}{\partial y^2} - L_s\frac{\partial^2 L_s}{\partial y^2} \right| + \left| (W-W_s) F \right| = O(\delta),$$ 
where $W$ and $W_s$ are defined by \Cref{eq:negcontactformcorrfunctionconditionW} with $L$ and $L_s$, respectively.
Hence \Cref{eq:negcontactformcorrfunctioncondition2} for $L_s$ follows from those for $L$, for small enough $\delta$.

We can thus define a family of negative contact forms $(\alpha_-)_s$ by applying \Cref{lemma:negcontactformcorrfunction} to $(L_s,F)$. By definition, $((\alpha_-)_s,\xi_+)$ is strongly adapted for every $s$.
Also, we have $\frac{\partial L_s}{\partial x}(0,0,z) = \frac{\partial^2 L_s}{\partial y \partial z}(0,0,z) = 0$ for every $s$.
Hence we have a homotopy of strongly adapted bicontact packages through strongly adapted bicontact packages where $\{0\} \times S^1$ is a closed orbit of the Reeb flow $(R_-)_s$.

Finally, we claim that any $\partial$-strong Birkhoff section $S$ for $R_-$ is automatically a Birkhoff section for $(R_-)_s$, for every $s$, provided that $\delta$ is small enough.
Indeed, from \Cref{eq:negcontactformcorrfunctionconditionW} and \Cref{eq:removecubictermestimate}, one can check that the transversal coefficients of $W_s$ have arbitrarily small variation in $C^1$ norm, while their tangential coefficient has arbitrarily small variation in $C^0$ norm, as $s$ varies over $[0,1]$.
Thus by \Cref{prop:C1closeimpliesblowupC0close}, the blown-up vector fields $(W_s)_\gamma$, which are positive multiples of the blown-up Reeb vector fields, have arbitrarily small variation in $C^0$ norm, which implies the claim. %
\end{constr}

\begin{constr}[Interpolating $Q$] \label{constr:interpolateQ}
Consider $\nu = [-\epsilon,\epsilon]_x \times [-\epsilon,\epsilon]_y  \times S^1_z$. Suppose $\alpha_-$ is a negative contact form defined by \Cref{eq:negcontactformcorrfunction} for some $(L,F)$ satisfying the items of \Cref{lemma:negcontactformcorrfunctiontaylorexpansion} and for which $L$ has no cubic term, i.e., such that
$$L(x,y,z) = \Lambda + \frac{1}{2} Q(z) x^2 - \frac{1}{2} y^2$$ for some constant $\Lambda$ and function $Q$ in $z$.
Let $Q(z,t)$ be a smooth function on $S^1 \times [0,1]$ where $Q(z,1) = Q(z)$.

As in the previous construction, we pick $\rho\colon \mathbb{R}_{\geq 0} \to [0,1]$ a function where $\rho(t) = 0$ for $t$ close to $0$ and $\rho(t) = 1$ for $t \geq 1$.
For $s \in [0,1]$ and $\epsilon, \delta > 0$, we define
$$L_s(x,y,z) = \Lambda + \frac{1}{2} \left( s Q(z) + (1-s) Q \left( z,\rho \left( \left( \frac{r}{\delta} \right)^\epsilon \right) \right) \right) x^2 - \frac{1}{2} y^2$$
where $r = \sqrt{x^2+y^2}$.
Then $L_1 = L$ and $L_0 = \Lambda + \frac{1}{2} Q(z,0) x^2 - \frac{1}{2} y^2$ in a smaller tubular neighborhood of $\{0\} \times \{0\} \times S^1$.

We claim that for sufficiently small $\epsilon, \delta$, $(L_s,F)$ satisfies \Cref{eq:negcontactformcorrfunctioncondition1} and \Cref{eq:negcontactformcorrfunctioncondition2} for all $s \in [0,1]$.
Indeed, let us denote $\widehat{Q}(x,y,z) = Q \left( z,\rho \left( \left( \frac{r}{\delta} \right)^\epsilon \right) \right)$, and compute
\begin{align*} 
\left| \frac{\partial \widehat{Q}}{\partial x} \right| &= \left| \frac{\partial Q}{\partial t} \rho' \left( \left( \frac{r}{\delta} \right)^\epsilon \right) \epsilon \frac{r^{\epsilon-1}}{\delta^\epsilon} \frac{\partial r}{\partial x} \right| \\
&\leq \left| \frac{\partial Q}{\partial t} \right| \left| \rho' \left( \left( \frac{r}{\delta} \right)^\epsilon \right) \left( \frac{r}{\delta} \right)^\epsilon \right| \frac{\epsilon}{r} \left| \frac{\partial r}{\partial x} \right| \\
&= 
\begin{cases}
0 & \text{if $r > \delta$} \\
O(1) \cdot O(1) \cdot O(\epsilon r^{-1}) \cdot O(1) = O(\epsilon r^{-1}) & \text{if $r \leq \delta$}
\end{cases} 
\end{align*}
More generally, similar computations show that 
\begin{equation} \label{eq:interpolateQestimate}
\frac{\partial^{i+j+k} \widehat{Q}}{\partial x^i \partial y^j \partial z^k} = 
\begin{cases}
0 & \text{if $r > \delta$} \\
O(\epsilon r^{-i-j}) & \text{if $r \leq \delta$}
\end{cases}
\end{equation}
for each $i,j,k$ with $i+j \geq 1$. %

Thus
\begin{align*}
\frac{\partial^{i+j+k} L_s}{\partial x^i \partial y^j \partial z^k} - \frac{\partial^{i+j+k} L}{\partial x^i \partial y^j \partial z^k} &= \frac{\partial^{i+j+k}}{\partial x^i \partial y^j \partial z^k} \left( \frac{1-s}{2} \left( \widehat{Q}(x,y,z) - Q(z) \right) x^2 \right) \\
&= 
\begin{cases}
O(\delta^{2-i}) & \text{if $j = 0$ and $0 \leq i \leq 2$} \\
O(\epsilon \delta^{2-i-j}) & \text{if $j > 0$ and $1 \leq i+j \leq 2$}
\end{cases}.
\end{align*} %
In particular, $\left| \frac{\partial^2 L_s}{\partial y^2} - \frac{\partial^2 L}{\partial y^2} \right| = O(\epsilon)$.
Hence \Cref{eq:negcontactformcorrfunctioncondition1} for $L_s$ follows from those for $L$, for small enough $\epsilon$.
Also, from \Cref{eq:negcontactformcorrfunctionconditionW}, one can check that 
$$\left| \left( -L\frac{\partial^2 L}{\partial y^2}+WF \right) - \left( -L_s\frac{\partial^2 L_s}{\partial y^2}+W_sF \right) \right| \leq \left|L\frac{\partial^2 L}{\partial y^2} - L_s\frac{\partial^2 L_s}{\partial y^2} \right| + \left| (W-W_s) F \right| = O(\delta) + O(\epsilon).$$ 
Hence \Cref{eq:negcontactformcorrfunctioncondition2} for $L_s$ follows from those for $L$, for small enough $\epsilon, \delta$.

We can thus define a family of negative contact forms $(\alpha_-)_s$ by applying \Cref{lemma:negcontactformcorrfunction} to $(L_s,F)$. By definition, $((\alpha_-)_s,\xi_+)$ is strongly adapted for every $s$.
Also, we have $\frac{\partial L_s}{\partial x}(0,0,z) = \frac{\partial^2 L_s}{\partial y \partial z}(0,0,z) = 0$ for every $s$.
Hence, as in the preceding construction, we have a homotopy through strongly adapted bicontact packages where $\{0\} \times\{0\} \times S^1$ is a closed orbit of the Reeb flow $(R_-)_s$.

Now suppose we are given a $\partial$-strong Birkhoff section $S$ for the original Reeb flow $R_-$. We will show that we can extend $S \cap \partial \nu$ into a $\partial$-strong Birkhoff section for the modified Reeb flow $(R_-)_0$, provided that a slope condition is satisfied.

To motivate the condition, we compute from \Cref{eq:negcontactformcorrfunctionconditionW} that
\begin{align*}
W_0 &= - y \left( 1 - \frac{1}{2} \frac{\partial^2 \widehat{Q}}{\partial y^2} x^2 \right) \frac{\partial}{\partial x} + \left( \frac{1}{2} \frac{\partial \widehat{Q}}{\partial x} x^2 + \widehat{Q} x + \frac{1}{2} \frac{\partial^2 \widehat{Q}}{\partial y \partial z} x^2 \right. \\
&~\left. - y \left( \frac{1}{2} \frac{\partial^2 \widehat{Q}}{\partial x \partial y} x^2 + \frac{\partial \widehat{Q}}{\partial y} x \right) \right) \frac{\partial}{\partial y} + \left( 1 - \frac{1}{2} \frac{\partial^2 \widehat{Q}}{\partial y^2} x^2 \right) \frac{\partial}{\partial z} \\
&= -y \frac{\partial}{\partial x} + \widehat{Q} x \frac{\partial}{\partial y} + \frac{\partial}{\partial z} \\
&~+ \frac{1}{2} \left( \frac{\partial^2 \widehat{Q}}{\partial y^2} x^2y \frac{\partial}{\partial x} + \left( \frac{\partial \widehat{Q}}{\partial x} x^2 + \frac{\partial^2 \widehat{Q}}{\partial y \partial z} x^2 - \frac{\partial^2 \widehat{Q}}{\partial x \partial y} x^2y - 2 \frac{\partial \widehat{Q}}{\partial y} xy \right) \frac{\partial}{\partial y} - \frac{\partial^2 \widehat{Q}}{\partial y^2} x^2 \frac{\partial}{\partial z} \right).
\end{align*}
We set
\begin{align*}
\widetilde{W_0} &:=  -y \frac{\partial}{\partial x} + \widehat{Q} x \frac{\partial}{\partial y} + \frac{\partial}{\partial z}, \text{ and }\\
U &:= W_0 - \widetilde{W_0}.
\end{align*}

From \Cref{eq:normalblowup} and \Cref{eq:interpolateQestimate}, we can compute that
\begin{align} \label{eq:interpolateQblowupestimate1}
U_\gamma &= \frac{1}{2} r \left( \frac{\partial^2 \widehat{Q}}{\partial y^2} xy \cos^2 \theta + \left( \frac{\partial \widehat{Q}}{\partial x} x + \frac{\partial^2 \widehat{Q}}{\partial y \partial z} x - \frac{\partial^2 \widehat{Q}}{\partial x \partial y} xy - 2 \frac{\partial \widehat{Q}}{\partial y} y \right) \cos \theta \sin \theta \right) \frac{\partial}{\partial r} \nonumber \\
&~+ \frac{1}{2} \left(- \frac{\partial^2 \widehat{Q}}{\partial y^2} xy \cos \theta \sin \theta + \left( \frac{\partial \widehat{Q}}{\partial x} x + \frac{\partial^2 \widehat{Q}}{\partial y \partial z} x - \frac{\partial^2 \widehat{Q}}{\partial x \partial y} xy - 2 \frac{\partial \widehat{Q}}{\partial y} y \right) \cos^2 \theta \right) \frac{\partial}{\partial \theta} \nonumber \\
&~- \frac{1}{2} \frac{\partial^2 \widehat{Q}}{\partial y^2} x^2 \frac{\partial}{\partial z} \nonumber \\
&= O(\epsilon r) \frac{\partial}{\partial r} + O(\epsilon) \frac{\partial}{\partial \theta} + O(\epsilon) \frac{\partial}{\partial z}
\end{align}
In particular, up to shrinking $\epsilon$, we can make $U_\gamma$ arbitrarily $C^0$-small, thus the blow-up $(W_0)_\gamma$, which is a positive multiple of the blown-up Reeb flow $((R_-)_0)_\gamma$, is arbitrarily $C^0$-close to the blow-up $(\widetilde{W}_0)_\gamma$.

Meanwhile, from \Cref{eq:normalblowup}, we can compute that
\begin{align} \label{eq:interpolateQblowupestimate2}
(\widetilde{W_0})_\gamma &= r (-1 + \widehat{Q}) \sin \theta \cos \theta \frac{\partial}{\partial r} + \left(\sin^2 \theta + \widehat{Q} \cos^2 \theta \right) \frac{\partial}{\partial \theta} + \frac{\partial}{\partial z} \nonumber \\
&= r \left( -1 + Q(z,\rho \left( \left( \frac{r}{\delta} \right)^\epsilon \right)) \right) \sin \theta \cos \theta \frac{\partial}{\partial r} \nonumber \\
&~+ \left(\sin^2 \theta + Q(z,\rho \left( \left( \frac{r}{\delta} \right)^\epsilon \right)) \cos^2 \theta \right) \frac{\partial}{\partial \theta} + \frac{\partial}{\partial z} \nonumber \\
&= O(r) \frac{\partial}{\partial r} + \left(\sin^2 \theta + Q(z,\rho \left( \left( \frac{r}{\delta} \right)^\epsilon \right)) \cos^2 \theta \right) \frac{\partial}{\partial \theta} + \frac{\partial}{\partial z}.
\end{align}
In particular, up to shrinking $\delta$, we can make the first term arbitrarily $C^0$-small.
This allows us to focus on
\[
\widetilde{W_\gamma}:= \left(\sin^2 \theta + Q(z,\rho \left( \left( \frac{r}{\delta} \right)^\epsilon \right)) \cos^2 \theta \right) \frac{\partial}{\partial \theta} + \frac{\partial}{\partial z}.
\]
In order to analyze $\widetilde{W_\gamma}$, we make the following definition:
Given a function $q(z)$ in $z$, we denote by $f(q)$ the slope of the foliation on $S^1_\theta \times S^1_z$ by flow lines of $\widetilde{W}(q) := (\sin^2 \theta + q(z) \cos^2 \theta) \frac{\partial}{\partial \theta} + \frac{\partial}{\partial z}$. %

The following lemma records the properties of $f(q)$ that we need.
\begin{lemma} \label{lemma:flowlineslope}
We identify the slopes with $\mathbb{R} \cup \{\infty\}$ so that $\frac{\partial}{\partial z}$ is $0$ and $\frac{\partial}{\partial \theta}$ is $\infty$. Then:
\begin{enumerate}
    \item $f(q)$ depends continuously on the function $q$.
    \item If $q_1(z) \leq q_2(z)$, then $f(q_1) \leq f(q_2)$.
    \item If $q(z)$ is constant equal to $1$, then $f(q) = 1$.
    \item If $q(z)$ is a nonpositive constant $q_0 \leq 0$, then $f(q) = 0$.
\end{enumerate}
\end{lemma}
\begin{proof}

(1) follows from standard ODE continuous dependence, and continuity of translation number. %
(2) follows from standard ODE comparison.

For (3), if $q(z) = 1$, then $\widetilde W(q) = \frac{\partial}{\partial \theta} + \frac{\partial}{\partial z}$, so it is clear that $f(q) = 1$.

For (4), suppose $q(z) = q_0 \leq 0$. Note that along $\theta = \frac{\pi}{2}$, the vector field has positive $\frac{\partial}{\partial \theta}$ component, and along $\theta = \pi$, the vector field has nonpositive $\frac{\partial}{\partial \theta}$ component. So there must be a closed orbit of slope $\frac{\partial}{\partial z}$ within the annulus $\{ \theta \in [\frac{\pi}{2},\pi] \}$, see \Cref{fig:flowlineslope0}. %

\begin{figure}
    \centering
    \begingroup%
  \makeatletter%
  \providecommand\color[2][]{%
    \errmessage{(Inkscape) Color is used for the text in Inkscape, but the package 'color.sty' is not loaded}%
    \renewcommand\color[2][]{}%
  }%
  \providecommand\transparent[1]{%
    \errmessage{(Inkscape) Transparency is used (non-zero) for the text in Inkscape, but the package 'transparent.sty' is not loaded}%
    \renewcommand\transparent[1]{}%
  }%
  \providecommand\rotatebox[2]{#2}%
  \newcommand*\fsize{\dimexpr\f@size pt\relax}%
  \newcommand*\lineheight[1]{\fontsize{\fsize}{#1\fsize}\selectfont}%
  \ifx\svgwidth\undefined%
    \setlength{\unitlength}{85.88976162bp}%
    \ifx\svgscale\undefined%
      \relax%
    \else%
      \setlength{\unitlength}{\unitlength * \real{\svgscale}}%
    \fi%
  \else%
    \setlength{\unitlength}{\svgwidth}%
  \fi%
  \global\let\svgwidth\undefined%
  \global\let\svgscale\undefined%
  \makeatother%
  \begin{picture}(1,1.14210302)%
    \lineheight{1}%
    \setlength\tabcolsep{0pt}%
    \put(0,0){\includegraphics[width=\unitlength,page=1]{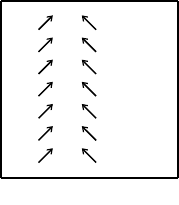}}%
    \put(0.21991432,0.02236741){\color[rgb]{0,0,0}\makebox(0,0)[lt]{\lineheight{1.25}\smash{\begin{tabular}[t]{l}$\frac{\pi}{2}$\end{tabular}}}}%
    \put(0.46588802,0.02235432){\color[rgb]{0,0,0}\makebox(0,0)[lt]{\lineheight{1.25}\smash{\begin{tabular}[t]{l}$\pi$\end{tabular}}}}%
  \end{picture}%
\endgroup%

    \caption{If $q(z)$ is a nonpositive constant $q_0 \leq 0$, then $f(q) = 0$.}
    \label{fig:flowlineslope0}
\end{figure}

\end{proof}

Now returning to the given Birkhoff section $S$, up to restricting to a smaller tubular neighborhood, we can first isotope $S$ along flow lines so that it intersects the tori $\{r=r_0\}$ in a multicurve $c \subset S^1_\theta \times S^1_z$ that is independent of $r_0$ and transverse to the meridional fibers $\{z=z_0\}$.
Since $(W_0)_\gamma$ is arbitrarily $C^0$-close to $\widetilde{W}_\gamma$, we can also assume that $c$ is transverse to $\widetilde{W}_\gamma$.

The slope condition we need for extending $S \cap \partial \nu$ into a Birkhoff section is for the slope $\sigma$ induced by $S$ to be different from $f(Q(\cdot,t))$ for every $t$. 

In this case, for every $t$, we can find a multicurve $c_t$ transverse to $\widetilde{W}(Q(\cdot,t))$, transverse to the meridional fibers $\{z=z_0\}$, and of slope $\sigma$. In fact, one can further make $c_t$ vary smoothly in $t$ and so that $c_t = c$ for $t$ close to $1$. We can then define a surface $A$ transverse to $\widetilde{W}_\gamma$ by taking the union of the curves $c_t$.
More precisely, parametrize $c_t(u) = (\theta(t,u),z(t,u))$, and set $A(r,u) = \left( r, \theta( \rho \left( \left( \frac{r}{\delta} \right)^\epsilon \right),u), z( \rho \left( \left( \frac{r}{\delta} \right)^\epsilon \right),u) \right)$.
This is a smooth extension of $S \cap \partial \nu$ into $\nu$. %

We can then compute that  
\begin{align*}
\frac{\partial A}{\partial r} \wedge \frac{\partial A}{\partial u} &= \left( \frac{\partial}{\partial r} + \rho' \left( \left( \frac{r}{\delta} \right)^\epsilon \right) \epsilon \left( \frac{r}{\delta} \right)^\epsilon \frac{1}{r} \left( \frac{\partial \theta}{\partial t} \frac{\partial}{\partial \theta} + \frac{\partial z}{\partial t} \frac{\partial}{\partial z} \right) \right) \wedge \left( \frac{\partial \theta}{\partial u} \frac{\partial}{\partial \theta} + \frac{\partial z}{\partial u} \frac{\partial}{\partial z} \right) \\
&= \frac{\partial \theta}{\partial u} \frac{\partial}{\partial r} \wedge \frac{\partial}{\partial \theta} + \frac{\partial z}{\partial u} \frac{\partial}{\partial r} \wedge \frac{\partial}{\partial z} + \rho' \left( \left( \frac{r}{\delta} \right)^\epsilon \right) \left( \frac{r}{\delta} \right)^\epsilon \frac{\epsilon}{r} \left( \frac{\partial \theta}{\partial t} \frac{\partial z}{\partial u} - \frac{\partial \theta}{\partial u} \frac{\partial z}{\partial t} \right) \frac{\partial}{\partial \theta} \wedge \frac{\partial}{\partial z} \\
&= O(1) \frac{\partial}{\partial r} \wedge \frac{\partial}{\partial \theta} + O(1) \frac{\partial}{\partial r} \wedge \frac{\partial}{\partial z} + O(\frac{\epsilon}{r}) \frac{\partial}{\partial \theta} \wedge \frac{\partial}{\partial z}
\end{align*} %

Thus, from \Cref{eq:interpolateQblowupestimate1} and \Cref{eq:interpolateQblowupestimate2},
\begin{align*}
\frac{\partial A}{\partial r} \wedge \frac{\partial A}{\partial u} \wedge ((W_0)_\gamma - \widetilde{W_\gamma}) &= \left( O(1) \frac{\partial}{\partial r} \wedge \frac{\partial}{\partial \theta} + O(1) \frac{\partial}{\partial r} \wedge \frac{\partial}{\partial z} + O(\frac{\epsilon}{r}) \frac{\partial}{\partial \theta} \wedge \frac{\partial}{\partial z} \right) \\
&~\wedge \left( (O(\epsilon r) + O(r)) \frac{\partial}{\partial r} + O(\epsilon) \frac{\partial}{\partial \theta} + O(\epsilon) \frac{\partial}{\partial z} \right) \\
&= O(\epsilon) \frac{\partial}{\partial r} \wedge \frac{\partial}{\partial \theta} \wedge \frac{\partial}{\partial z}
\end{align*}
can be made arbitrarily small.

Meanwhile,
\begin{align*}
\frac{\partial A}{\partial r} \wedge \frac{\partial A}{\partial u} \wedge \widetilde{W_\gamma} &= \left( \frac{\partial \theta}{\partial u} \frac{\partial}{\partial r} \wedge \frac{\partial}{\partial \theta} + \frac{\partial z}{\partial u} \frac{\partial}{\partial r} \wedge \frac{\partial}{\partial z} \right) \wedge \widetilde{W_\gamma} \\
&= c'_{\rho \left( \left( \frac{r}{\delta} \right)^\epsilon \right)} \wedge \widetilde{W} \left( Q \left( \cdot, \rho \left( \left( \frac{r}{\delta} \right)^\epsilon \right) \right) \right) \wedge \frac{\partial}{\partial r}
\end{align*}
is bounded away from $0$.
Thus $\frac{\partial A}{\partial r} \wedge \frac{\partial A}{\partial u} \wedge (W_0)_\gamma$ is also bounded away from zero. That is, $A$ is transverse to $(W_0)_\gamma$, thus gives rise to a Birkhoff section under the collapsing map $M_\gamma \to M$.

\end{constr}

\begin{constr}[Interpolating $F$] \label{constr:interpolateF}
Consider $\nu = [-\epsilon,\epsilon]_x \times [-\epsilon,\epsilon]_y  \times S^1_z$. Suppose $(\alpha_-)_1$ is a negative contact form defined by \Cref{eq:negcontactformcorrfunction} for some $(L,F_1)$ satisfying the items of \Cref{lemma:negcontactformcorrfunctiontaylorexpansion}.
Let $F_0$ be an arbitrary function on $D^2 \times S^1$ satisfying item (\ref{item:negcontactformcorrfunctiontaylorexpansionF=0}) of \Cref{lemma:negcontactformcorrfunctiontaylorexpansion}, i.e., that $F_0 = 0$ along $\{0\} \times S^1$.

We once again pick $\rho\colon \mathbb{R}_{\geq 0} \to [0,1]$ a function such that $\rho(t) = 0$ for $t$ close to $0$ and $\rho(t) = 1$ for $t \geq 1$.
For $s \in [0,1]$ and $\delta > 0$, we define
$$\widehat{F}_s(x,y,z) = F_1(x,y,z) + (1-s) \left( 1- \rho \left( \frac{r}{\delta} \right) \right) (F_0(x,y,z) - F_1(x,y,z))$$
where $r = \sqrt{x^2+y^2}$.
Then $\widehat{F}_1 = F_1$ and $\widehat{F}_0 = F_0$ in a smaller tubular neighborhood of $\{0\} \times S^1$.

We claim that for sufficiently small $\delta$, $(L,\widehat{F}_s)$ satisfies \Cref{eq:negcontactformcorrfunctioncondition1} and \Cref{eq:negcontactformcorrfunctioncondition2} for all $s \in [0,1]$.
Indeed, \Cref{eq:negcontactformcorrfunctioncondition1} is clear since $L$ is unchanged. 

For \Cref{eq:negcontactformcorrfunctioncondition2}, we have
\begin{align*}
& \left( -L\frac{\partial^2 L}{\partial y^2}+W\widehat{F}_s \right) - \left( -L\frac{\partial^2 L}{\partial y^2}+WF_1 \right) \\
&= W \left( \widehat{F}_s - F_1 \right) \\
&= y \frac{\partial^2 L}{\partial y^2} \frac{\partial}{\partial x} \left( \widehat{F}_s - F_1 \right) + \left( \frac{\partial L}{\partial x} + \frac{\partial^2 L}{\partial y \partial z} - y \frac{\partial^2 L}{\partial x \partial y} \right) \frac{\partial}{\partial y} \left( \widehat{F}_s - F_1 \right) - \frac{\partial^2 L}{\partial y^2} \frac{\partial}{\partial z} \left( \widehat{F}_s - F_1 \right)
\end{align*} %
From item (\ref{item:negcontactformcorrfunctiontaylorexpansionF=0}) of \Cref{lemma:negcontactformcorrfunctiontaylorexpansion}, we have $F_0 = F_1 = 0$ along $\{0\} \times S^1$, thus $F_0$ and $F_1$ are $O(r)$.
Thus 
\begin{align*}
& \left| y \frac{\partial^2 L}{\partial y^2} \frac{\partial}{\partial x} \left( \widehat{F}_s - F_1 \right) \right| \\
&\leq \left| - y \frac{\partial^2 L}{\partial y^2} (1-s) \rho' \left( \frac{r}{\delta} \right) \frac{1}{\delta} \frac{\partial r}{\partial x} (F_0 - F_1) + y \frac{\partial^2 L}{\partial y^2} (1-s) \left( 1- \rho \left( \frac{r}{\delta} \right) \right) \left( \frac{\partial F_0}{\partial x} - \frac{\partial F_1}{\partial x} \right) \right| \\
&\leq |y| \left| \frac{\partial^2 L}{\partial y^2} \right| \left| \rho' \left( \frac{r}{\delta} \right) \right| \frac{1}{\delta} (|F_0| + |F_1|) + |y| \left| \frac{\partial^2 L}{\partial y^2} \right| \left| 1- \rho \left( \frac{r}{\delta} \right) \right| \left( \left| \frac{\partial F_0}{\partial x} \right| + \left| \frac{\partial F_1}{\partial x} \right| \right) \\
&= O(\delta) \cdot O(1) \cdot O(1) \cdot O(\delta^{-1}) \cdot O(\delta) + O(\delta) \cdot O(1) \cdot O(1) \cdot O(1) = O(\delta)
\end{align*}

Meanwhile, from item (\ref{item:negcontactformcorrfunctiontaylorexpansionLfirstorder}) of \Cref{lemma:negcontactformcorrfunctiontaylorexpansion}, we have $\frac{\partial L}{\partial x} = \frac{\partial^2 L}{\partial y \partial z} = 0$ along $\{0\} \times S^1$, thus $\frac{\partial L}{\partial x} + \frac{\partial^2 L}{\partial y \partial z} - y \frac{\partial^2 L}{\partial x \partial y} = O(r)$.
Computing as above, we then have 
$$\left( \frac{\partial L}{\partial x} + \frac{\partial^2 L}{\partial y \partial z} - y \frac{\partial^2 L}{\partial x \partial y} \right) \frac{\partial}{\partial y} \left( \widehat{F}_s - F_1 \right) = O(\delta).$$

Finally, by item (\ref{item:negcontactformcorrfunctiontaylorexpansionF=0}) of \Cref{lemma:negcontactformcorrfunctiontaylorexpansion}, we also have $\frac{\partial F_0}{\partial z} = \frac{\partial F_1}{\partial z} = 0$ along $\{0\} \times S^1$, thus $\frac{\partial F_0}{\partial z}$ and $\frac{\partial F_1}{\partial z}$ are $O(r)$.
\begin{align*}
\left| \frac{\partial^2 L}{\partial y^2} \frac{\partial}{\partial z} \left( \widehat{F}_s - F_1 \right) \right| &= \left| \frac{\partial^2 L}{\partial y^2} (1-s) \left( 1 - \rho \left( \frac{r}{\delta} \right) \right) \left( \frac{\partial F_0}{\partial z} - \frac{\partial F_1}{\partial z} \right) \right| \\
&\leq \left| \frac{\partial^2 L}{\partial y^2} \right| \left| 1 - \rho \left( \frac{r}{\delta} \right) \right| \left( \left| \frac{\partial F_0}{\partial z} \right| + \left| \frac{\partial F_1}{\partial z} \right| \right) \\
&= O(1) \cdot O(1) \cdot O(\delta) = O(\delta).
\end{align*}

This verifies that $(-L\frac{\partial^2 L}{\partial y^2}+W\widehat{F}_s) - (-L\frac{\partial^2 L}{\partial y^2}+WF_1) = O(\delta)$, thus for small enough $\delta$, \Cref{eq:negcontactformcorrfunctioncondition2} for $\widehat{F}_s$ follows from that for $F_1$.

We can thus define a family of negative contact forms $(\alpha_-)_s$ by applying \Cref{lemma:negcontactformcorrfunction} to $(L,\widehat{F}_s)$. By definition, $((\alpha_-)_s,\xi_+)$ is strongly adapted for every $s$.
Since $W$ is unchanged, $\{0\} \times S^1$ is a closed orbit of the Reeb flow $(R_-)_s$ for every $s$.
Also from the fact that $W$ is unchanged, we see that any $\partial$-strong Birkhoff section $S$ for $R_-$ is automatically a Birkhoff section for $(R_-)_s$, for every $s$.

\end{constr}

\subsection{Proof of \Cref{prop:standardizeorbit}} \label{subsec:standardizeorbitproof}

We now collect the ingredients to present a complete proof of \Cref{prop:standardizeorbit}.

\begin{proof}[Proof of \Cref{prop:standardizeorbit}]
We first perturb $\alpha_-$ so that $\frac{1}{2\pi} \int_\gamma \alpha_-$ is rational.
Then we apply \Cref{lemma:negcontactformcorrfunctiontaylorexpansion} to obtain a tubular neighborhood $\nu$ of $\gamma$ where we can express $\alpha_-$ as \Cref{eq:negcontactformcorrfunction} for some functions $(L,F)$ satisfying the items in \Cref{lemma:negcontactformcorrfunctiontaylorexpansion}.
Moreover, note that $\Lambda = L(0,0,z) = \frac{1}{2\pi} \int_\gamma \alpha_-$ is a positive rational number.

Next, we apply \Cref{constr:removecubicterm} to homotope the bicontact package so that $L$ has no cubic term, i.e., so that $L = \Lambda + \frac{1}{2} Q(z) x^2 - \frac{1}{2} y^2$ in a smaller tubular neighborhood of $\gamma$.
After this construction, we can assume that the given $\partial$-strong Birkhoff section $S$ to the initial Reeb flow remains a $\partial$-strong Birkhoff section to the modified Reeb flow.

Then, we apply \Cref{constr:interpolateQ} to homotope the bicontact package so that $L$ is standard, i.e., so that it is given by \Cref{eq:stdL}, in a smaller tubular neighborhood of $\gamma$.
If this is all we wish to accomplish, then we can take $Q(z,t)$ in the construction to be any smooth function on $S^1 \times [0,1]$ where $Q(z,1) = Q(z)$ and $Q(z,0) = \omega_0^2$ for any positive rational number $\omega_0$.

If we want to maintain a Birkhoff section, we have to work harder.
We identify the slopes on $\partial \nu$ so that $\frac{\partial}{\partial z}$ is $0$ and $\frac{\partial}{\partial \theta}$ is $\infty$. 
Let $\sigma$ be the slope of intersection between the Birkhoff section $S$ and $\partial \nu$.
From \Cref{eq:normalblowup}, the blown-up Reeb flow is currently a positive multiple of
$$r ( -1 + Q(z) ) \sin \theta \cos \theta \frac{\partial}{\partial r} +(\sin^2 \theta + Q(z) \cos^2 \theta ) \frac{\partial}{\partial \theta} + \frac{\partial}{\partial z}.$$
In particular, the foliation by flow lines on the boundary torus has slope $f(Q)$.
Since $\gamma$ is a negative boundary orbit, and since $S$ is $\partial$-strong, we have $0 \leq f(Q) < \sigma$.

We take a constant $K$ such that $K < Q(z)$ for all $z$.
By \Cref{lemma:flowlineslope}(1,3,4), we can pick a positive rational number $\omega_0$ such that $f(\omega_0^2) \in (0,\sigma)$.
We then pick $Q(z,t)$ so that
\begin{enumerate}[label=(\roman*)]
    \item $Q(z,1) = Q(z)$,
    \item $Q(z,t)$ is increasing in $t$ for $t \in [\frac{1}{2},1]$,
    \item $Q(z,\frac{1}{2}) = K$ for all $z$,
    \item $Q(z,t)$ is independent of $z$ for $t \in [0,\frac{1}{2}]$, 
    \item $Q(z,t)$ is decreasing in $t$ for $t \in [0,\frac{1}{2}]$, and %
    \item $Q(z,0) = \omega_0^2$.
\end{enumerate}
Then by (ii) and \Cref{lemma:flowlineslope}(2), we have $f(Q(\cdot,t)) < \sigma$ for $t \in [\frac{1}{2},1]$, and by (iv), (v), \Cref{lemma:flowlineslope}(2), and the definition of $\omega_0$, we have $f(Q(\cdot,t)) < \sigma$ for $t \in [0,\frac{1}{2}]$.
Thus as argued for in \Cref{constr:interpolateQ}, we can redefine $S$ to be a $\partial$-strong Birkhoff section for the modified Reeb flow.

Finally, we apply \Cref{constr:interpolateF} to homotope the bicontact package so that $F$ is standard, i.e., so that it is given by \Cref{eq:stdF}, in a smaller tubular neighborhood of $\gamma$.
After this construction, we can assume that the given $\partial$-strong Birkhoff section $S$ to the initial Reeb flow remains a $\partial$-strong Birkhoff section to the modified Reeb flow.

Once \Cref{eq:stdL} and \Cref{eq:stdF} are satisfied, $\gamma$ is a standard Reeb orbit of the Reeb flow, as explained in \Cref{rem_standardReeb_L_F}.
\end{proof}

We take the occasion to supply a proof of \Cref{thm:paddingaddwinding0twists}.

\begin{proof}[Proof of \Cref{thm:paddingaddwinding0twists}]
We first apply \Cref{prop:standardizeorbit} to replace the given bicontact package by one that has a Birkhoff section with the same monodromy and induced line field, up to homotopy, but with boundary orbits being standard Reeb orbits.
Then we apply \Cref{thm:paddingaddwinding0twistsstdboundary}. Here we use the fact that the homotopy class of the induced line field is preserved, so that the notion of `winding number zero curves' stays unchanged.
\end{proof}

\section{Application to the Christy--Fried--Ghys conjecture} \label{sec:mainthmproof}

In this section, we take the final steps towards proving our main theorems.
In \Cref{subsec:birkhoffsectionanosovflow}, we will construct an appropriate Birkhoff section to the Anosov flow.
In \Cref{subsec:transferbirkhoffsection}, we will show that such a Birkhoff section can be transferred to the Reeb flow $R_-$.
We can then apply \Cref{thm:paddingaddwinding0twists} and \Cref{cor:factorizationthmoneboundary} to such a Birkhoff section.

\subsection{Birkhoff sections of Anosov flows} \label{subsec:birkhoffsectionanosovflow}

A theorem of Marty \cite[Theorem F]{Mar25} states that every skew Anosov flow admits an embedded Birkhoff section with only one boundary component. 
The main result of this subsection is a slight refinement that says that the genus of such a Birkhoff section can be made arbitrarily large.

\begin{thm} \label{thm:constructbirkhoffoneboundary}
Let $\phi$ be a negative skew Anosov flow.
For every $N > 0$, there exists an embedded Birkhoff section $S$ to $\phi$ with genus $\geq N$ and only one boundary component, which is negative.
\end{thm}

Our method of proof is essentially just to retrace the argument of Marty while keeping track of boundary slopes in order to compute the genus of the final Birkhoff section.

To this end, we set up some notation.
Recall from \Cref{rmk:partialBirkhoffsection} that a partial Birkhoff section is a Birkhoff section that does not necessarily meet every forward and backward trajectory.
Let $S$ be a partial Birkhoff section to an Anosov flow $\phi$ with orientable foliations.
Then the weak stable plane field $E^s \oplus \langle X \rangle$ induces a nonsingular orientable line field $(E^s \oplus \langle X \rangle) \cap TS$ on $\intr(S)$.
In particular, we can talk about the winding number of a boundary component $b$ of $S$ with respect to this line field. We denote this number by $\wind_\phi(S,b)$. 

Suppose $\gamma$ is a closed orbit of $\phi$. 
Let $\nu$ be a small tubular neighborhood of $\gamma$.
We define the \textbf{boundary class} of $S$ at $\gamma$ to be the homology class of the intersection $S \cap \partial \nu$, and we denote it by $\partial_\gamma S$.
Note that we do not require that $\gamma$ be a boundary orbit of $S$ in this definition.

More quantitatively:
Let $\mu$ be the oriented meridian on $\partial \nu$, i.e., the meridian that bounds a disc that is positively transverse to $\phi$. Let $\lambda$ be the degeneracy slope on $\partial \nu$, oriented in the direction of the flow.
Note that $(\mu, \lambda)$ is an oriented basis of $H_1(\partial \nu)$.
Writing $\partial_\gamma S = p \mu + q \lambda$, we refer to $q$ as the \textbf{multiplicity} of $S$ at $\gamma$.
Note that $p \leq 0$. If $S$ is a Birkhoff section, then $p < 0$. 
In such a case, we also have that $\gamma$ is a boundary orbit of $S$ if and only if $q \neq 0$, and $q = \pm 1$ if and only if $S$ is embedded at $\gamma$.

Refining these notions, let $b$ be a boundary component of $S$ lying along $\gamma$.
We define the \textbf{boundary class} of $S$ at $b$ to be the homology class of the component of $S \cap \partial \nu$ lying near $b$, and we denote it by $\partial_b S$.
Writing $\partial_b S = p \mu + q \lambda$, we refer to $q$ as the \textbf{multiplicity} of $S$ at $b$.
The winding number $\wind_\phi(S,b)$ can be computed as $p$. In particular, we always have $\wind_\phi(S,b) \leq 0$.
We say that $b$ is a \textbf{positive} boundary component if $q > 0$, and a \textbf{negative} boundary component if $q < 0$.
If $S$ is a Birkhoff section, then we have $p < 0$ and $q \neq 0$, and the above definition of sign agrees with that given in \Cref{subsec:birkhoffsection}.

\begin{figure}
    \centering
    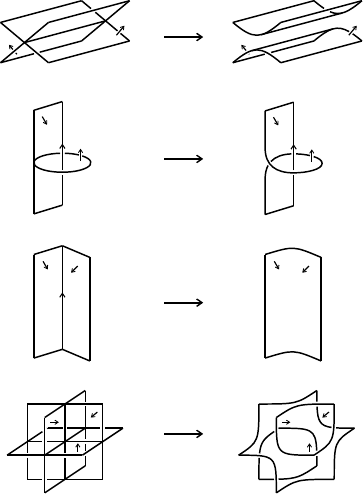
    \caption{Two Birkhoff sections are Fried summed by taking their union and resolving the curves and arcs of intersection.}
    \label{fig:friedsum}
\end{figure}

The union of two (partial) Birkhoff sections can be turned into a (partial) Birkhoff section by resolving the curves and arcs of self-intersection as in \Cref{fig:friedsum}. 
We refer to this as \textbf{Fried summing} the (partial) Birkhoff sections, as it was first introduced in this context in \cite{Fri83}.
The boundary class at a fixed orbit is additive under Fried sum. In particular, the multiplicity at a fixed orbit is also additive under Fried sum.

Generalizing, we can Fried sum any finite collection of (partial) Birkhoff sections.

For the rest of this subsection, we fix a negative skew Anosov flow $\phi$ with orientable foliations. 
We also fix an orientation on the line fields $E^s$ and $E^u$ so that if $e^{s/u}$ is a positive vector in $E^{s/u}$, then $(e^s, X, e^u)$ is a positive basis of $TM$.

Let $D$ be a transverse disc to $\phi$. 
The weak stable and unstable foliations of $\phi$ induce transverse 1-dimensional foliations $D^s$ and $D^u$ on $D$.
A \textbf{rectangle} on $D$ is a subsurface $R \subset D$ bounded by two segments lying on leaves of $D^s$ and two segments lying on leaves of $D^u$. 
Each rectangle can be written as $I^s \times I^u$, where $I^{s/u}$ is an interval in the space of leaves of $D^{u/s}$ respectively.
Up to shrinking $D$, we can assume that $D$ is itself a rectangle.

Meanwhile, the orientations on $E^{s/u}$ induce orientations on $D^{s/u}$. Given two points $x,y \in D$, we write $x >_s y$ if the $D^s$-leaf containing $x$ is on the positive side of the $D^s$-leaf containing $y$. We define $x >_u y$ similarly.

\begin{proof}[Proof of \Cref{thm:constructbirkhoffoneboundary}]
As mentioned at the start of the subsection, we simply retrace the argument in \cite{Mar25}, so we will be brief on the points that are already explained there.

From the proof of \cite[Theorem F]{Mar25}, there exists a Birkhoff section $A$ (obtained by desingularizing an immersed Birkhoff annulus) with two boundary components lying along $\gamma_1$ and $\gamma_2$, with multiplicity $-n_1$ and $-n_2$ respectively, where $n_1, n_2 \in \mathbb{Z}_+$, so that $\gamma_1$ and $\gamma_2$ intersect some transverse rectangle $R$ at some points $x_1$ and $x_2$ with $x_1 <_s x_2$ and $x_1 >_u x_2$.

By \cite[Lemma 32]{Mar25}, for every $k \in \mathbb{Z}_+$, there exists a partial Birkhoff section $P$ (obtained by desingularizing an immersed Fried pants) with one boundary component lying along $\gamma_1$ with multiplicity $2kn_1$, one boundary component lying along $\gamma_2$ with multiplicity $2kn_2$, and one boundary component lying along some $\gamma_3$ with multiplicity $-1$.

We take the Fried sum $S = 2kA + P$.
Then for $i=1,2$, the multiplicity of $S$ at $\gamma_i$ is $2k(-n_i) + 2kn_i = 0$, so $\gamma_i$ is not a boundary orbit of $S$.
Meanwhile, the boundary class of $S$ at $\gamma_3$ is 
\begin{align*}
\partial_{\gamma_3} S &= 2k \partial_{\gamma_3} A + \partial_{\gamma_3} P \\
&= 2k (- \langle A, \gamma_3 \rangle \mu) + (\wind_\phi(P,\gamma_3) \mu - \lambda) \\
&= (- 2k \langle A, \gamma_3 \rangle + \wind_\phi(P,\gamma_3)) \mu - \lambda.
\end{align*}
In particular, the multiplicity of $S$ at $\gamma_3$ is $-1$, so $S$ is embedded and has only one negative boundary component along $\gamma_3$.
Moreover, we have $\wind_\phi(S,\gamma_3) = - 2k \langle A, \gamma_3 \rangle + \wind_\phi(P,\gamma_3)$, so by \Cref{prop:windpoincarehopf}, we have $\genus(S) = k \langle A, \gamma_3 \rangle - \frac{\wind_\phi(P,\gamma_3)+1}{2} \geq k - \frac{1}{2} \to \infty$ as $k \to \infty$.
\end{proof}

\subsection{Transferring the Birkhoff section to the Reeb flow} \label{subsec:transferbirkhoffsection}

In this subsection, we explain how, for a negative skew Anosov flow $\phi$, one can transfer a Birkhoff section for $\phi$ to a Birkhoff section for the Reeb flow $R_-$.

\begin{prop} \label{prop:swapanosovR_-}
Suppose $\phi$ is a negative skew Anosov flow.
Then up to isotopic equivalence, $\phi$ admits a strongly adapted bicontact package $(\alpha_-,\xi_+)$ where the Reeb flow $R_-$ is isotopically equivalent to $\phi$.
\end{prop}
\begin{proof}
By \cite[Theorem A]{Mar25}, up to isotopic equivalence we can assume that $\phi$ is the Reeb flow of a negative contact form $\alpha_-$.
Let $X$ be the generating vector field of $\phi$.
Pick a strongly adapted bicontact package $(\widehat{\alpha}_-,\xi_+)$ supporting $X$. 
Let $Y$ be a nonvanishing vector field lying in the intersection of planes $\ker \alpha_- \cap \xi_+$.
It is argued in \cite[Proposition 4.15]{Bar25} that $Y$ must be isotopically equivalent to $\phi$.

For the reader's convenience we briefly sketch the argument: %
Since the Reeb vector field to $\alpha_-$ is $X$, which lies in $\xi_+$, $(\alpha_-,\xi_+)$ is a strongly adapted bicontact package supporting $Y$. 
By \Cref{thm:hozoori}, $Y$ is an Anosov flow. 
Meanwhile, since $Y$ lies in $\xi_+$, it is transverse to the stable foliation of $\phi$. By results of Fenley \cite[Main theorem]{Fen05}, $Y$ must be isotopically equivalent to $\phi$.

Thus to arrive at the statement of the proposition, we consider the flow generated by $Y$. This flow is isotopically equivalent to $\phi$ and supported by $(\alpha_-,\xi_+)$, and the Reeb flow of $\alpha_-$ is exactly $\phi$.
\end{proof}

\begin{cor} \label{cor:transferbirkhoffsection}
Suppose $\phi$ is a negative skew Anosov flow.
Then up to isotopic equivalence, $\phi$ admits a strongly adapted bicontact package $(\alpha_-,\xi_+)$ where the Reeb flow $R_-$ admits an embedded Birkhoff section $S$ with genus $\geq 5$ and only one boundary component $b$, which is negative, and so that $\wind_{\xi_+}(S,b) < 0$.
\end{cor}
\begin{proof}
We first pick a strongly adapted bicontact package as in \Cref{prop:swapanosovR_-}, and apply \Cref{thm:constructbirkhoffoneboundary} to construct a Birkhoff section for $\phi$ with genus $\geq 5$ and only one boundary component $b$, which is negative.
Finally, since by \Cref{prop:contactstructurerudder}, $\xi_+$ has to stay within the negative quadrants of $R_-$, we have $\wind_{\xi_+}(S,b) = \wind_{R_-}(S,b) < 0$.
\end{proof}

\subsection{Proof of main results} \label{subsec:mainthmproof}

\begin{proof}[Proof of \Cref{thm:flowreducetogeneralizedht}]
Let $\phi$ be a skew Anosov flow. Up to reversing the orientation of the manifold, we can assume that it is a negative skew Anosov flow.
By \Cref{cor:transferbirkhoffsection}, up to isotopic equivalence, $\phi$ admits a strongly adapted bicontact package $(\alpha_-,\xi_+)$ where the Reeb flow $R_-$ admits an embedded Birkhoff section $S$ with genus $\geq 5$ and one negative boundary component $b$ lying along some closed orbit $\gamma$, and so that $\wind_{\xi_+}(S,b) < 0$.

Let $f$ be the monodromy of $S$ and let $\ell$ be the line field induced by $\xi_+$ on $\intr(S)$.
As explained in \Cref{subsec:factoroutwindzerocurve}, we have $f \in \PMod(\intr(S))[\ell]$.
Thus by \Cref{cor:factorizationthmoneboundary}, there exist curves $c_i$ that have zero winding number with respect to $\ell$ such that $f = \prod \tau_{c_i}^{n_i}$.

We apply \Cref{prop:standardizeorbit} to make $\gamma$ a standard Reeb orbit.
Here we use the fact that $\gamma$ is a negative boundary orbit.

Then by \Cref{thm:paddingaddwinding0twistsstdboundary}, we can transform $\phi$ to an Anosov flow $\widehat{\phi}$ on $\widehat{M}$ by bicontact padding surgery then Goodman surgeries, so that $\widehat{\phi}$ is supported by a strongly adapted bicontact package $(\widehat{\alpha}_-,\widehat{\xi}_+)$ where $\widehat{R}_-$ admits a Birkhoff section $\widehat{S} \cong S$ with trivial monodromy.

Then $\widehat{M}$ must be a Seifert fibered space.
By \cite[Theorem 4.9]{Bar96}, $\widehat{\phi}$ is a fiberwise finite cover of a geodesic flow. 
Moreover, $\widehat{M}$ has at most one singular fiber, so the base orbifold of the geodesic flow has at most one cone point.
\end{proof}

\begin{proof}[Proof of \Cref{cor:dehnsurgerytosuspension}]
By \cite[Corollary B]{ABM24}, every $(M,\phi) = (M_0,\phi_0)$ as in the corollary can be transformed into a skew Anosov flow $(M_1,\phi_1)$ via Goodman surgery, which is topologically a Dehn surgery along closed orbits.

We then enter the jurisdiction of \Cref{thm:flowreducetogeneralizedht}; we use the notation in its proof. 
We standardize boundary orbits (\Cref{prop:standardizeorbit}) and perform bicontact padding surgery to transform $(M_1,\phi_1)$ to $(M_2,\phi_2)$.
The surgery is topologically Dehn surgery along closed orbits of the Reeb flow $R_-$ (both before and after standardization). But since $R_-$ (before standardization) is isotopically equivalent to $\phi_1$, these orbits are isotopic to a collection of orbits of $\phi_1$. %
We then perform Goodman surgery to transform $(M_2,\phi_2)$ into a finite cover of a geodesic flow $(M_3,\phi_3)$.

In particular, $\phi_3$ is a transitive Anosov flow defined on a graph manifold, so we can apply \cite[Theorem 1.2]{Tsa26} and \cite[Theorem 1.4]{Tsa24b} to transform $(M_3,\phi_3)$ into an Anosov suspension flow $(M_4,\phi_4)$ via Goodman surgery.

\end{proof}

\section{Discussion and further questions} \label{sec:furtherquestions}

\subsection{Variations of \Cref{thm:flowreducetogeneralizedht}} \label{subsec:variation}

We expect the techniques of this paper to apply to more general Anosov flows. As an illustration, we show the following theorem.

\begin{thm} \label{thm:reducetototallyfree}
Let $\phi$ be a, possibly non-transitive, Anosov flow with orientable foliations. If $\phi$ is supported by a strongly adapted bicontact package $(\alpha_-,\xi_+)$ where the Reeb flow $R_-$ admits an embedded Birkhoff section $S$ with genus $\geq 5$, all boundary components being negative, and at most two boundary components having even winding numbers, then $\phi$ can be transformed into a totally free Anosov flow by bicontact padding and Goodman surgery. In particular, $\phi$ can be transformed into a suspension.
\end{thm}

Here, a \textbf{totally free Anosov flow} is an Anosov flow on a graph manifold for which the restriction to each Seifert fibered piece is a finite fiberwise cover of a geodesic flow. (These flows were originally called \emph{generalized Handel--Thurston flows} in \cite{Bar96}. We follow instead the terminology used in \cite{BarMan25}.)  Note that we do not know whether the assumptions of \Cref{thm:reducetototallyfree} imply that $\phi$ must be transitive, but the conclusion is that after bicontact padding surgery it will become transitive.

\begin{proof}[Proof of \Cref{thm:reducetototallyfree}]
Let $f$ be the monodromy of $S$ and let $\ell$ be the line field induced by $\xi_+$ on $\intr(S)$. 
Since $S$ is embedded, $f$ is an element of $\PMod(\intr(S))[\ell]$.
By \Cref{cor:factorizationthmtwoevenboundary}, there exist winding number zero curves $c_i$ and a disjoint collection of curves $\{d_1,...,d_k\}$ such that $f = \prod \tau_{c_i}^{n_i} \prod \tau_{d_j}^{m_j}$ for some integers $n_i, m_j$, and such that each complementary region of $\bigcup_{j=1}^k d_j$ has genus at least one.

We apply \Cref{prop:standardizeorbit} to make all boundary orbits of $S$ standard.
Here we use the fact that the boundary components of $S$ are all negative.

Then by \Cref{thm:paddingaddwinding0twists}, we can transform $\phi$ to an Anosov flow $\widehat{\phi}$ on $\widehat{M}$ by bicontact padding surgery then Goodman surgeries, so that $\widehat{\phi}$ is supported by a strongly adapted bicontact package $(\widehat{\alpha}_-,\widehat{\xi}_+)$ where $\widehat{R}_-$ admits a Birkhoff section $\widehat{S} \cong S$ with monodromy $\prod \tau_{d_j}^{m_j}$.
Right away, this implies that $\widehat{M}$ is a graph manifold.

To argue that $\widehat{\phi}$ can be assumed to be totally free, we have to work slightly harder:
The base orbifolds of the JSJ pieces of $\widehat{M}$ are the complementary regions of $\bigcup_{j=1}^k d_j$, with boundary components of $\widehat{S}$ capped off, possibly with a cone point.%
It can be shown that for every nonsingular line field on a surface of genus $\geq 1$, there is always a winding number zero curve. See for example \cite[Proposition 2.15]{CS23}.
By including such a curve when applying \Cref{prop:paddingrealizeanosovorbit}, we can assume that each JSJ piece of $\widehat{M}$ contains a closed orbit of $\widehat{\phi}$ that is not homotopic to the Seifert fiber.

Now, we have a very good understanding of Anosov flows on graph manifolds due to work of Barbot and Fenley \cite{BF13, BF15, BF21}:
Say that a Seifert fibered piece is \textbf{periodic} if its Seifert fibers are freely homotopic to a closed orbit, otherwise it is \textbf{free}.
Among the homotopy classes contained in the fundamental group of a periodic piece, the only ones that represent a closed orbit contained in that piece are the Seifert fiber and its multiples.
Thus in our manifold $\widehat{M}$, all JSJ pieces are free.
By \cite[Théorème B]{Bar96}, this implies that $\widehat{\phi}$ is totally free, in the sense defined above.

\cite[Théorème B]{Bar96} also states that totally free Anosov flows are skew, thus transitive. So one can apply \cite[Theorem 1.2]{Tsa26} and \cite[Theorem 1.4]{Tsa24b} to transform an Anosov flow output by \Cref{thm:reducetototallyfree} into an Anosov suspension flow via Goodman surgery, as in \Cref{cor:flowreducetosusp}.
\end{proof}

To ascertain the applicability of \Cref{thm:reducetototallyfree}, we ask the following question.

\begin{quest} \label{quest:twoevenboundary}
Is every Anosov flow supported by a strongly adapted bicontact package $(\alpha_-,\xi_+)$ where the Reeb flow $R_-$ admits an embedded Birkhoff section $S$ with genus $\geq 5$, all boundary components being negative, and at most two boundary components having even winding numbers?
\end{quest}

We remark that Kim--Kim--van Koert \cite[Theorem 1.1]{KKvK22} showed that given any $N \in \mathbb{Z}_+$, there exist Reeb flows that do not admit embedded Birkhoff sections with $\leq N$ boundary components. As a point of comparison, Tsang \cite{Tsa24} showed that every transitive pseudo-Anosov flow admits a Birkhoff section with two boundary components.

Related to the thread of non-transitive Anosov flows, it might also be interesting to generalize \Cref{thm:flowreducetogeneralizedht} to \textit{hyperbolic plugs}, which are essentially Anosov flows on manifolds with boundary, and which can be glued to form Anosov flows on closed manifolds. See \cite{BBY17, Pau25}.
This should require a suitable generalization of strongly adapted bicontact packages to manifolds with boundary, which has been recently developed in \cite{Hoz26}.
On the other hand, it is less clear to us what a suitable generalization of Birkhoff sections should be in this context.

Finally, it might be interesting to generalize the ideas of this paper to Anosov flows with non-orientable foliations.
For simplicity, we keep the manifold oriented. Then the bicontact machinery should still work in this case: One has to consider contact forms that are only locally well-defined up to a sign. The new \textbf{non-orientable} contact forms will then have Reeb vector fields that are also only locally well-defined up to a sign, but the strongly adapted condition still makes sense.
This should culminate in the following analogue of \Cref{thm:hozoori}, although we haven't checked all the details.

\begin{conj} \label{conj:hozoorinonorientable}
A flow on a closed oriented 3-manifold is Anosov with non-orientable foliations if and only if it is supported by a strongly adapted non-orientable bicontact package.
\end{conj}

For the purposes of this discussion, let us take \Cref{conj:hozoorinonorientable} for granted. 

Since the Reeb vector field is only locally well-defined up to a sign, we no longer have a Reeb flow. Instead, we only have a 1-dimensional foliation. Nevertheless, the notion of a Birkhoff section carries over, by replacing \Cref{defn:birkhoffsection}(3) by `For every $p \in M$, choose a branch of the Reeb vector field $R_-$ along the leaf containing $p$, then there exists $T_- < 0 < T_+$ such that $R_-^{T_-}(p), R_-^{T_+}(p) \in S$.' %
In the non-orientable case, a Birkhoff section would now determine a semi-fibration, i.e., a fibration over the orbifold $[0,1]$ with reflector boundary, away from the boundary orbits. %
Nevertheless, it should still be possible to define bicontact padding surgery along such Birkhoff sections.

Instead, the main difficulty lies in supplying the appropriate mapping class group factorization result.
We would like to perform Dehn twists along curves on $S$ so that the resulting manifold is, say, a graph manifold. In particular, the double cover will also be a graph manifold, hence we would like for its monodromy to be a product of Dehn twists along disjoint curves. %
But if we do a single Dehn twist along $c \subset S$, then the monodromy of the double cover is actually changed by two Dehn twists, along the two lifts of $c$ on the two lifts of $S$.
This seems to suggest that the factorization result is more delicate in this case, since one is forced to add Dehn twists two-by-two, instead of one-by-one.

\subsection{Some ideas towards proving \Cref{conj:friedchristyghys}} \label{subsec:friedchristyghysideas}

\begin{figure}
    \centering
    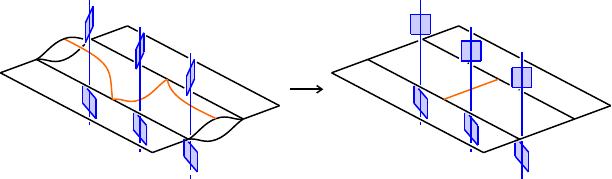
    \caption{Performing a positive Dehn twist adds left-handed rotation.}
    \label{fig:dehntwistsaddrotation}
\end{figure}

Given \Cref{cor:dehnsurgerytosuspension}, an approach to proving \Cref{conj:friedchristyghys} is to show that bicontact padding surgery can be realized by a sequence of Goodman surgeries.
Here is a preliminary idea for this: 
Let $(\alpha_-,\xi_+)$ be a strongly adapted bicontact package. Suppose $S$ is a Birkhoff section to the Reeb flow $R_-$.
Morally, performing a positive Dehn twist along a curve $c$ on $S$ adds (left-handed) rotation of $\xi_+$. See \Cref{fig:dehntwistsaddrotation}.

Now, if $S$ has, say, one boundary component, then we have a chain relation in the mapping class group of $S$, which expresses a positive Dehn twist along $\partial S$ as a product of positive Dehn twists along non-isolating curves in $\intr(S)$. See \cite[Section 9.4.2]{FM12}.
This should allow one to `insert rotation' while only changing the monodromy by boundary twists.
Thus at least on $\xi_-$, this has the same effect as bicontact padding surgery.
Together with some manner of control on $\xi_+$, such as the ideas explained in \Cref{subsec:poscontactstructureopenbook}, it is conceivable that such an operation can be shown to be equivalent to bicontact padding surgery.

However, there are still a number of obstacles to overcome in pursuing this idea:
For example, we need to make sure that we can maintain a strongly adapted bicontact package when performing these twists.
As in this paper, this is possible if the chain of curves are isotopic to orbits of the Anosov flow, since then we can just do Goodman surgery along the orbits.
However, as explained in \Cref{subsec:intropaddingmotivation}, we do not know how to arrange for this without doing padding in the first place.

A weaker hypothesis that could work here is if the chain of curves are isotopic to $\xi_-$-Legendrian-$\xi_+$-transverse curves, since then we can do Legendrian-transverse surgery.
Regardless, we still face the issue of making the curves $\xi_-$-Legendrian. The only method we know to do this is to apply the Legendrian realization principle, but in general, we need to do padding in the first place to ensure there is enough room for this.

Also, in both of these approaches, the chain of curves will necessarily have zero winding number. 
In general, there is an obstruction to this being possible, in terms of an Arf invariant associated to the nonsingular line field on $S$ induced by $\xi_+$. See \cite[Section 2.2 and Proposition 2.15]{CS23}.
It is unclear to us how to control this Arf invariant when applying the usual methods of constructing Birkhoff sections, for example summing with Fried pants.

\subsection{Open books supporting $\ker \alpha_-$ and the isotopy class of $\xi_+$}
\label{subsec:poscontactstructureopenbook}

Let $S$ be a negative embedded Birkhoff section for $R_-$ in a
strongly adapted bicontact package $(\alpha_-,\xi_+)$. The section $S$,
together with its monodromy $h$, defines an open book decomposition
supporting $\ker\alpha_-$; see \cite[Section 15]{Hon04}. By Giroux's
correspondence, this open book determines the isotopy class of
$\ker\alpha_-$. 
 A natural question from the point of view of bicontact structures, is whether one can incorporate additional information to the pair $(S,h)$ to not only recover $\ker \alpha_-$, as given by Giroux's correspondence, but also recover $\xi_+$ up to isotopy. This question is particularly interesting in light of the recent work of Bowden--Colin \cite{BC26} showing that, for bicontact structures that support a transitive Anosov flow, knowing the two contact structures up to isotopy independently, is (almost always) enough to recover the Anosov flow up to isotopy equivalence.

One potential candidate for this additional data that carries some information about $\xi_+$ is the line field $\ell_+$ induced by $\xi_+$ on $S$, or even just its homotopy class $[\ell_+]$. Thus we ask:
\begin{quest}
Does the open book $(S,h)$, together with the (homotopy class of the) line field $\ell_+$ induced by $\xi_+$ on the chosen fixed section $S$, determine
the isotopy class of $\xi_+$ among positive contact structures
arising from strongly adapted bicontact packages and realizing these
data?
\end{quest}

Note that the monodromy $h$ must preserve the induced line field on $S$. Moreover, because of strong adaptedness of the package, then, at least heuristically, the freedom one should have in building $\xi_+$ from $(S,h)$ and $\ell_+$, should only be a choice of ``how many turns'' does $\xi_+$ realizes along the segments $\{x\}\times [0,1]$ in the mapping torus of $h$. This number of turns cannot be completely free, as we have seen in the padding surgery; adding turns in $\intr(S)\times (0,1)$ require a Dehn surgery on the boundary orbits, so in this context, adding turns change the monodromy by some Dehn twists along the boundary components of $S$.

\subsection{Relation with \cite{BBMY26}} \label{subsec:relationwithbbmy}

In \cite{BBMY26}, B\'eguin--Bonatti--Ma--Yu construct Anosov flows on many fibered hyperbolic 3-manifolds. 
We review a prototypical example of their construction: 
Let $\Sigma_1$ and $\Sigma_2$ be two one-holed tori, equipped with hyperbolic metrics with geodesic boundary components.
For each $i=1,2$, consider the unit tangent bundle $T^1 \Sigma_i$ over $\Sigma_i$. Let $S_i$ be a section to the fibering $T^1 \Sigma_i \to \Sigma_i$. There is a natural identification $S_i \cong \Sigma_i$. On the boundary torus $\partial T^1 \Sigma_i$, we have the homology class $l_i$ of the closed orbit lying over $\partial \Sigma_i$ with the same orientation, the homology class $f_i$ of the Seifert fibers, and the homology class of $\partial S_i$, which is given by $l_i-f_i$.

It is then possible to glue the geodesic flows on $T^1 \Sigma_1$ and $T^1 \Sigma_2$ with a map on the boundary tori that identifies
$l_1 \leftrightarrow l_2$ and $l_1-f_1 \leftrightarrow -(l_2-f_2)$
in order to obtain a totally free Anosov flow $\phi$ on a closed oriented manifold $M$.
The surfaces $S_1$ and $S_2$ glue up to a closed fiber surface $S$ in $M$. The monodromy on $S$ is $\tau_d^{-2}$ where $d$ is the glued common boundary of $S_1$ and $S_2$.

Now, for each (essential) curve $c \subset \Sigma_i$, there is a closed orbit $\gamma$ of $\phi$ in $T^1 \Sigma_i$ lying over $c$. 
Moreover, by picking $S_i$ appropriately, $\gamma$ can be isotoped to the curve $c \subset S_i \subset S$ in $M$, and the framing on $\gamma$ induced by $S$ under such an isotopy is the same as the degeneracy slope. Thus by performing Goodman surgery along $\gamma$, one can modify the monodromy on $S$ by adding any power of a Dehn twist along $c$.

Perhaps more surprisingly, for a curve $c \subset S$ not contained in any $\Sigma_i$, it is still true that there is a closed orbit $\gamma$ of $\phi$ which can be isotoped to $c$. However, it is no longer true that the framing on $\gamma$ induced by $S$ is the same as the degeneracy slope. 
B\'eguin--Bonatti--Ma--Yu showed via an explicit computation that if the geometric intersection number between $c$ and the curve $d$ above is $2k$, then the two framings on $\gamma$ differ by $k$. 
For $k=1,2$, one can still perform \emph{certain} Goodman surgeries along $\gamma$ to modify the monodromy on $S$ by \emph{certain} powers of a Dehn twist along $c$. 

We now attempt to recast this construction in a bicontact framework:
Recall the standard bicontact package for geodesic flows, as described in \Cref{subsec:introgeodflow}, where the Reeb flow $R_-$ is exactly the Seifert fibering.
As shown in \cite{PinSal26}, it is possible to glue together these bicontact packages on $T^1\Sigma_1$ and $T^1\Sigma_2$ by interpolating near the gluing torus. 
In the current situation, the glued Reeb flow will be positively transverse to $S_1$, negatively transverse to $S_2$, and tangent to $d$. 
From this, one can build a contact vector field that is positively transverse to $S$ everywhere, and for which $S_+ = S_1$, $S_- = S_2$, and $\Gamma = d$.

Under this setup, one can attempt to run the strategy as described in \Cref{subsec:intropaddingmotivation}:
Given a curve $c \subset S$, try to apply Legendrian realization principle, in the form of our \Cref{thm:legendrianrealization}, to isotope $S$ so that $c$ is Legendrian with respect to $\xi_-$. 
If this can be accomplished, then one can talk about the Legendrian framing on $c$ induced by $\xi_-$.
A computation, standard in contact topology, would show that this framing differs from that induced by $S$ by half the number of times $c$ intersects $d$.
This would match up with the computation in \cite{BBMY26}.
Next, try to isotope $\xi_+$ along flow lines until $c$ is Legendrian with respect to both $\xi_+$ and $\xi_-$, hence is a closed orbit.
Then one would be able to perform Goodman surgery along $c$. %

One runs into the same three issues pointed out in \Cref{subsec:intropaddingmotivation}:
\begin{enumerate}
    \item One needs enough `height' of $S$ with respect to a contact vector field $V$ to embed the isotopy provided by Legendrian realization principle.
    \item One needs the winding number of $c$ with respect to $\xi_+$ to be zero. Here, as long as $c$ intersects $d$ transversely, $c$ should be transverse to $R_-$ everywhere, so the notion of winding number makes sense.
    \item One needs enough rotation of $\xi_+$ along $R_-$ in order to embed the `pullback annulus'.
\end{enumerate}

We suspect that these issues might not come up in the setting of \cite{BBMY26}:
Using the fact that $\Sigma_i$ are one-holed tori, they choose the sections $S_i$ to $T^1\Sigma_i \to \Sigma_i$ in a specific way. They then make specific choices for the curves $c$. It might be the case that $S$ is already `sufficiently Legendrian' near those curves so that one doesn't need much height to isotope it, avoiding (1). One can also check that the winding numbers of those curves are indeed zero, which addresses (2).
Finally, the fact that $\xi_+$ already completes a full rotation before $S_i$ bumps into itself under $R_-$ in $T^1 \Sigma_i$ might be enough to avoid (3).

However, in general, it would be interesting to define a padding surgery that works in this context. 
We point out that the differences between this context and that in our paper are that
\begin{enumerate}[label=(\roman*)]
    \item in our paper, we have a Birkhoff section, which in particular has boundary, while here we have a closed surface, and
    \item in our paper, the Reeb vector field is transverse to the Birkhoff section, while here the Reeb vector field is tangent to the surface along a curve,
\end{enumerate}
so \Cref{constr:paddingsurgery} does not directly apply.

The generalized padding operation should accomplish the same thing as taking fiberwise covers of $T^1 \Sigma_i$ before gluing, which as in \cite{BBMY26}, allows one to locate more orbits parallel to each curve $c$ on $S$, so that one can change the monodromy by more Dehn twist terms. 

An optimistic goal is to use such an operation, as well as other ideas in this paper and \cite{BBMY26}, to give an affirmative answer to the following question of Potrie.

\begin{quest}[{Potrie, \cite[Question 8]{Pot18}}]
Does every hyperbolic 3-manifold admit a finite cover supporting an Anosov flow?
\end{quest}

\bibliographystyle{alphaurl}

\bibliography{bib.bib}

@article {Iri15,
    AUTHOR = {Irie, Kei},
     TITLE = {Dense existence of periodic {R}eeb orbits and {ECH} spectral
              invariants},
   JOURNAL = {J. Mod. Dyn.},
  FJOURNAL = {Journal of Modern Dynamics},
    VOLUME = {9},
      YEAR = {2015},
     PAGES = {357--363},
      ISSN = {1930-5311,1930-532X},
       DOI = {10.3934/jmd.2015.9.357},
       URL = {https://doi.org/10.3934/jmd.2015.9.357},
}

@article {Ano67,
    AUTHOR = {Anosov, D. V.},
     TITLE = {Geodesic flows on closed {R}iemannian manifolds of negative curvature},
   JOURNAL = {Trudy Mat. Inst. Steklov.},
  FJOURNAL = {Akademiya Nauk SSSR. Trudy Matematicheskogo Instituta imeni V. A. Steklova},
    VOLUME = {90},
      YEAR = {1967},
     PAGES = {209},
      ISSN = {0371-9685},
   MRCLASS = {57.36 (28.00)},
  MRNUMBER = {224110},
MRREVIEWER = {M.\ Eisenberg},
}

@article {ABM24,
    AUTHOR = {Asaoka, Masayuki and Bonatti, Christian and Marty, Th\'eo},
     TITLE = {Oriented {B}irkhoff sections of {A}nosov flows},
   JOURNAL = {J. Topol.},
  FJOURNAL = {Journal of Topology},
    VOLUME = {17},
      YEAR = {2024},
    NUMBER = {4},
     PAGES = {Paper No. e12356, 45},
      ISSN = {1753-8416,1753-8424},
   MRCLASS = {37D20 (37C86)},
  MRNUMBER = {4806779},
       DOI = {10.1112/topo.12356},
       URL = {https://doi.org/10.1112/topo.12356},
}

@article {Bar96,
    AUTHOR = {Barbot, Thierry},
     TITLE = {Flots d'{A}nosov sur les vari\'{e}t\'{e}s graph\'{e}es au sens
              de {W}aldhausen},
   JOURNAL = {Ann. Inst. Fourier (Grenoble)},
  FJOURNAL = {Universit\'{e} de Grenoble. Annales de l'Institut Fourier},
    VOLUME = {46},
      YEAR = {1996},
    NUMBER = {5},
     PAGES = {1451--1517},
      ISSN = {0373-0956,1777-5310},
   MRCLASS = {57R30 (57N10 58F15 58F18)},
  MRNUMBER = {1427133},
MRREVIEWER = {Rafael\ Oswaldo\ Ruggiero},
       DOI = {10.5802/aif.1556},
       URL = {https://doi.org/10.5802/aif.1556},
}

@article {BF13,
    AUTHOR = {Barbot, Thierry and Fenley, S\'{e}rgio R.},
     TITLE = {Pseudo-{A}nosov flows in toroidal manifolds},
   JOURNAL = {Geom. Topol.},
  FJOURNAL = {Geometry \& Topology},
    VOLUME = {17},
      YEAR = {2013},
    NUMBER = {4},
     PAGES = {1877--1954},
      ISSN = {1465-3060},
   MRCLASS = {37D20 (37D50 57R30)},
  MRNUMBER = {3109861},
MRREVIEWER = {Boris Hasselblatt},
       DOI = {10.2140/gt.2013.17.1877},
       URL = {https://doi.org/10.2140/gt.2013.17.1877},
}

@article {BF15,
    AUTHOR = {Barbot, Thierry and Fenley, S\'{e}rgio R.},
     TITLE = {Classification and rigidity of totally periodic
              pseudo-{A}nosov flows in graph manifolds},
   JOURNAL = {Ergodic Theory Dynam. Systems},
  FJOURNAL = {Ergodic Theory and Dynamical Systems},
    VOLUME = {35},
      YEAR = {2015},
    NUMBER = {6},
     PAGES = {1681--1722},
      ISSN = {0143-3857},
   MRCLASS = {37D20 (54H20 57N12)},
  MRNUMBER = {3377280},
MRREVIEWER = {Krzysztof Ciesielski},
       DOI = {10.1017/etds.2014.9},
       URL = {https://doi.org/10.1017/etds.2014.9},
}

@article {BF21,
    AUTHOR = {Barbot, Thierry and Fenley, S\'{e}rgio R.},
     TITLE = {Free {S}eifert pieces of pseudo-{A}nosov flows},
   JOURNAL = {Geom. Topol.},
  FJOURNAL = {Geometry \& Topology},
    VOLUME = {25},
      YEAR = {2021},
    NUMBER = {3},
     PAGES = {1331--1440},
      ISSN = {1465-3060},
   MRCLASS = {37D05 (34C25 34C45 34D23 37C83 37D20 57M60 57R30)},
  MRNUMBER = {4268164},
MRREVIEWER = {Thomas Barthelm\'{e}},
       DOI = {10.2140/gt.2021.25.1331},
       URL = {https://doi.org/10.2140/gt.2021.25.1331},
}

@article {BBowM24,
    AUTHOR = {Barthelm\'e, Thomas and Mann, Kathryn},
     TITLE = {Orbit equivalences of {$\Bbb{R}$}-covered {A}nosov flows and
              hyperbolic-like actions on the line},
      NOTE = {Appendix written jointly with Jonathan Bowden},
   JOURNAL = {Geom. Topol.},
  FJOURNAL = {Geometry \& Topology},
    VOLUME = {28},
      YEAR = {2024},
    NUMBER = {2},
     PAGES = {867--899},
      ISSN = {1465-3060,1364-0380},
   MRCLASS = {37D20 (37A20 57M60)},
  MRNUMBER = {4718129},
       DOI = {10.2140/gt.2024.28.867},
       URL = {https://doi.org/10.2140/gt.2024.28.867},
}

@misc{BarMan25,
      title={Pseudo-{A}nosov Flows: A Plane Approach}, 
      author={Thomas Barthelmé and Kathryn Mann},
      year={2026},
      eprint={2509.15375},
      archivePrefix={arXiv},
      primaryClass={math.DS},
      url={https://arxiv.org/abs/2509.15375}, 
}

@article {Bar25,
    AUTHOR = {Barthelm\'{e}, Thomas},
     TITLE = {A sm\"{o}rg\aa sbord of (bi)contact structures, {R}eeb flows
              and pseudo-{A}nosov flows},
   JOURNAL = {J. Fixed Point Theory Appl.},
  FJOURNAL = {Journal of Fixed Point Theory and Applications},
    VOLUME = {27},
      YEAR = {2025},
    NUMBER = {2},
     PAGES = {Paper No. 23, 24},
      ISSN = {1661-7738,1661-7746},
   MRCLASS = {53D10 (37D20 37J55 57K33)},
  MRNUMBER = {4875002},
MRREVIEWER = {Boris\ Hasselblatt},
       DOI = {10.1007/s11784-025-01176-5},
       URL = {https://doi.org/10.1007/s11784-025-01176-5},
}

@misc{BP26,
      title={Uniqueness of gluings and virtual finiteness of pseudo-{A}nosov flows on graph manifolds}, 
      author={Thomas Barthelmé and Neige Paulet},
      year={2026},
      eprint={2608.13526},
      archivePrefix={arXiv},
      primaryClass={math.DS},
      url={https://arxiv.org/abs/2608.13526}, 
}

@article {BBY17,
    AUTHOR = {B\'{e}guin, Fran\c{c}ois and Bonatti, Christian and Yu, Bin},
     TITLE = {Building {A}nosov flows on 3-manifolds},
   JOURNAL = {Geom. Topol.},
  FJOURNAL = {Geometry \& Topology},
    VOLUME = {21},
      YEAR = {2017},
    NUMBER = {3},
     PAGES = {1837--1930},
      ISSN = {1465-3060},
   MRCLASS = {37D20 (57M50)},
  MRNUMBER = {3650083},
MRREVIEWER = {Rafael Oswaldo Ruggiero},
       DOI = {10.2140/gt.2017.21.1837},
       URL = {https://doi.org/10.2140/gt.2017.21.1837},
}

@misc{BY24,
      title={Existence of arbitrary large numbers of non-$\mathbb{R}$-covered {A}nosov flows on hyperbolic $3$-manifolds}, 
      author={Francois Béguin and Bin Yu},
      year={2024},
      eprint={2402.06551},
      archivePrefix={arXiv},
      primaryClass={math.DS},
      url={https://arxiv.org/abs/2402.06551}, 
}

@misc{BBMY26,
      title={Construction of Anosov flows on fibered hyperbolic 3-manifolds}, 
      author={François Béguin and Christian Bonatti and Biao Ma and Bin Yu},
      year={2026},
      eprint={2603.06105},
      archivePrefix={arXiv},
      primaryClass={math.DS},
      url={https://arxiv.org/abs/2603.06105}, 
}

@article {Bru93,
    AUTHOR = {Brunella, Marco},
     TITLE = {Separating the basic sets of a nontransitive {A}nosov flow},
   JOURNAL = {Bull. London Math. Soc.},
  FJOURNAL = {The Bulletin of the London Mathematical Society},
    VOLUME = {25},
      YEAR = {1993},
    NUMBER = {5},
     PAGES = {487--490},
      ISSN = {0024-6093},
   MRCLASS = {58F15},
  MRNUMBER = {1233413},
MRREVIEWER = {Lawrence Conlon},
       DOI = {10.1112/blms/25.5.487},
       URL = {https://doi.org/10.1112/blms/25.5.487},
}

@article {BI23,
    AUTHOR = {Bonatti, Christian and Iakovoglou, Ioannis},
     TITLE = {Anosov flows on 3-manifolds: the surgeries and the foliations},
   JOURNAL = {Ergodic Theory Dynam. Systems},
  FJOURNAL = {Ergodic Theory and Dynamical Systems},
    VOLUME = {43},
      YEAR = {2023},
    NUMBER = {4},
     PAGES = {1129--1188},
      ISSN = {0143-3857,1469-4417},
   MRCLASS = {37D20 (37D40 57M10 57R30)},
  MRNUMBER = {4555824},
       DOI = {10.1017/etds.2021.170},
       URL = {https://doi.org/10.1017/etds.2021.170},
}

@article {BM22,
    AUTHOR = {Bowden, Jonathan and Mann, Kathryn},
     TITLE = {{$C^0$} stability of boundary actions and inequivalent
              {A}nosov flows},
   JOURNAL = {Ann. Sci. \'{E}c. Norm. Sup\'{e}r. (4)},
  FJOURNAL = {Annales Scientifiques de l'\'{E}cole Normale Sup\'{e}rieure.
              Quatri\`eme S\'{e}rie},
    VOLUME = {55},
      YEAR = {2022},
    NUMBER = {4},
     PAGES = {1003--1046},
      ISSN = {0012-9593,1873-2151},
   MRCLASS = {53C22 (37C15 37C85 37D20)},
  MRNUMBER = {4468857},
MRREVIEWER = {Boris\ Hasselblatt},
       DOI = {10.24033/asens.2512},
       URL = {https://doi.org/10.24033/asens.2512},
}

@misc{BC26,
      title={On the Finiteness of {A}nosov flows on $3$-manifolds via Contact Geometry}, 
      author={Jonathan Bowden and Vincent Colin},
      year={2026},
      eprint={2609.03700},
      archivePrefix={arXiv},
      primaryClass={math.DS},
      url={https://arxiv.org/abs/2609.03700}, 
}

@article {CS23,
    AUTHOR = {Calderon, Aaron and Salter, Nick},
     TITLE = {Framed mapping class groups and the monodromy of strata of
              abelian differentials},
   JOURNAL = {J. Eur. Math. Soc. (JEMS)},
  FJOURNAL = {Journal of the European Mathematical Society (JEMS)},
    VOLUME = {25},
      YEAR = {2023},
    NUMBER = {12},
     PAGES = {4719--4790},
      ISSN = {1435-9855,1435-9863},
   MRCLASS = {57K20 (30F30)},
  MRNUMBER = {4662301},
       DOI = {10.4171/jems/1290},
       URL = {https://doi.org/10.4171/jems/1290},
}

@article {CP25,
    AUTHOR = {Clay, Adam and Pinsky, Tali},
     TITLE = {Graph manifolds that admit arbitrarily many {A}nosov flows},
   JOURNAL = {Math. Ann.},
  FJOURNAL = {Mathematische Annalen},
    VOLUME = {391},
      YEAR = {2025},
    NUMBER = {4},
     PAGES = {5375--5390},
      ISSN = {0025-5831,1432-1807},
   MRCLASS = {37D40 (37D20 57M50)},
  MRNUMBER = {4884550},
       DOI = {10.1007/s00208-024-03048-8},
       URL = {https://doi.org/10.1007/s00208-024-03048-8},
}

@article {CDHR24,
    AUTHOR = {Colin, Vincent and Dehornoy, Pierre and Hryniewicz, Umberto
              and Rechtman, Ana},
     TITLE = {Generic properties of 3-dimensional {R}eeb flows: {B}irkhoff
              sections and entropy},
   JOURNAL = {Comment. Math. Helv.},
  FJOURNAL = {Commentarii Mathematici Helvetici. A Journal of the Swiss
              Mathematical Society},
    VOLUME = {99},
      YEAR = {2024},
    NUMBER = {3},
     PAGES = {557--611},
      ISSN = {0010-2571,1420-8946},
   MRCLASS = {53E50 (37C27 37C35 57K33)},
  MRNUMBER = {4771036},
MRREVIEWER = {Boris\ Hasselblatt},
       DOI = {10.4171/cmh/573},
       URL = {https://doi.org/10.4171/cmh/573},
}

@misc{Deh26,
      title={Genus one {B}irkhoff sections for geodesic flows on orbifolds}, 
      author={Pierre Dehornoy},
      year={2026},
      eprint={2603.23000},
      archivePrefix={arXiv},
      primaryClass={math.DS},
      url={https://arxiv.org/abs/2603.23000}, 
}

@misc{Deh,
      title={Birkhoff sections for {A}nosov flows: from genus two to genus one}, 
      author={Pierre Dehornoy},
      note={In preparation},
}

@book {ET98,
    AUTHOR = {Eliashberg, Yakov M. and Thurston, William P.},
     TITLE = {Confoliations},
    SERIES = {University Lecture Series},
    VOLUME = {13},
 PUBLISHER = {American Mathematical Society, Providence, RI},
      YEAR = {1998},
     PAGES = {x+66},
      ISBN = {0-8218-0776-5},
   MRCLASS = {53C15 (57N10 57R30 58F05)},
  MRNUMBER = {1483314},
MRREVIEWER = {Hansj\"{o}rg Geiges},
       DOI = {10.1090/ulect/013},
       URL = {https://doi.org/10.1090/ulect/013},
}

@misc{ET,
	Author = {John B. Etnyre and Bülent Tosun},
	Title = {Low-{D}imensional {C}ontact {G}eometry},
	Note = {In preparation},
URL = {https://etnyre.math.gatech.edu/preprints/ContactBookV0.1.pdf}}

@book {FM12,
    AUTHOR = {Farb, Benson and Margalit, Dan},
     TITLE = {A primer on mapping class groups},
    SERIES = {Princeton Mathematical Series},
    VOLUME = {49},
 PUBLISHER = {Princeton University Press, Princeton, NJ},
      YEAR = {2012},
     PAGES = {xiv+472},
      ISBN = {978-0-691-14794-9},
   MRCLASS = {57M50 (20F36 20F65 57M07 57N05)},
  MRNUMBER = {2850125},
MRREVIEWER = {Stephen P. Humphries},
}

@article {Fen05,
    AUTHOR = {Fenley, S\'{e}rgio R.},
     TITLE = {Regulating flows, topology of foliations and rigidity},
   JOURNAL = {Trans. Amer. Math. Soc.},
  FJOURNAL = {Transactions of the American Mathematical Society},
    VOLUME = {357},
      YEAR = {2005},
    NUMBER = {12},
     PAGES = {4957--5000},
      ISSN = {0002-9947},
   MRCLASS = {37D20 (37C85 53C12 53C23 57R30)},
  MRNUMBER = {2165394},
MRREVIEWER = {Thierry Barbot},
       DOI = {10.1090/S0002-9947-05-03644-5},
       URL = {https://doi.org/10.1090/S0002-9947-05-03644-5},
}

@article {Fri83,
	AUTHOR = {Fried, David},
	TITLE = {Transitive {A}nosov flows and pseudo-{A}nosov maps},
	JOURNAL = {Topology},
	FJOURNAL = {Topology. An International Journal of Mathematics},
	VOLUME = {22},
	YEAR = {1983},
	NUMBER = {3},
	PAGES = {299--303},
	ISSN = {0040-9383},
	MRCLASS = {58F15 (57R30)},
	MRNUMBER = {710103},
	MRREVIEWER = {Michael Handel},
	DOI = {10.1016/0040-9383(83)90015-0},
	URL = {https://doi.org/10.1016/0040-9383(83)90015-0},
}

@book {FH19,
    AUTHOR = {Fisher, Todd and Hasselblatt, Boris},
     TITLE = {Hyperbolic flows},
    SERIES = {Zurich Lectures in Advanced Mathematics},
 PUBLISHER = {EMS Publishing House, Berlin},
      YEAR = {[2019] \copyright 2019},
     PAGES = {xiv+723},
      ISBN = {978-3-03719-200-9},
   MRCLASS = {37-02 (37A30 37A35 37D20 37D40)},
  MRNUMBER = {3972204},
MRREVIEWER = {Miguel Paternain},
       DOI = {10.4171/200},
       URL = {https://doi.org/10.4171/200},
}

@misc{GL26,
      title={Finiteness of {A}nosov and pseudo-{A}nosov flows on hyperbolic 3-manifolds}, 
      author={David Gabai and Lily Qiao Li},
      year={2026},
      eprint={2610.24190},
      archivePrefix={arXiv},
      primaryClass={math.GT},
      url={https://arxiv.org/abs/2610.24190}, 
}

@article {Ghy84,
    AUTHOR = {Ghys, \'{E}tienne},
     TITLE = {Flots d'{A}nosov sur les {$3$}-vari\'{e}t\'{e}s fibr\'{e}es en
              cercles},
   JOURNAL = {Ergodic Theory Dynam. Systems},
  FJOURNAL = {Ergodic Theory and Dynamical Systems},
    VOLUME = {4},
      YEAR = {1984},
    NUMBER = {1},
     PAGES = {67--80},
      ISSN = {0143-3857,1469-4417},
   MRCLASS = {58F15 (58F17)},
  MRNUMBER = {758894},
MRREVIEWER = {Alberto\ Verjovsky},
       DOI = {10.1017/S0143385700002273},
       URL = {https://doi.org/10.1017/S0143385700002273},
}

@incollection {Goo83,
	AUTHOR = {Goodman, Sue},
	TITLE = {Dehn surgery on {A}nosov flows},
	BOOKTITLE = {Geometric dynamics ({R}io de {J}aneiro, 1981)},
	SERIES = {Lecture Notes in Math.},
	VOLUME = {1007},
	PAGES = {300--307},
	PUBLISHER = {Springer, Berlin},
	YEAR = {1983},
	MRCLASS = {58F15 (57R65)},
	MRNUMBER = {1691596},
	DOI = {10.1007/BFb0061421},
	URL = {https://doi.org/10.1007/BFb0061421},
}

@misc{Hon04,
      title={Notes for Math 599: Contact geometry}, 
      author={Ko Honda},
      year={2004},
      url={https://www.math.ucla.edu/~honda/math599/notes.pdf}, 
}

@article {Hoz25,
    AUTHOR = {Hozoori, Surena},
     TITLE = {Strongly adapted contact geometry of {A}nosov 3-flows},
   JOURNAL = {J. Fixed Point Theory Appl.},
  FJOURNAL = {Journal of Fixed Point Theory and Applications},
    VOLUME = {27},
      YEAR = {2025},
    NUMBER = {2},
     PAGES = {Paper No. 37, 37},
      ISSN = {1661-7738,1661-7746},
   MRCLASS = {53E50 (37D20 37J55)},
  MRNUMBER = {4886075},
MRREVIEWER = {Boris\ Hasselblatt},
       DOI = {10.1007/s11784-025-01180-9},
       URL = {https://doi.org/10.1007/s11784-025-01180-9},
}

@misc{Hoz26,
      title={Bi-contact plugs and dynamical hyperbolicity}, 
      author={Surena Hozoori},
      year={2026},
      eprint={2610.11243},
      archivePrefix={arXiv},
      primaryClass={math.SG},
      url={https://arxiv.org/abs/2610.11243}, 
}

@article {KKvK22,
    AUTHOR = {Kim, Juno and Kim, Yonghwan and van Koert, Otto},
     TITLE = {Reeb flows without simple global surfaces of section},
   JOURNAL = {Involve},
  FJOURNAL = {Involve. A Journal of Mathematics},
    VOLUME = {15},
      YEAR = {2022},
    NUMBER = {5},
     PAGES = {813--842},
      ISSN = {1944-4176,1944-4184},
   MRCLASS = {37J55 (53D10)},
  MRNUMBER = {4555146},
       DOI = {10.2140/involve.2022.15.813},
       URL = {https://doi.org/10.2140/involve.2022.15.813},
}

@incollection {Kir97,
     TITLE = {Problems in low-dimensional topology},
 BOOKTITLE = {Geometric topology ({A}thens, {GA}, 1993)},
    SERIES = {AMS/IP Stud. Adv. Math.},
    VOLUME = {2},
     PAGES = {35--473},
    EDITOR = {Kirby, Rob},
 PUBLISHER = {Amer. Math. Soc., Providence, RI},
      YEAR = {1997},
      ISBN = {0-8218-0653-X},
   MRCLASS = {57-02},
  MRNUMBER = {1470751},
       DOI = {10.1090/amsip/002.2/02},
       URL = {https://doi.org/10.1090/amsip/002.2/02},
}

@book {Kru97,
    AUTHOR = {Krupkov\'{a}, Olga},
     TITLE = {The geometry of ordinary variational equations},
    SERIES = {Lecture Notes in Mathematics},
    VOLUME = {1678},
 PUBLISHER = {Springer-Verlag, Berlin},
      YEAR = {1997},
     PAGES = {x+251},
      ISBN = {3-540-63832-6},
   MRCLASS = {58F05 (34A26 49J10 58E30 70Hxx)},
  MRNUMBER = {1484970},
MRREVIEWER = {Frans\ Cantrijn},
       DOI = {10.1007/BFb0093438},
       URL = {https://doi.org/10.1007/BFb0093438},
}

@article{Mar25,
   title={Skewed {A}nosov flows in dimension 3 are Reeb-like},
   ISSN={1435-9863},
   url={http://dx.doi.org/10.4171/jems/1619},
   DOI={10.4171/jems/1619},
   journal={Journal of the European Mathematical Society},
   publisher={European Mathematical Society - EMS - Publishing House GmbH},
   author={Marty, Théo Anaël},
   year={2025},
   month=Mar }

@misc{ParSeg26,
      title={Drilling veering triangulations with applications to pseudo-Anosov flows}, 
      author={Anna Parlak and Henry Segerman},
      year={2026},
      eprint={2610.12302},
      archivePrefix={arXiv},
      primaryClass={math.GT},
      url={https://arxiv.org/abs/2610.12302}, 
}

@article {Pau25,
    AUTHOR = {Paulet, Neige},
     TITLE = {Anosov flows in dimension 3 from gluing building blocks with
              quasi-transverse boundary},
   JOURNAL = {J. Mod. Dyn.},
  FJOURNAL = {Journal of Modern Dynamics},
    VOLUME = {21},
      YEAR = {2025},
     PAGES = {21--240},
      ISSN = {1930-5311,1930-532X},
   MRCLASS = {37D20 (37C10 37D05 57K30)},
  MRNUMBER = {4879873},
MRREVIEWER = {Thilo\ Kuessner},
       DOI = {10.3934/jmd.2025002},
       URL = {https://doi.org/10.3934/jmd.2025002},
}

@article {PT72,
    AUTHOR = {Plante, J. F. and Thurston, W. P.},
     TITLE = {Anosov flows and the fundamental group},
   JOURNAL = {Topology},
  FJOURNAL = {Topology. An International Journal of Mathematics},
    VOLUME = {11},
      YEAR = {1972},
     PAGES = {147--150},
      ISSN = {0040-9383},
   MRCLASS = {58F10},
  MRNUMBER = {295389},
MRREVIEWER = {John\ M.\ Franks},
       DOI = {10.1016/0040-9383(72)90002-X},
       URL = {https://doi.org/10.1016/0040-9383(72)90002-X},
}

@misc{PinSal26,
      title={New {Anosov} flows via bicontact structures}, 
      author={Tali Pinsky and Federico Salmoiraghi},
      year={2026},
      eprint={2603.22250},
      archivePrefix={arXiv},
      primaryClass={math.GT},
      url={https://arxiv.org/abs/2603.22250},
}

@inproceedings {Pot18,
    AUTHOR = {Potrie, Rafael},
     TITLE = {Robust dynamics, invariant structures and topological
              classification},
 BOOKTITLE = {Proceedings of the {I}nternational {C}ongress of
              {M}athematicians---{R}io de {J}aneiro 2018. {V}ol. {III}.
              {I}nvited lectures},
     PAGES = {2063--2085},
 PUBLISHER = {World Sci. Publ., Hackensack, NJ},
      YEAR = {2018},
      ISBN = {978-981-3272-92-7; 978-981-3272-87-3},
   MRCLASS = {37D30 (37C20 57R30)},
  MRNUMBER = {3966842},
}

@article {Sal25,
    AUTHOR = {Salmoiraghi, Federico},
     TITLE = {Surgery on {A}nosov flows using bi-contact geometry},
   JOURNAL = {Ergodic Theory Dynam. Systems},
  FJOURNAL = {Ergodic Theory and Dynamical Systems},
    VOLUME = {45},
      YEAR = {2025},
    NUMBER = {12},
     PAGES = {3832--3864},
      ISSN = {0143-3857,1469-4417},
   MRCLASS = {37D20 (53D25 57K33)},
  MRNUMBER = {4986560},
       DOI = {10.1017/etds.2025.10196},
       URL = {https://doi.org/10.1017/etds.2025.10196},
}

@article {Tsa24,
    AUTHOR = {Tsang, Chi Cheuk},
     TITLE = {Constructing {B}irkhoff sections for pseudo-{A}nosov flows
              with controlled complexity},
   JOURNAL = {Ergodic Theory Dynam. Systems},
  FJOURNAL = {Ergodic Theory and Dynamical Systems},
    VOLUME = {44},
      YEAR = {2024},
    NUMBER = {8},
     PAGES = {2308--2360},
      ISSN = {0143-3857},
   MRCLASS = {37D20 (57K35)},
  MRNUMBER = {4798798},
MRREVIEWER = {Boris Hasselblatt},
       DOI = {10.1017/etds.2023.105},
       URL = {https://doi.org/10.1017/etds.2023.105},
}

@misc{Tsa24a,
      title={Horizontal {G}oodman surgery and almost equivalence of pseudo-{A}nosov flows}, 
      author={Chi Cheuk Tsang},
      year={2024},
      eprint={2401.01847},
      archivePrefix={arXiv},
      primaryClass={math.DS},
      url={https://arxiv.org/abs/2401.01847}, 
}

@misc{Tsa24b,
      title={Examples of {A}nosov flows with genus one {B}irkhoff sections}, 
      author={Chi Cheuk Tsang},
      year={2024},
      eprint={2402.00229},
      archivePrefix={arXiv},
      primaryClass={math.DS},
      url={https://arxiv.org/abs/2402.00229}, 
}

@misc{Tsa26,
      title={From pseudo-Anosov flows on graph manifolds to totally periodic flows}, 
      author={Chi Cheuk Tsang},
      year={2026},
      eprint={2603.27362},
      archivePrefix={arXiv},
      primaryClass={math.DS},
      url={https://arxiv.org/abs/2603.27362}, 
}

@article {Yu23,
    AUTHOR = {Yu, Bin},
     TITLE = {Anosov flows on {D}ehn surgeries on the figure-eight knot},
   JOURNAL = {Duke Math. J.},
  FJOURNAL = {Duke Mathematical Journal},
    VOLUME = {172},
      YEAR = {2023},
    NUMBER = {11},
     PAGES = {2195--2240},
      ISSN = {0012-7094,1547-7398},
   MRCLASS = {37D40 (37D20 57K32)},
  MRNUMBER = {4627250},
MRREVIEWER = {Thilo\ Kuessner},
       DOI = {10.1215/00127094-2022-0079},
       URL = {https://doi.org/10.1215/00127094-2022-0079},
}

\end{document}